\voffset=-.5in
\magnification=1200
\vsize=7.5in\hsize=7.00 true in\hoffset=-.25in
\tolerance = 100000
$\ $\vskip .3in

\centerline{``On the Strong Equality between Supercompactness and
Strong Compactness''}
\vskip .25in
\centerline{by}
\vskip .25in
\centerline{Arthur W. Apter*}
\centerline{Department of Mathematics}
\centerline{Baruch College of CUNY}
\centerline{New York, New York 10010}
\vskip .125in
\centerline{and}
\vskip .125in
\centerline{Saharon Shelah**}
\centerline{Department of Mathematics}
\centerline{The Hebrew University}
\centerline{Jerusalem, Israel}
\vskip .125in\centerline{and}\vskip .125in
\centerline{Department of Mathematics}
\centerline{Rutgers University}
\centerline{New Brunswick, New Jersey 08904}
\vskip .250in
\centerline{February 19, 1995}
\vskip .25in
\noindent Abstract: We show that supercompactness and strong
compactness can be equivalent even as properties of pairs of regular
cardinals.  
Specifically, we show that if $V \models$ ZFC + GCH is a given
model (which in interesting cases contains instances of
supercompactness), then there is some cardinal and cofinality
preserving generic extension $V[{G}] \models$ ZFC + GCH in
which, (a) (preservation) for $\kappa \le \lambda$ regular, if
$V      \models ``\kappa$ is $\lambda$ supercompact'', then
$V[G] \models ``\kappa$ is $\lambda$ supercompact''
and so
that, (b) (equivalence) for $\kappa \le \lambda$ regular,
$V[{G}] \models ``\kappa$ is $\lambda$ strongly compact'' iff
$V[{G}] \models ``\kappa$ is $\lambda$ supercompact'', except
possibly if $\kappa$ is a measurable limit of
cardinals which are $\lambda $ supercompact.
\hfil
\vskip .25in
\noindent *The research of the first author was partially supported
by PSC-CUNY Grant 662341 and a salary grant from Tel Aviv
University.  In addition, the first author wishes to thank the
Mathematics Departments of Hebrew University and Tel Aviv
University for the hospitality shown him during his
sabbatical in Israel.
\hfil\vskip .125in\noindent
**Publication 495.  The second author wishes to thank the
Basic Research Fund of the Israeli Academy of Sciences for
partially supporting this research.
\hfil\break

% VANILLA.STY
% COPYRIGHT (C) 1985, 1986 BY MICHAEL SPIVAK
% version date 1/1/86
\catcode`\@=11
\font\tensmc=cmcsc10      %change to CM fonts 3-31-87
%\font\tensmc=amcsc10
\def\smc{\tensmc}

\def\hcorrection#1{\advance\hoffset by #1 }
\def\vcorrection#1{\advance\voffset by #1 }
\def\wlog#1{}
\newif\iftitle@
\outer\def\title{\title@true\vglue 24\p@ plus 12\p@ minus 12\p@
   \bgroup\let\\=\cr\tabskip\centering
   \halign to \hsize\bgroup\tenbf\hfill\ignorespaces##\unskip\hfill\cr}
\def\endtitle{\cr\egroup\egroup\vglue 18\p@ plus 12\p@ minus 6\p@}
\outer\def\author{\iftitle@\vglue -18\p@ plus -12\p@ minus -6\p@\fi\vglue
    12\p@ plus 6\p@ minus 3\p@\bgroup\let\\=\cr\tabskip\centering
    \halign to \hsize\bgroup\smc\hfill\ignorespaces##\unskip\hfill\cr}
\def\endauthor{\cr\egroup\egroup\vglue 18\p@ plus 12\p@ minus 6\p@}
\outer\def\heading{\bigbreak\bgroup\let\\=\cr\tabskip\centering
    \halign to \hsize\bgroup\smc\hfill\ignorespaces##\unskip\hfill\cr}
\def\endheading{\cr\egroup\egroup\nobreak\medskip}

\outer\def\proclaim#1{\medbreak\noindent\smc\ignorespaces
    #1\unskip.\enspace\sl\ignorespaces}
\outer\def\endproclaim{\par\ifdim\lastskip<\medskipamount\removelastskip
  \penalty 55 \fi\medskip\rm}
\outer\def\demo#1{\par\ifdim\lastskip<\smallskipamount\removelastskip
    \smallskip\fi\noindent{\smc\ignorespaces#1\unskip:\enspace}\rm
      \ignorespaces}

\newcount\footmarkcount@
\footmarkcount@=1
\def\makefootnote@#1#2{\insert\footins{\interlinepenalty=100
  \splittopskip=\ht\strutbox \splitmaxdepth=\dp\strutbox
  \floatingpenalty=\@MM
  \leftskip=\z@\rightskip=\z@\spaceskip=\z@\xspaceskip=\z@
  \noindent{#1}\footstrut\rm\ignorespaces #2\strut}}
\def\footnote{\let\@sf=\empty\ifhmode\edef\@sf{\spacefactor
   =\the\spacefactor}\/\fi\futurelet\next\footnote@}
\def\footnote@{\ifx"\next\let\next\footnote@@\else
    \let\next\footnote@@@\fi\next}
\def\footnote@@"#1"#2{#1\@sf\relax\makefootnote@{#1}{#2}}
\def\footnote@@@#1{$^{\number\footmarkcount@}$\makefootnote@
   {$^{\number\footmarkcount@}$}{#1}\global\advance\footmarkcount@ by 1 }

\hyphenation{man-u-script man-u-scripts ap-pen-dix ap-pen-di-ces}
\hyphenation{data-base data-bases}
\ifx\amstexloaded@\relax\catcode`\@=13
  \endinput\else\let\amstexloaded@=\relax\fi
\newlinechar=`\^^J
\def\eat@#1{}
\def\Space@.{\futurelet\Space@\relax}
\Space@. %
\newhelp\athelp@
{Only certain combinations beginning with @ make sense to me.^^J
Perhaps you wanted \string\@\space for a printed @?^^J
I've ignored the character or group after @.}
\def\futureletnextat@{\futurelet\next\at@}
{\catcode`\@=\active
\lccode`\Z=`\@ \lowercase
{\gdef@{\expandafter\csname futureletnextatZ\endcsname}
\expandafter\gdef\csname atZ\endcsname
   {\ifcat\noexpand\next a\def\next{\csname atZZ\endcsname}\else
   \ifcat\noexpand\next0\def\next{\csname atZZ\endcsname}\else
    \def\next{\csname atZZZ\endcsname}\fi\fi\next}
\expandafter\gdef\csname atZZ\endcsname#1{\expandafter
   \ifx\csname #1Zat\endcsname\relax\def\next
     {\errhelp\expandafter=\csname athelpZ\endcsname
      \errmessage{Invalid use of \string@}}\else
       \def\next{\csname #1Zat\endcsname}\fi\next}
\expandafter\gdef\csname atZZZ\endcsname#1{\errhelp
    \expandafter=\csname athelpZ\endcsname
      \errmessage{Invalid use of \string@}}}}
\def\atdef@#1{\expandafter\def\csname #1@at\endcsname}
\newhelp\defahelp@{If you typed \string\define\space cs instead of
\string\define\string\cs\space^^J
I've substituted an inaccessible control sequence so that your^^J
definition will be completed without mixing me up too badly.^^J
If you typed \string\define{\string\cs} the inaccessible control sequence^^J
was defined to be \string\cs, and the rest of your^^J
definition appears as input.}
\newhelp\defbhelp@{I've ignored your definition, because it might^^J
conflict with other uses that are important to me.}
\def\define{\futurelet\next\define@}
\def\define@{\ifcat\noexpand\next\relax
  \def\next{\define@@}%
  \else\errhelp=\defahelp@
  \errmessage{\string\define\space must be followed by a control
     sequence}\def\next{\def\garbage@}\fi\next}
\def\undefined@{}
\def\preloaded@{}
\def\define@@#1{\ifx#1\relax\errhelp=\defbhelp@
   \errmessage{\string#1\space is already defined}\def\next{\def\garbage@}%
   \else\expandafter\ifx\csname\expandafter\eat@\string
	 #1@\endcsname\undefined@\errhelp=\defbhelp@
   \errmessage{\string#1\space can't be defined}\def\next{\def\garbage@}%
   \else\expandafter\ifx\csname\expandafter\eat@\string#1\endcsname\relax
     \def\next{\def#1}\else\errhelp=\defbhelp@
     \errmessage{\string#1\space is already defined}\def\next{\def\garbage@}%
      \fi\fi\fi\next}
\def\famzero{\fam\z@}

\def\lim{\mathop{\famzero lim}}

\def\max{\mathop{\famzero max}}
\def\min{\mathop{\famzero min}}

\def\sup{\mathop{\famzero sup}}

\def\textfont@#1#2{\def#1{\relax\ifmmode
    \errmessage{Use \string#1\space only in text}\else#2\fi}}
\textfont@\rm\tenrm
\textfont@\it\tenit
\textfont@\sl\tensl
\textfont@\bf\tenbf
\textfont@\smc\tensmc
\let\ic@=\/
\def\/{\unskip\ic@}
\def\textfonti{\the\textfont1 }
\def\t#1#2{{\edef\next{\the\font}\textfonti\accent"7F \next#1#2}}
\let\B=\=
\let\D=\.
\def~{\unskip\nobreak\ \ignorespaces}
{\catcode`\@=\active
\gdef\@{\char'100 }}
\atdef@-{\leavevmode\futurelet\next\athyph@}
\def\athyph@{\ifx\next-\let\next=\athyph@@
  \else\let\next=\athyph@@@\fi\next}
\def\athyph@@@{\hbox{-}}
\def\athyph@@#1{\futurelet\next\athyph@@@@}
\def\athyph@@@@{\if\next-\def\next##1{\hbox{---}}\else
    \def\next{\hbox{--}}\fi\next}
\def\.{.\spacefactor=\@m}
\atdef@.{\null.}
\atdef@,{\null,}
\atdef@;{\null;}
\atdef@:{\null:}
\atdef@?{\null?}
\atdef@!{\null!}
\def\srdr@{\thinspace}
\def\drsr@{\kern.02778em}
\def\sldl@{\kern.02778em}
\def\dlsl@{\thinspace}
\atdef@"{\unskip\futurelet\next\atqq@}
\def\atqq@{\ifx\next\Space@\def\next. {\atqq@@}\else
	 \def\next.{\atqq@@}\fi\next.}
\def\atqq@@{\futurelet\next\atqq@@@}
\def\atqq@@@{\ifx\next`\def\next`{\atqql@}\else\def\next'{\atqqr@}\fi\next}
\def\atqql@{\futurelet\next\atqql@@}
\def\atqql@@{\ifx\next`\def\next`{\sldl@``}\else\def\next{\dlsl@`}\fi\next}
\def\atqqr@{\futurelet\next\atqqr@@}
\def\atqqr@@{\ifx\next'\def\next'{\srdr@''}\else\def\next{\drsr@'}\fi\next}
\def\flushpar{\par\noindent}
\def\textfontii{\the\textfont2 }
\def\{{\relax\ifmmode\lbrace\else
    {\textfontii f}\spacefactor=\@m\fi}
\def\}{\relax\ifmmode\rbrace\else
    \let\@sf=\empty\ifhmode\edef\@sf{\spacefactor=\the\spacefactor}\fi
      {\textfontii g}\@sf\relax\fi}
\def\nonhmodeerr@#1{\errmessage
     {\string#1\space allowed only within text}}
\def\linebreak{\relax\ifhmode\unskip\break\else
    \nonhmodeerr@\linebreak\fi}
\def\allowlinebreak{\relax
   \ifhmode\allowbreak\else\nonhmodeerr@\allowlinebreak\fi}
\newskip\saveskip@
\def\nolinebreak{\relax\ifhmode\saveskip@=\lastskip\unskip
  \nobreak\ifdim\saveskip@>\z@\hskip\saveskip@\fi
   \else\nonhmodeerr@\nolinebreak\fi}
\def\newline{\relax\ifhmode\null\hfil\break
    \else\nonhmodeerr@\newline\fi}
\def\nonmathaerr@#1{\errmessage
     {\string#1\space is not allowed in display math mode}}
\def\nonmathberr@#1{\errmessage{\string#1\space is allowed only in math mode}}
\def\mathbreak{\relax\ifmmode\ifinner\break\else
   \nonmathaerr@\mathbreak\fi\else\nonmathberr@\mathbreak\fi}
\def\nomathbreak{\relax\ifmmode\ifinner\nobreak\else
    \nonmathaerr@\nomathbreak\fi\else\nonmathberr@\nomathbreak\fi}
\def\allowmathbreak{\relax\ifmmode\ifinner\allowbreak\else
     \nonmathaerr@\allowmathbreak\fi\else\nonmathberr@\allowmathbreak\fi}
\def\pagebreak{\relax\ifmmode
   \ifinner\errmessage{\string\pagebreak\space
     not allowed in non-display math mode}\else\postdisplaypenalty-\@M\fi
   \else\ifvmode\penalty-\@M\else\edef\spacefactor@
       {\spacefactor=\the\spacefactor}\vadjust{\penalty-\@M}\spacefactor@
	\relax\fi\fi}
\def\nopagebreak{\relax\ifmmode
     \ifinner\errmessage{\string\nopagebreak\space
    not allowed in non-display math mode}\else\postdisplaypenalty\@M\fi
    \else\ifvmode\nobreak\else\edef\spacefactor@
	{\spacefactor=\the\spacefactor}\vadjust{\penalty\@M}\spacefactor@
	 \relax\fi\fi}
\def\newpage{\relax\ifvmode\vfill\penalty-\@M\else\nonvmodeerr@\newpage\fi}
\def\nonvmodeerr@#1{\errmessage
    {\string#1\space is allowed only between paragraphs}}
\def\smallpagebreak{\relax\ifvmode\smallbreak
      \else\nonvmodeerr@\smallpagebreak\fi}
\def\medpagebreak{\relax\ifvmode\medbreak
       \else\nonvmodeerr@\medpagebreak\fi}
\def\bigpagebreak{\relax\ifvmode\bigbreak
      \else\nonvmodeerr@\bigpagebreak\fi}
\newdimen\captionwidth@
\captionwidth@=\hsize
\advance\captionwidth@ by -1.5in
\def\caption#1{}
\def\topspace#1{\gdef\thespace@{#1}\ifvmode\def\next
    {\futurelet\next\topspace@}\else\def\next{\nonvmodeerr@\topspace}\fi\next}
\def\topspace@{\ifx\next\Space@\def\next. {\futurelet\next\topspace@@}\else
     \def\next.{\futurelet\next\topspace@@}\fi\next.}
\def\topspace@@{\ifx\next\caption\let\next\topspace@@@\else
    \let\next\topspace@@@@\fi\next}
 \def\topspace@@@@{\topinsert\vbox to
       \thespace@{}\endinsert}
\def\topspace@@@\caption#1{\topinsert\vbox to
    \thespace@{}\nobreak
      \smallskip
    \setbox\z@=\hbox{\noindent\ignorespaces#1\unskip}%
   \ifdim\wd\z@>\captionwidth@
   \centerline{\vbox{\hsize=\captionwidth@\noindent\ignorespaces#1\unskip}}%
   \else\centerline{\box\z@}\fi\endinsert}
\def\midspace#1{\gdef\thespace@{#1}\ifvmode\def\next
    {\futurelet\next\midspace@}\else\def\next{\nonvmodeerr@\midspace}\fi\next}
\def\midspace@{\ifx\next\Space@\def\next. {\futurelet\next\midspace@@}\else
     \def\next.{\futurelet\next\midspace@@}\fi\next.}
\def\midspace@@{\ifx\next\caption\let\next\midspace@@@\else
    \let\next\midspace@@@@\fi\next}
 \def\midspace@@@@{\midinsert\vbox to
       \thespace@{}\endinsert}
\def\midspace@@@\caption#1{\midinsert\vbox to
    \thespace@{}\nobreak
      \smallskip
      \setbox\z@=\hbox{\noindent\ignorespaces#1\unskip}%
      \ifdim\wd\z@>\captionwidth@
    \centerline{\vbox{\hsize=\captionwidth@\noindent\ignorespaces#1\unskip}}%
    \else\centerline{\box\z@}\fi\endinsert}
\mathchardef\prime@="0230
\def\prime{{{}\prime@{}}}
\def\prim@s{\prime@\futurelet\next\pr@m@s}
\let\dsize=\displaystyle

\def\,{\relax\ifmmode\mskip\thinmuskip\else\thinspace\fi}
\def\!{\relax\ifmmode\mskip-\thinmuskip\else\negthinspace\fi}
\def\frac#1#2{{#1\over#2}}

\def\:{\nobreak\hskip.1111em{:}\hskip.3333em plus .0555em\relax}
\def\intic@{\mathchoice{\hskip5\p@}{\hskip4\p@}{\hskip4\p@}{\hskip4\p@}}
\def\negintic@
 {\mathchoice{\hskip-5\p@}{\hskip-4\p@}{\hskip-4\p@}{\hskip-4\p@}}
\def\intkern@{\mathchoice{\!\!\!}{\!\!}{\!\!}{\!\!}}
\def\intdots@{\mathchoice{\cdots}{{\cdotp}\mkern1.5mu
    {\cdotp}\mkern1.5mu{\cdotp}}{{\cdotp}\mkern1mu{\cdotp}\mkern1mu
      {\cdotp}}{{\cdotp}\mkern1mu{\cdotp}\mkern1mu{\cdotp}}}
\newcount\intno@
\def\iint{\intno@=\tw@\futurelet\next\ints@}
\def\iiint{\intno@=\thr@@\futurelet\next\ints@}
\def\iiiint{\intno@=4 \futurelet\next\ints@}
\def\idotsint{\intno@=\z@\futurelet\next\ints@}
\def\ints@{\findlimits@\ints@@}
\newif\iflimtoken@
\newif\iflimits@
\def\findlimits@{\limtoken@false\limits@false\ifx\next\limits
 \limtoken@true\limits@true\else\ifx\next\nolimits\limtoken@true\limits@false
    \fi\fi}
\def\multintlimits@{\intop\ifnum\intno@=\z@\intdots@
  \else\intkern@\fi
    \ifnum\intno@>\tw@\intop\intkern@\fi
     \ifnum\intno@>\thr@@\intop\intkern@\fi\intop}
\def\multint@{\int\ifnum\intno@=\z@\intdots@\else\intkern@\fi
   \ifnum\intno@>\tw@\int\intkern@\fi
    \ifnum\intno@>\thr@@\int\intkern@\fi\int}
\def\ints@@{\iflimtoken@\def\ints@@@{\iflimits@
   \negintic@\mathop{\intic@\multintlimits@}\limits\else
    \multint@\nolimits\fi\eat@}\else
     \def\ints@@@{\multint@\nolimits}\fi\ints@@@}
\def\Sb{_\bgroup\vspace@
	\baselineskip=\fontdimen10 \scriptfont\tw@
	\advance\baselineskip by \fontdimen12 \scriptfont\tw@
	\lineskip=\thr@@\fontdimen8 \scriptfont\thr@@
	\lineskiplimit=\thr@@\fontdimen8 \scriptfont\thr@@
	\Let@\vbox\bgroup\halign\bgroup \hfil$\scriptstyle
	    {##}$\hfil\cr}
\def\endSb{\crcr\egroup\egroup\egroup}
\def\Sp{^\bgroup\vspace@
	\baselineskip=\fontdimen10 \scriptfont\tw@
	\advance\baselineskip by \fontdimen12 \scriptfont\tw@
	\lineskip=\thr@@\fontdimen8 \scriptfont\thr@@
	\lineskiplimit=\thr@@\fontdimen8 \scriptfont\thr@@
	\Let@\vbox\bgroup\halign\bgroup \hfil$\scriptstyle
	    {##}$\hfil\cr}
\def\endSp{\crcr\egroup\egroup\egroup}
\def\Let@{\relax\iffalse{\fi\let\\=\cr\iffalse}\fi}
\def\vspace@{\def\vspace##1{\noalign{\vskip##1 }}}
\def\aligned{\,\vcenter\bgroup\vspace@\Let@\openup\jot\m@th\ialign
  \bgroup \strut\hfil$\displaystyle{##}$&$\displaystyle{{}##}$\hfil\crcr}
\def\endaligned{\crcr\egroup\egroup}
\def\matrix{\,\vcenter\bgroup\Let@\vspace@
    \normalbaselines
  \m@th\ialign\bgroup\hfil$##$\hfil&&\quad\hfil$##$\hfil\crcr
    \mathstrut\crcr\noalign{\kern-\baselineskip}}
\def\endmatrix{\crcr\mathstrut\crcr\noalign{\kern-\baselineskip}\egroup
		\egroup\,}
\newtoks\hashtoks@
\hashtoks@={#}
\def\format{\crcr\egroup\iffalse{\fi\ifnum`}=0 \fi\format@}
\def\format@#1\\{\def\preamble@{#1}%
  \def\c{\hfil$\the\hashtoks@$\hfil}%
  \def\r{\hfil$\the\hashtoks@$}%
  \def\l{$\the\hashtoks@$\hfil}%
  \setbox\z@=\hbox{\xdef\Preamble@{\preamble@}}\ifnum`{=0 \fi\iffalse}\fi
   \ialign\bgroup\span\Preamble@\crcr}

\def\cases{\left\{\,\vcenter\bgroup\vspace@
     \normalbaselines\openup\jot\m@th
       \Let@\ialign\bgroup$##$\hfil&\quad$##$\hfil\crcr
      \mathstrut\crcr\noalign{\kern-\baselineskip}}

\newif\iftagsleft@
\tagsleft@true
\def\TagsOnRight{\global\tagsleft@false}
\def\tag#1$${\iftagsleft@\leqno\else\eqno\fi
 \hbox{\def\pagebreak{\global\postdisplaypenalty-\@M}%
 \def\nopagebreak{\global\postdisplaypenalty\@M}\rm(#1\unskip)}%
  $$\postdisplaypenalty\z@\ignorespaces}
\interdisplaylinepenalty=\@M
\def\allowdisplaybreak@{\def\allowdisplaybreak{\noalign{\allowbreak}}}
\def\displaybreak@{\def\displaybreak{\noalign{\break}}}
\def\align#1\endalign{\def\tag{&}\vspace@\allowdisplaybreak@\displaybreak@
  \iftagsleft@\lalign@#1\endalign\else
   \ralign@#1\endalign\fi}
\def\ralign@#1\endalign{\displ@y\Let@\tabskip\centering\halign to\displaywidth
     {\hfil$\displaystyle{##}$\tabskip=\z@&$\displaystyle{{}##}$\hfil
       \tabskip=\centering&\llap{\hbox{(\rm##\unskip)}}\tabskip\z@\crcr
	     #1\crcr}}
\def\lalign@
 #1\endalign{\displ@y\Let@\tabskip\centering\halign to \displaywidth
   {\hfil$\displaystyle{##}$\tabskip=\z@&$\displaystyle{{}##}$\hfil
   \tabskip=\centering&\kern-\displaywidth
	\rlap{\hbox{(\rm##\unskip)}}\tabskip=\displaywidth\crcr
	       #1\crcr}}
\def\overrightarrow{\mathpalette\overrightarrow@}
\def\overrightarrow@#1#2{\vbox{\ialign{$##$\cr
    #1{-}\mkern-6mu\cleaders\hbox{$#1\mkern-2mu{-}\mkern-2mu$}\hfill
     \mkern-6mu{\to}\cr
     \noalign{\kern -1\p@\nointerlineskip}
     \hfil#1#2\hfil\cr}}}
\def\overleftarrow{\mathpalette\overleftarrow@}
\def\overleftarrow@#1#2{\vbox{\ialign{$##$\cr
     #1{\leftarrow}\mkern-6mu\cleaders\hbox{$#1\mkern-2mu{-}\mkern-2mu$}\hfill
      \mkern-6mu{-}\cr
     \noalign{\kern -1\p@\nointerlineskip}
     \hfil#1#2\hfil\cr}}}
\def\overleftrightarrow{\mathpalette\overleftrightarrow@}
\def\overleftrightarrow@#1#2{\vbox{\ialign{$##$\cr
     #1{\leftarrow}\mkern-6mu\cleaders\hbox{$#1\mkern-2mu{-}\mkern-2mu$}\hfill
       \mkern-6mu{\to}\cr
    \noalign{\kern -1\p@\nointerlineskip}
      \hfil#1#2\hfil\cr}}}
\def\underrightarrow{\mathpalette\underrightarrow@}
\def\underrightarrow@#1#2{\vtop{\ialign{$##$\cr
    \hfil#1#2\hfil\cr
     \noalign{\kern -1\p@\nointerlineskip}
    #1{-}\mkern-6mu\cleaders\hbox{$#1\mkern-2mu{-}\mkern-2mu$}\hfill
     \mkern-6mu{\to}\cr}}}
\def\underleftarrow{\mathpalette\underleftarrow@}
\def\underleftarrow@#1#2{\vtop{\ialign{$##$\cr
     \hfil#1#2\hfil\cr
     \noalign{\kern -1\p@\nointerlineskip}
     #1{\leftarrow}\mkern-6mu\cleaders\hbox{$#1\mkern-2mu{-}\mkern-2mu$}\hfill
      \mkern-6mu{-}\cr}}}
\def\underleftrightarrow{\mathpalette\underleftrightarrow@}
\def\underleftrightarrow@#1#2{\vtop{\ialign{$##$\cr
      \hfil#1#2\hfil\cr
    \noalign{\kern -1\p@\nointerlineskip}
     #1{\leftarrow}\mkern-6mu\cleaders\hbox{$#1\mkern-2mu{-}\mkern-2mu$}\hfill
       \mkern-6mu{\to}\cr}}}
\def\sqrt#1{\radical"270370 {#1}}
\def\dots{\relax\ifmmode\let\next=\ldots\else\let\next=\tdots@\fi\next}
\def\tdots@{\unskip\ \tdots@@}
\def\tdots@@{\futurelet\next\tdots@@@}
\def\tdots@@@{$\mathinner{\ldotp\ldotp\ldotp}\,
   \ifx\next,$\else
   \ifx\next.\,$\else
   \ifx\next;\,$\else
   \ifx\next:\,$\else
   \ifx\next?\,$\else
   \ifx\next!\,$\else
   $ \fi\fi\fi\fi\fi\fi}
\def\text{\relax\ifmmode\let\next=\text@\else\let\next=\text@@\fi\next}
\def\text@@#1{\hbox{#1}}
\def\text@#1{\mathchoice
 {\hbox{\everymath{\displaystyle}\def\textfonti{\the\textfont1 }%
    \def\textfontii{\the\textfont2 }\textdef@@ T#1}}
 {\hbox{\everymath{\textstyle}\def\textfonti{\the\textfont1 }%
    \def\textfontii{\the\textfont2 }\textdef@@ T#1}}
 {\hbox{\everymath{\scriptstyle}\def\textfonti{\the\scriptfont1 }%
   \def\textfontii{\the\scriptfont2 }\textdef@@ S\rm#1}}
 {\hbox{\everymath{\scriptscriptstyle}\def\textfonti{\the\scriptscriptfont1 }%
   \def\textfontii{\the\scriptscriptfont2 }\textdef@@ s\rm#1}}}
\def\textdef@@#1{\textdef@#1\rm \textdef@#1\bf
   \textdef@#1\sl \textdef@#1\it}

\def\textdef@#1#2{\def\next{\csname\expandafter\eat@\string#2fam\endcsname}%
\if S#1\edef#2{\the\scriptfont\next\relax}%
 \else\if s#1\edef#2{\the\scriptscriptfont\next\relax}%
 \else\edef#2{\the\textfont\next\relax}\fi\fi}
\scriptfont\itfam=\tenit \scriptscriptfont\itfam=\tenit
\scriptfont\slfam=\tensl \scriptscriptfont\slfam=\tensl
\mathcode`\0="0030
\mathcode`\1="0031
\mathcode`\2="0032
\mathcode`\3="0033
\mathcode`\4="0034
\mathcode`\5="0035
\mathcode`\6="0036
\mathcode`\7="0037
\mathcode`\8="0038
\mathcode`\9="0039
\def\Cal{\relax\ifmmode\let\next=\Cal@\else
     \def\next{\errmessage{Use \string\Cal\space only in math mode}}\fi\next}
\def\Cal@#1{{\fam2 #1}}
\def\bold{\relax\ifmmode\let\next=\bold@\else
   \def\next{\errmessage{Use \string\bold\space only in math
      mode}}\fi\next}\def\bold@#1{{\fam\bffam #1}}
\mathchardef\Gamma="0000
\mathchardef\Delta="0001
\mathchardef\Theta="0002
\mathchardef\Lambda="0003
\mathchardef\Xi="0004
\mathchardef\Pi="0005
\mathchardef\Sigma="0006
\mathchardef\Upsilon="0007
\mathchardef\Phi="0008
\mathchardef\Psi="0009
\mathchardef\Omega="000A
\mathchardef\varGamma="0100
\mathchardef\varDelta="0101
\mathchardef\varTheta="0102
\mathchardef\varLambda="0103
\mathchardef\varXi="0104
\mathchardef\varPi="0105
\mathchardef\varSigma="0106
\mathchardef\varUpsilon="0107
\mathchardef\varPhi="0108
\mathchardef\varPsi="0109
\mathchardef\varOmega="010A
\font\dummyft@=dummy
\fontdimen1 \dummyft@=\z@
\fontdimen2 \dummyft@=\z@
\fontdimen3 \dummyft@=\z@
\fontdimen4 \dummyft@=\z@
\fontdimen5 \dummyft@=\z@
\fontdimen6 \dummyft@=\z@
\fontdimen7 \dummyft@=\z@
\fontdimen8 \dummyft@=\z@
\fontdimen9 \dummyft@=\z@
\fontdimen10 \dummyft@=\z@
\fontdimen11 \dummyft@=\z@
\fontdimen12 \dummyft@=\z@
\fontdimen13 \dummyft@=\z@
\fontdimen14 \dummyft@=\z@
\fontdimen15 \dummyft@=\z@
\fontdimen16 \dummyft@=\z@
\fontdimen17 \dummyft@=\z@
\fontdimen18 \dummyft@=\z@
\fontdimen19 \dummyft@=\z@
\fontdimen20 \dummyft@=\z@
\fontdimen21 \dummyft@=\z@
\fontdimen22 \dummyft@=\z@
\def\fontlist@{\\{\tenrm}\\{\sevenrm}\\{\fiverm}\\{\teni}\\{\seveni}%
 \\{\fivei}\\{\tensy}\\{\sevensy}\\{\fivesy}\\{\tenex}\\{\tenbf}\\{\sevenbf}%
 \\{\fivebf}\\{\tensl}\\{\tenit}\\{\tensmc}}
\def\dodummy@{{\def\\##1{\global\let##1=\dummyft@}\fontlist@}}
\newif\ifsyntax@
\newcount\countxviii@
\def\newtoks@{\alloc@5\toks\toksdef\@cclvi}
\def\nopages@{\output={\setbox\z@=\box\@cclv \deadcycles=\z@}\newtoks@\output}
\def\syntax{\syntax@true\dodummy@\countxviii@=\count18
\loop \ifnum\countxviii@ > \z@ \textfont\countxviii@=\dummyft@
   \scriptfont\countxviii@=\dummyft@ \scriptscriptfont\countxviii@=\dummyft@
     \advance\countxviii@ by-\@ne\repeat
\dummyft@\tracinglostchars=\z@
  \nopages@\frenchspacing\hbadness=\@M}
\def\magstep#1{\ifcase#1 1000\or
 1200\or 1440\or 1728\or 2074\or 2488\or
 \errmessage{\string\magstep\space only works up to 5}\fi\relax}
{\lccode`\2=`\p \lccode`\3=`\t
 \lowercase{\gdef\tru@#123{#1truept}}}

\def\scaletype#1{\mag=#1\relax
 \hsize=\expandafter\tru@\the\hsize
 \vsize=\expandafter\tru@\the\vsize
 \dimen\footins=\expandafter\tru@\the\dimen\footins}

\def\scalefont#1#2\andcallit#3{\edef\font@{\the\font}#1\font#3=
  \fontname\font\space scaled #2\relax\font@}
\def\Mag@#1#2{\ifdim#1<1pt\multiply#1 #2\relax\divide#1 1000 \else
  \ifdim#1<10pt\divide#1 10 \multiply#1 #2\relax\divide#1 100\else
  \divide#1 100 \multiply#1 #2\relax\divide#1 10 \fi\fi}
\def\scalelinespacing#1{\Mag@\baselineskip{#1}\Mag@\lineskip{#1}%
  \Mag@\lineskiplimit{#1}}
\def\wlog#1{\immediate\write-1{#1}}
\catcode`\@=\active

\catcode`@=11
\def\binrel@#1{\setbox\z@\hbox{\thinmuskip0mu
\medmuskip\m@ne mu\thickmuskip\@ne mu$#1\m@th$}%
\setbox\@ne\hbox{\thinmuskip0mu\medmuskip\m@ne mu\thickmuskip
\@ne mu${}#1{}\m@th$}%
\setbox\tw@\hbox{\hskip\wd\@ne\hskip-\wd\z@}}
\def\overset#1\to#2{\binrel@{#2}\ifdim\wd\tw@<\z@
\mathbin{\mathop{\kern\z@#2}\limits^{#1}}\else\ifdim\wd\tw@>\z@
\mathrel{\mathop{\kern\z@#2}\limits^{#1}}\else
{\mathop{\kern\z@#2}\limits^{#1}}{}\fi\fi}
\def\underset#1\to#2{\binrel@{#2}\ifdim\wd\tw@<\z@
\mathbin{\mathop{\kern\z@#2}\limits_{#1}}\else\ifdim\wd\tw@>\z@
\mathrel{\mathop{\kern\z@#2}\limits_{#1}}\else
{\mathop{\kern\z@#2}\limits_{#1}}{}\fi\fi}
\def\circle#1{\leavevmode\setbox0=\hbox{h}\dimen@=\ht0
\advance\dimen@ by-1ex\rlap{\raise1.5\dimen@\hbox{\char'27}}#1}
\def\sqr#1#2{{\vcenter{\hrule height.#2pt
     \hbox{\vrule width.#2pt height#1pt \kern#1pt
       \vrule width.#2pt}
     \hrule height.#2pt}}}
\def\square{\mathchoice\sqr34\sqr34\sqr{2.1}3\sqr{1.5}3}
\def\force{\hbox{$\|\hskip-2pt\hbox{--}$\hskip2pt}}
 
\catcode`@=\active
%.....
%mathchar.tex
%.....
 
\catcode`\@=11
\def\bold{\relaxnext@\ifmmode\let\next\bold@\else
 \def\next{\Err@{Use \string\bold\space only in math mode}}\fi\next}
\def\bold@#1{{\bold@@{#1}}}
\def\bold@@#1{\fam\bffam#1}
\def\hexnumber@#1{\ifnum#1<10 \number#1\else
 \ifnum#1=10 A\else\ifnum#1=11 B\else\ifnum#1=12 C\else
 \ifnum#1=13 D\else\ifnum#1=14 E\else\ifnum#1=15 F\fi\fi\fi\fi\fi\fi\fi}
\def\bffam@{\hexnumber@\bffam}
%\mathchardef\boldGamma="0\bffam@00
%\mathchardef\boldDelta="0\bffam@01
%\mathchardef\boldTheta="0\bffam@02
%\mathchardef\boldLambda="0\bffam@03
%\mathchardef\boldXi="0\bffam@04
%\mathchardef\boldPi="0\bffam@05
%\mathchardef\boldSigma="0\bffam@06
%\mathchardef\boldUpsilon="0\bffam@07
%\mathchardef\boldPhi="0\bffam@08
%\mathchardef\boldPsi="0\bffam@09
%\mathchardef\boldOmega="0\bffam@0A
 
                              % change to msam and msbm files 1.1.92
 
\font\tenmsx=msam10
\font\sevenmsx=msam7
\font\fivemsx=msam5
\font\tenmsy=msbm10
\font\sevenmsy=msbm7
\font\fivemsy=msbm7
 
\newfam\msxfam
\newfam\msyfam
\textfont\msxfam=\tenmsx
\scriptfont\msxfam=\sevenmsx
\scriptscriptfont\msxfam=\fivemsx
\textfont\msyfam=\tenmsy
\scriptfont\msyfam=\sevenmsy
\scriptscriptfont\msyfam=\fivemsy
\def\msx@{\hexnumber@\msxfam}
\def\msy@{\hexnumber@\msyfam}
\mathchardef\boxdot="2\msx@00
\mathchardef\boxplus="2\msx@01
\mathchardef\boxtimes="2\msx@02
\mathchardef\square="0\msx@03
\mathchardef\blacksquare="0\msx@04
\mathchardef\centerdot="2\msx@05
\mathchardef\lozenge="0\msx@06
\mathchardef\blacklozenge="0\msx@07
\mathchardef\circlearrowright="3\msx@08
\mathchardef\circlearrowleft="3\msx@09
\mathchardef\rightleftharpoons="3\msx@0A
\mathchardef\leftrightharpoons="3\msx@0B
\mathchardef\boxminus="2\msx@0C
\mathchardef\Vdash="3\msx@0D
\mathchardef\Vvdash="3\msx@0E
\mathchardef\vDash="3\msx@0F
\mathchardef\twoheadrightarrow="3\msx@10
\mathchardef\twoheadleftarrow="3\msx@11
\mathchardef\leftleftarrows="3\msx@12
\mathchardef\rightrightarrows="3\msx@13
\mathchardef\upuparrows="3\msx@14
\mathchardef\downdownarrows="3\msx@15
\mathchardef\upharpoonright="3\msx@16

\mathchardef\downharpoonright="3\msx@17
\mathchardef\upharpoonleft="3\msx@18
\mathchardef\downharpoonleft="3\msx@19
\mathchardef\rightarrowtail="3\msx@1A
\mathchardef\leftarrowtail="3\msx@1B
\mathchardef\leftrightarrows="3\msx@1C
\mathchardef\rightleftarrows="3\msx@1D
\mathchardef\Lsh="3\msx@1E
\mathchardef\Rsh="3\msx@1F
\mathchardef\rightsquigarrow="3\msx@20
\mathchardef\leftrightsquigarrow="3\msx@21
\mathchardef\looparrowleft="3\msx@22
\mathchardef\looparrowright="3\msx@23
\mathchardef\circeq="3\msx@24
\mathchardef\succsim="3\msx@25
\mathchardef\gtrsim="3\msx@26
\mathchardef\gtrapprox="3\msx@27
\mathchardef\multimap="3\msx@28
\mathchardef\therefore="3\msx@29
\mathchardef\because="3\msx@2A
\mathchardef\doteqdot="3\msx@2B

\mathchardef\triangleq="3\msx@2C
\mathchardef\precsim="3\msx@2D
\mathchardef\lesssim="3\msx@2E
\mathchardef\lessapprox="3\msx@2F
\mathchardef\eqslantless="3\msx@30
\mathchardef\eqslantgtr="3\msx@31
\mathchardef\curlyeqprec="3\msx@32
\mathchardef\curlyeqsucc="3\msx@33
\mathchardef\preccurlyeq="3\msx@34
\mathchardef\leqq="3\msx@35
\mathchardef\leqslant="3\msx@36
\mathchardef\lessgtr="3\msx@37
\mathchardef\backprime="0\msx@38
\mathchardef\risingdotseq="3\msx@3A
\mathchardef\fallingdotseq="3\msx@3B
\mathchardef\succcurlyeq="3\msx@3C
\mathchardef\geqq="3\msx@3D
\mathchardef\geqslant="3\msx@3E
\mathchardef\gtrless="3\msx@3F
\mathchardef\sqsubset="3\msx@40
\mathchardef\sqsupset="3\msx@41
\mathchardef\vartriangleright="3\msx@42
\mathchardef\vartriangleleft ="3\msx@43
\mathchardef\trianglerighteq="3\msx@44
\mathchardef\trianglelefteq="3\msx@45
\mathchardef\bigstar="0\msx@46
\mathchardef\between="3\msx@47
\mathchardef\blacktriangledown="0\msx@48
\mathchardef\blacktriangleright="3\msx@49
\mathchardef\blacktriangleleft="3\msx@4A
\mathchardef\vartriangle="3\msx@4D
\mathchardef\blacktriangle="0\msx@4E
\mathchardef\triangledown="0\msx@4F
\mathchardef\eqcirc="3\msx@50
\mathchardef\lesseqgtr="3\msx@51
\mathchardef\gtreqless="3\msx@52
\mathchardef\lesseqqgtr="3\msx@53
\mathchardef\gtreqqless="3\msx@54
\mathchardef\Rrightarrow="3\msx@56
\mathchardef\Lleftarrow="3\msx@57
\mathchardef\veebar="2\msx@59
\mathchardef\barwedge="2\msx@5A
\mathchardef\doublebarwedge="2\msx@5B
\mathchardef\angle="0\msx@5C
\mathchardef\measuredangle="0\msx@5D
\mathchardef\sphericalangle="0\msx@5E
\mathchardef\varpropto="3\msx@5F
\mathchardef\smallsmile="3\msx@60
\mathchardef\smallfrown="3\msx@61
\mathchardef\Subset="3\msx@62
\mathchardef\Supset="3\msx@63
\mathchardef\Cup="2\msx@64

\mathchardef\Cap="2\msx@65

\mathchardef\curlywedge="2\msx@66
\mathchardef\curlyvee="2\msx@67
\mathchardef\leftthreetimes="2\msx@68
\mathchardef\rightthreetimes="2\msx@69
\mathchardef\subseteqq="3\msx@6A
\mathchardef\supseteqq="3\msx@6B
\mathchardef\bumpeq="3\msx@6C
\mathchardef\Bumpeq="3\msx@6D
\mathchardef\lll="3\msx@6E

\mathchardef\ggg="3\msx@6F

\mathchardef\circledS="0\msx@73
\mathchardef\pitchfork="3\msx@74
\mathchardef\dotplus="2\msx@75
\mathchardef\backsim="3\msx@76
\mathchardef\backsimeq="3\msx@77
\mathchardef\complement="0\msx@7B
\mathchardef\intercal="2\msx@7C
\mathchardef\circledcirc="2\msx@7D
\mathchardef\circledast="2\msx@7E
\mathchardef\circleddash="2\msx@7F
\def\ulcorner{\delimiter"4\msx@70\msx@70 }
\def\urcorner{\delimiter"5\msx@71\msx@71 }
\def\llcorner{\delimiter"4\msx@78\msx@78 }
\def\lrcorner{\delimiter"5\msx@79\msx@79 }
\def\yen{{\mathhexbox@\msx@55 }}
\def\checkmark{{\mathhexbox@\msx@58 }}
\def\circledR{{\mathhexbox@\msx@72 }}
\def\maltese{{\mathhexbox@\msx@7A }}
\mathchardef\lvertneqq="3\msy@00
\mathchardef\gvertneqq="3\msy@01
\mathchardef\nleq="3\msy@02
\mathchardef\ngeq="3\msy@03
\mathchardef\nless="3\msy@04
\mathchardef\ngtr="3\msy@05
\mathchardef\nprec="3\msy@06
\mathchardef\nsucc="3\msy@07
\mathchardef\lneqq="3\msy@08
\mathchardef\gneqq="3\msy@09
\mathchardef\nleqslant="3\msy@0A
\mathchardef\ngeqslant="3\msy@0B
\mathchardef\lneq="3\msy@0C
\mathchardef\gneq="3\msy@0D
\mathchardef\npreceq="3\msy@0E
\mathchardef\nsucceq="3\msy@0F
\mathchardef\precnsim="3\msy@10
\mathchardef\succnsim="3\msy@11
\mathchardef\lnsim="3\msy@12
\mathchardef\gnsim="3\msy@13
\mathchardef\nleqq="3\msy@14
\mathchardef\ngeqq="3\msy@15
\mathchardef\precneqq="3\msy@16
\mathchardef\succneqq="3\msy@17
\mathchardef\precnapprox="3\msy@18
\mathchardef\succnapprox="3\msy@19
\mathchardef\lnapprox="3\msy@1A
\mathchardef\gnapprox="3\msy@1B
\mathchardef\nsim="3\msy@1C
\mathchardef\napprox="3\msy@1D
\mathchardef\varsubsetneq="3\msy@20
\mathchardef\varsupsetneq="3\msy@21
\mathchardef\nsubseteqq="3\msy@22
\mathchardef\nsupseteqq="3\msy@23
\mathchardef\subsetneqq="3\msy@24
\mathchardef\supsetneqq="3\msy@25
\mathchardef\varsubsetneqq="3\msy@26
\mathchardef\varsupsetneqq="3\msy@27
\mathchardef\subsetneq="3\msy@28
\mathchardef\supsetneq="3\msy@29
\mathchardef\nsubseteq="3\msy@2A
\mathchardef\nsupseteq="3\msy@2B
\mathchardef\nparallel="3\msy@2C
\mathchardef\nmid="3\msy@2D
\mathchardef\nshortmid="3\msy@2E
\mathchardef\nshortparallel="3\msy@2F
\mathchardef\nvdash="3\msy@30
\mathchardef\nVdash="3\msy@31
\mathchardef\nvDash="3\msy@32
\mathchardef\nVDash="3\msy@33
\mathchardef\ntrianglerighteq="3\msy@34
\mathchardef\ntrianglelefteq="3\msy@35
\mathchardef\ntriangleleft="3\msy@36
\mathchardef\ntriangleright="3\msy@37
\mathchardef\nleftarrow="3\msy@38
\mathchardef\nrightarrow="3\msy@39
\mathchardef\nLeftarrow="3\msy@3A
\mathchardef\nRightarrow="3\msy@3B
\mathchardef\nLeftrightarrow="3\msy@3C
\mathchardef\nleftrightarrow="3\msy@3D
\mathchardef\divideontimes="2\msy@3E
\mathchardef\varnothing="0\msy@3F
\mathchardef\nexists="0\msy@40
\mathchardef\mho="0\msy@66
\mathchardef\thorn="0\msy@67
\mathchardef\beth="0\msy@69
\mathchardef\gimel="0\msy@6A
\mathchardef\daleth="0\msy@6B
\mathchardef\lessdot="3\msy@6C
\mathchardef\gtrdot="3\msy@6D
\mathchardef\ltimes="2\msy@6E
\mathchardef\rtimes="2\msy@6F
\mathchardef\shortmid="3\msy@70
\mathchardef\shortparallel="3\msy@71
\mathchardef\smallsetminus="2\msy@72
\mathchardef\thicksim="3\msy@73
\mathchardef\thickapprox="3\msy@74
\mathchardef\approxeq="3\msy@75
\mathchardef\succapprox="3\msy@76
\mathchardef\precapprox="3\msy@77
\mathchardef\curvearrowleft="3\msy@78
\mathchardef\curvearrowright="3\msy@79
\mathchardef\digamma="0\msy@7A
\mathchardef\varkappa="0\msy@7B
\mathchardef\hslash="0\msy@7D
\mathchardef\hbar="0\msy@7E
\mathchardef\backepsilon="3\msy@7F
\def\Bbb{\relaxnext@\ifmmode\let\next\Bbb@\else
 \def\next{\Err@{Use \string\Bbb\space only in math mode}}\fi\next}
\def\Bbb@#1{{\Bbb@@{#1}}}
\def\Bbb@@#1{\noaccents@\fam\msyfam#1}
\catcode`\@=12
%\font\teneuf=eufm10
%\font\fiveeuf=eufm5
%\font\seveneuf=eufm7
%\font\sc=cmcsc10
                       % msxm and msym fonts changed for msam and msbm fonts
                       %  01.01.92

\font\tenmsy=msbm10
\font\sevenmsy=msbm7
\font\fivemsy=msbm5
\font\tenmsx=msam10
\font\sevenmsx=msam7
\font\fivemsx=msam5
%\font\sevenmcyr=mcyr7
%\font\tenmcyr=mcyr10
%%\font\tenmcyb=mcyb10
%%\font\eightmcyb=mcyb8
\newfam\msyfam

\textfont\msyfam=\tenmsy
\scriptfont\msyfam=\sevenmsy
\scriptscriptfont\msyfam=\fivemsy
\newfam\msxfam

\textfont\msxfam=\tenmsx
\scriptfont\msxfam=\sevenmsx
\scriptscriptfont\msxfam=\fivemsx
%\newfam\mcyrfam
%\def\mcyr{\fam\mcyrfam\tenmcyr}
%\textfont\mcyrfam=\tenmcyr
%\scriptfont\mcyrfam=\sevenmcyr
%\newfam\euffam
%\def\euf{\fam\euffam\teneuf}
%\textfont\euffam=\teneuf
%\scriptfont\euffam=\seveneuf
%\scriptscriptfont\euffam=\fiveeuf

%
%this is def.tex
 
\magnification 1200
\baselineskip=24pt plus 6pt
\def\today{\ifcase\month\or January\or February\or March\or April\or
May\or June\or July\or August\or September\or October\or November\or
December \fi\space\number \day, \number\year}
\define\a{\alpha}

\define\egg{\Relbar\kern-.30em\Relbar\kern-.30em\Relbar}

\define\p1{P^1_{\delta, \lambda}}
\define\dom{\hbox{\rm dom}}
\define\pbf{\par\bigpagebreak\flushpar}
\define\g0{G^0_{\delta, \lambda}}
\define\cof{\hbox{\rm cof}}

\define\dell{\delta, \lambda}
\define\ao{_{\alpha_{0}} }
\define\pj{^{p_j}}

\def\b{\beta}
\def\a{\alpha}
\def\l{\lambda}

\def\gam{\gamma}

\def\del{\delta}
%\font\frac=eurb10
%\font\frak=eufm10

\def\k{\kappa}
\pageno=1
\magnification 1200
\baselineskip=24pt plus 6pt

\define\cU{\Cal U}

%\define\bA{\bold{A}}
\define\cD{\Cal{D}}

\def\b{\beta}
\def\a{\alpha}
\def\l{\lambda}

\def\gam{\gamma}

\def\del{\delta}
%\font\frac=eurb10
%\font\frak=eufm10
%\def\gB{{\frak\$B}}
%\def\gV{{\frak\$V}}
%\def\gI{{\frak\$I}}
%\def\lorop{\mathop\lor}
\def\k{\kappa}
\def\d{\delta}
\def\no{\noindent}
\def\g{\gamma}
\def\la{\langle}
\def\ra{\rangle}

\def\A{{\cal A}}
\def\B{{\cal B}}
\def\supp{{\hbox{\rm support}}}
\def\dom{{\hbox{\rm dom}}}
\def\min{{\hbox{\rm min}}}
\def\max{{\hbox{\rm max}}}
\def\cof{{\hbox{\rm cof}}}

\def\sqr#1#2{{\vcenter{\vbox{\hrule height.#2pt
       \hbox{\vrule width.#2pt height#1pt \kern#1pt
         \vrule width.#2pt}
        \hrule height.#2pt}}}}
\def\square{\mathchoice\sqr34\sqr34\sqr{2.1}3\sqr{1.5}3}
\def\finpf{\no\hfill $\square$ }

\voffset=-.5in\vsize=7.5in
 %%%%%%%%%%%%%%%%%%%%%%%%%%%%%%%%%%%%%%%%%%%%%%%%%%%%%%%%%%%%%%%%%%%%

\S0 Introduction and Preliminaries
 
It is a well known fact that the notion  of strongly compact
cardinal represents a singularity in the hierarchy of large cardinals.
The work of Magidor [Ma1] shows that the least strongly compact cardinal
and the least supercompact cardinal can coincide, but also, the least strongly
compact cardinal and the least measurable cardinal can coincide.
The work of Kimchi and Magidor [KiM] generalizes this, showing that
the class of 
strongly compact cardinals and the class of supercompact cardinals
can coincide (except by results of Menas [Me] and [A] at certain measurable
limits of supercompact cardinals), and the first $n$ strongly
compact cardinals (for $n$ a natural number) and the first $n$ measurable
cardinals can coincide.
Thus, the precise identity of certain members of the class of strongly compact
cardinals cannot be ascertained vis \`a vis the class of measurable
cardinals or the class of supercompact cardinals.
 
An interesting aspect of the proofs of both [Ma1] and [KiM] is that in
each result, all ``bad" instances of strong compactness are not obliterated.
Specifically, in each model, since the strategy employed in destroying strongly
 compact
 cardinals which aren't also supercompact is to make them non-strongly compact
after a certain point either by adding a  Prikry sequence or a non-reflecting
stationary set of ordinals of the appropriate cofinality, there may
be cardinals  $\k$ and $\l$ so that $\k$ is $\l$ strongly compact
yet $\k$ isn't $\l$ supercompact.
Thus, whereas it was proven by Kimchi and Magidor that the classes
of strongly compact and supercompact cardinals can coincide (with
the exceptions noted above), it was not known whether a ``local''
version of this were possible, i.e., if it were possible to
obtain a model in which, for the class of pairs $(\k,\l)$,
$\k$ is $\l$ strongly compact iff $\k$ is $\l$
supercompact. This is more delicate.
 
The purpose of this paper is to answer the above
question in the affirmative.
Specifically, we prove the following
 
\proclaim{Theorem} Suppose $ V \models $  ZFC + GCH is a
given model (which in interesting cases contains instances of
supercompactness). There is then some cardinal and cofinality
preserving generic extension $V[G] \models$ ZFC + GCH
in which: \hfil\break
(a) (Preservation) For $\k \le \l$ regular, if $V
\models ``\k$ is $\l$ supercompact'', then $V[G] \models ``\k$ is
$\l$ supercompact''. The converse implication holds except
possibly when $\k = \sup\{ \delta < \k : \delta$ is $\lambda$
supercompact$\}$.
                              \hfil\break
(b) (Equivalence) For $\k \le \l$ regular, $V[G]
\models ``\k$ is $\l$ strongly compact'' iff $V[G] \models
``\k$ is $\l$ supercompact'', except possibly if $\k$ is a
measurable limit of cardinals which are $\l  $ supercompact.
\endproclaim
 
Note that the limitation given in (b) above is reasonable, since
trivially, if $\k$ is measurable, $\k < \lambda$, and
$\k = \sup\{\delta < \k : \delta$ is either $\lambda$
supercompact or $\lambda$ strongly compact$\}$, then $\k$ is
$\lambda$ strongly compact. Further,
it is a theorem of Menas [Me] that under GCH, for $\k$ the first,
second, third, or ${\a}$th for $\a < \k$ measurable limit
of cardinals which are $\k^+$ strongly compact or
$\k^+$ supercompact, $\k$ is $\k^+$ strongly compact yet
$\k$ isn't
$\k^+$ supercompact.
Thus, if there are
sufficiently large cardinals in the universe, it will never be
possible to have a complete coincidence between the notions of
$\k$ being $\l$ strongly compact and $\k$ being $\l$
supercompact for $\l$ a regular cardinal.
 
Note that in the statement of our Theorem, we do not mention
what happens if $\l > \k$ is a singular cardinal. This is since
the behavior when $\l > \k$ is a singular cardinal is provable
in ZFC + GCH (which implies any limit cardinal is a strong
limit cardinal). Specifically, if $\l > \k$ is so that
cof($\l) < \k$, then by a theorem of Magidor [Ma3], $\k$ is
$\l$ supercompact iff $\k$ is $\l^+$ supercompact, so
automatically, by clause (a) of our Theorem, $\l$
supercompactness is preserved between $V$ and $V[G]$.
Also, if $\l > \k$ is so that cof($\l) < \k$, then by a
theorem of Solovay [SRK], $\k$ is $\l$ strongly compact iff
$\k$ is $\l^+$ strongly compact, so by clause (b) of our
Theorem, it can never be the case that $V[G]
\models ``\k$ is $\l$ strongly compact'' unless
$V[G]         \models ``\k$ is $\l$ supercompact'' as well.
Further, if $\l > \k$ is so that $\l >$ cof$(\l) \ge \k$,
then it is not too difficult to see (and will be shown in
Section 2) that if $\k$ is $\l'$ strongly compact or $\l'$
supercompact for all $\l' < \l$, then $\k$ is $\l$ strongly
compact, and there is no reason to believe $\k$ must be
$\l$ supercompact. In fact,
it is a theorem of Magidor [Ma4] (irrespective of
GCH) that if $\mu$ is a
supercompact cardinal, there will always be many cardinals
$\k,\l < \mu$ so that $\l > \k$ is a singular cardinal of
cofinality $\ge \k$, $\k$ is $\l$ strongly compact,
$\k$ is $\l'$ supercompact for all $\l' < \l$, yet
$\k$ isn't $\l$ supercompact. Thus, there can never be a
complete coincidence between the notions of $\k$ being
$\l$ strongly compact and $\k$ being $\l$ supercompact if
$\l > \k$ is an arbitrary cardinal, assuming there are
supercompact cardinals in the universe.

The structure of this paper is as follows.
Section 0 contains our introductory comments and preliminary material
concerning notation, terminology, etc.
Section 1 defines and discusses the basic properties of the
forcing notion used in the iteration we employ to construct
our final model.
Section 2 gives a complete statement and proof of the theorem of Magidor
mentioned in the above  paragraph
 and proves  our Theorem in the case for which there is one
supercompact cardinal $\k$ in the universe which contains
no strongly inaccessible cardinals above it.
Section 3 shows how the ideas of Section 2 can be used to prove the
Theorem in the general case.
Section 4 contains our concluding remarks.
 
Before beginning the material of Section 1, we briefly mention some
preliminary information.  Essentially, our notation and terminology are
standard, and when this is not the case, this will be clearly noted.  We take
this opportunity to mention we will be assuming GCH throughout the course of
this paper.  For $\a < \b$ ordinals, $[\a, \b], [\a, \b), (\a, \b]$, and $(\a,
\b) $ are as in standard interval notation.  If $f$ is the characteristic
function of a set $x \subseteq \a$, then $x = \{\b : f(\b) = 1 \}$.
 
When forcing, $q \ge p$ will mean that $q$ is stronger than $p$.  For $P$ a
partial ordering, $\varphi$ a formula in the forcing language with respect to
$P$, and $ p \in P$, $ p \| \varphi$ will mean $p$ decides $\varphi$.  For $G$
$V$-generic over $P$, we will use both $V[G]$ and $V^{P}$ to indicate the
universe obtained by forcing with $P$.  If $x \in V[G]$, then $\dot x$ will be
a term in $V$ for $x$.  We may, from time to time, confuse terms with the sets
they denote and write $x$ when we actually mean $\dot x$, especially when $x$
is some variant of the generic set $G$.
 
If $\k$ is a cardinal, then for
$P$ a partial ordering, $P$ is $(\k, \infty)$-distributive if
for any sequence $\la D_\a : \a < \k \ra$ of dense open
subsets of $P$, $D = \underset \a < \k \to{\cap} D_\a$
is a dense open subset of $P$.
$P$ is $\k$-closed if given a sequence
$\langle p_\a: \a < \k \rangle$ of elements of $P$ so that
$\beta < \gamma < \k$ implies $p_\beta \le p_\gamma$ (an increasing chain of
 length
 $\k$), then there is some $p \in P$ (an upper bound to this chain) so that
$p_\a \le p$ for all $\a < \k$.
$P$ is $<\k$-closed if $P$ is $\delta$-closed for all cardinals $\delta <
\k$.
$P$ is $\k$-directed closed if for every cardinal $\delta < \k$ and every
 directed
 set $\langle p_\a : \a < \delta \rangle $ of elements of $P$
(where $\langle p_\alpha : \alpha < \delta \rangle$ is directed if
for every two distinct elements $p_\rho, p_\nu \in
\langle p_\alpha : \alpha < \delta \rangle$, $p_\rho$ and
$p_\nu$ have a common upper bound) there is an
upper bound $p \in P$. $P$ is $\k$-strategically closed if in the
two person game in which the players construct an increasing sequence
 $\langle p_\a: \a \le\k\rangle$, where player I plays odd stages and player
II plays even and limit stages, then player II has a strategy which ensures the
 game
 can always be continued.
$ P$ is $< \k$-strategically closed if $P$ is $\delta$-strategically
 closed for all cardinals $\delta < \k$.
$P$ is $\prec   \k$-strategically closed if in the two
person game in which the players construct an increasing
sequence $\langle p_\alpha : \alpha < \k \rangle$, where
player I plays odd stages and player II plays even and limit
stages, then player II has a strategy which ensures the game
can always be continued.
Note that trivially, if $P$ is $\k$-closed, then $P$ is $\k$-strategically
closed and $\prec \k^+$-strategically closed. The converse of
both of these facts is false.
 
For $\k$ a regular cardinal, two partial orderings to  which we will refer
quite a bit are the standard partial orderings $Q^0_\k$ for adding a Cohen
 subset to
$\k^+$ using conditions having support $\k$ and $Q^1_\k$  for adding
$\k^+$ many Cohen subsets to $\k$ using conditions having
support $<\k$.
The basic properties and explicit definitions of these partial orderings
may be found in [J].
 
Finally, we mention that we are assuming complete familiarity with the notions
 of
strong compactness and supercompactness.
Interested readers may consult [SRK] or [KaM] for further details.
We note only that all elementary embeddings witnessing the $\lambda$
supercompactness  of $\k$ are presumed to come from some
fine, $\k$-complete, normal
ultrafilter $\cU$ over $P_\k (\l) = \{ x \subseteq \l: | x| < \k \}$.
Also, where appropriate, all ultrapowers via a supercompact ultrafilter over
$P_\k (\l)$ will be confused with their transitive isomorphs.
 
\S 1 The Forcing Conditions
 
In this section, we describe and prove the basic properties of the forcing
conditions we shall use in our later iteration.  Let $\delta < \l, \, \l \ge
\aleph_1$ be regular cardinals in our ground model $V$.  We define three
notions of forcing.  Our first notion of forcing $P^0_{\delta, \l}$ is just
the standard notion of forcing for adding a non-reflecting stationary set of
ordinals of cofinality $\delta $ to $\l^+$.  Specifically, $P^0_{\delta,\l} =
\{ p$ : For some $\a < \l^+$, $p : \a \to \{0,1\}$ is a characteristic
function of $S_p$, a subset of $\a$ not stationary at its supremum nor having
any initial segment which is stationary at its supremum, so that $\b \in S_p$
implies $\b > \d$ and $\cof(\b) = \d \}$, ordered by $q \ge p$ iff $q
\supseteq p$ and $S_p = S_q \cap \sup (S_p)$, i.e., $S_q$ is an end extension
of $S_p$. It is well-known that for $G$ $V$-generic over $P^0_{\delta, \l}$
(see [Bu] or [KiM]), in $V[G]$, a non-reflecting stationary set
$S=S[G]=\cup\{S_p:p\in G \} \subseteq \l^+$ of ordinals of cofinality $\delta
$ has been introduced, the bounded subsets of $\l^+$ are the same as those in
$V$, and cardinals, cofinalities, and GCH have been preserved.  It is also
virtually immediate that $P^0_{\d, \l}$ is $\d$-directed closed.
 
Work now in $V_1 = V^{P^0_{\delta, \l}}$, letting $\dot S$ be a term always
forced to denote the above set $S$. \ $P^2_{\delta, \l}[S]$ is the standard
notion of forcing for introducing a club set $C$ which is disjoint to $S$ (and
therefore makes $S$ non-stationary).  Specifically, $P^2_{\delta, \l} [S] = \{
p$ : For some successor ordinal $\a < \l^+$, $p : \a \to \{0,1\}$ is a
characteristic function of $C_p$, a club subset of $\a$, so that $C_p \cap S =
\emptyset \}$, ordered by $ q \ge p $ iff $C_q$ is an end extension of $C_p$.
It is again well-known (see [MS]) that for $H$ $V_1$-generic over
$P^2_{\delta, \l}[S]$, a club set $C = C[H] = \cup \{C_p : p \in H \}\subseteq
\l^+$ which is disjoint to $S$ has been introduced, the bounded subsets 
 of $\l^+$ are the same as those in $V_1$, and cardinals, cofinalities, and
GCH have been preserved.
 
Before defining in $V_1$ the partial ordering
$P^1_{\delta, \l}[S]$ which will be used to destroy strong
compactness, we first prove two preliminary lemmas.
 
\proclaim{Lemma 1} $\force_{P^0_{\delta, \l}} ``
\clubsuit({\dot S})$'', i.e., $V_1 \models ``$There is a
sequence $\langle x_\alpha : \alpha \in S \rangle$ so that
for each $\alpha \in S$, $ x_\alpha \subseteq \alpha$ is
cofinal in
$\alpha$, and for any $A \in
{[\l^+]}^{\l^+}$, $\{\alpha \in S : x_\alpha
\subseteq A      \}$ is stationary''.
\endproclaim
 
\demo{Proof of Lemma 1} Since $V \models$ GCH and $V$ and
$V_1$ contain the same bounded subsets of $\l^+$, we can let
$\langle y_\alpha : \alpha < \l^+ \rangle \in V$ be a listing
of all elements $x \in {({[\l^+]}^\delta)}^V =
{({[\l^+]}^\delta)}^{V_1}$ so that each $x \in
{[\l^+]}^\delta$ appears on this list $\l^+$ times at ordinals
of cofinality $\delta$, i.e., for any $x \in {[\l^+]}^\delta$,
$\l^+ = \sup \{\a < \l^+$ : cof$(\a) = \delta$ and
$y_\a = x \}$. This then
allows us to define $\langle x_\alpha : \alpha \in
S \rangle$ by letting $x_\alpha$ be $y_\beta$ for the least
$\beta \in S - (\alpha + 1)$ so that $y_\beta \subseteq
\alpha$ and $y_\beta$ is unbounded in $\alpha$. By genericity,
each $x_\alpha$ is well-defined.
 
Let now $p \in P^0_{\delta, \l}$ be so that $p \force
``{\dot A} \in {[\l^+]}^{\l^+}$ and ${\dot K} \subseteq
\l^+$ is club''. We show that for some $r \ge p$ and some
$\zeta < \l^+$, $r \force ``\zeta \in {\dot K} \cap
{\dot S}$ and ${\dot x}_\zeta \subseteq {\dot A}$''. To do
this, we inductively define an increasing sequence $\langle
p_\alpha : \alpha < \delta \rangle$ of elements of
$P^0_{\delta, \l}$ and increasing sequences $\langle
\beta_\alpha : \alpha < \delta \rangle$ and $\langle
\gamma_\alpha : \alpha < \delta \rangle$ of ordinals
$< \l^+$ so that $\beta_0 \le \gamma_0 \le \beta_1 \le
\gamma_1 \le \cdots \le \beta_\alpha \le \gamma_\alpha
\le \cdots$ $(\alpha < \delta)$.
We begin by letting $p_0 = p$ and $\beta_0 = \gamma_0 =
0$. For $\eta = \alpha + 1 < \delta$ a successor, let
$p_\eta \ge p_\alpha$ and $\beta_\eta \le \gamma_\eta$,
$\beta_\eta \ge$max$(\beta_\alpha, \gamma_\alpha,
\sup($dom$(p_\alpha))) + 1$ be so that $p_\eta \force
``\beta_\eta \in {\dot A}$ and $\gamma_\eta \in
{\dot K}$''. For $\rho < \delta$ a limit, let $p_\rho =
\underset \alpha < \rho \to {\cup} p_\alpha$,
$\beta_\rho = \underset \alpha < \rho \to {\cup} \beta_\alpha$,
and $\gamma_\rho = \underset \alpha < \rho \to {\cup}
\gamma_\alpha$. Note that since $\rho < \delta$, $p_\rho$ is
well-defined, and since $\delta < \l^+$, $\beta_\rho,
\gamma_\rho < \l^+$. Also, by construction,
$\underset \alpha < \delta \to {\cup} \beta_\alpha =
\underset \alpha < \delta \to {\cup} \gamma_\alpha =
\underset \alpha < \delta \to {\cup} \sup($dom$(p_\alpha))
< \l^+$.
Call $\zeta$ this common sup. We thus have that $q =
\underset \alpha < \delta \to {\bigcup} p_\alpha \cup
\{\zeta\}$ is a well-defined condition so that $q \force
``\{\beta_\alpha : \alpha \in \delta - \{0\} \}
\subseteq \dot A$ and $\zeta \in {\dot K} \cap
{\dot S}$''.
 
To complete the proof of Lemma 1, we know that as $\langle
\beta_\alpha : \alpha \in \delta - \{0\} \rangle \in V$ and
as each $y \in \langle y_\a : \a < \l^+ \rangle        $ must
appear $\l^+$ times at ordinals of cofinality $\delta$,
we can find some $\eta \in (\zeta,
\l^+)$ so that cof$(\eta) = \delta$ and
$\langle \beta_\alpha : \alpha \in \delta -
\{0\} \rangle = y_\eta$. If we let $r \ge q$ be so that
$r \force ``{\dot S} \cap [\zeta, \eta] = \{\zeta, \eta\}$'',
then $r \force ``{\dot x}_\zeta = y_\eta =
\langle \beta_\alpha : \alpha \in \delta - \{0\}
\rangle$''. This proves Lemma 1.
 
\noindent \hfill $\square$ Lemma 1
 
We fix now in $V_1$ a $\clubsuit(S)$ sequence
$X = \la x_\a : \a \in S \ra$.
 
\proclaim{Lemma 2}
Let $S'$ be an initial segment of $S$ so that $S'$
is not stationary at its supremum nor has any initial
segment which is stationary at its supremum. There is then
a sequence $\la y_\a : \a \in S' \ra$ so that for every
$\a \in S'$, $y_\a \subseteq x_\a$, $x_\a - y_\a$ is
bounded in $\a$, and if $\a_1 \neq \a_2 \in S'$, then
$y_{\a_1} \cap y_{\a_2} = \emptyset$.
\endproclaim
 
\demo{Proof of Lemma 2}
We define by induction on $\a \le \a_0 = \sup S' + 1$
a function $h_\a$ so that $\dom(h_\a) = S' \cap \a$,
$h_\a(\b) < \b$, and $\la x_\b - h_\a(\b) :
\b \in S' \cap \a \ra$ is pairwise disjoint. The sequence
$\la x_\b - h_{\a_0}(\b) : \b \in S' \ra$ will be our
desired sequence.
 
If $\a = 0$, then we take $h_\a$ to be the empty function.
If $\a = \b + 1$ and $\b \not\in S'$, then we take
$h_\a = h_\b$. If $\a = \b + 1$ and $\b \in S'$, then
we notice that since each $x_\g \in X$ has order type
$\d$ and is cofinal in $\g$, for all $\g \in S' \cap
\b$, $x_\b \cap \g$ is bounded in $\g$. This allows us to
define a function $h_\a$ having domain $S' \cap \a$ by
$h_\a(\b) = 0$, and for $\g \in S' \cap \b$,
$h_\a(\g) = \min(\{\rho : \rho < \g$, $\rho \ge
h_\b(\g)$, and $x_\b \cap \g \subseteq \rho \})$. By the
next to last sentence and the induction hypothesis on
$h_\b$, $h_\a(\g) < \g$. And, if $\g_1 < \g_2 \in
S' \cap \a$, then if $\g_2 < \b$, $(x_{\g_1} -
h_\a(\g_1)) \cap (x_{\g_2} - h_\a(\g_2)) \subseteq
(x_{\g_1} - h_\b(\g_1)) \cap (x_{\g_2} -
h_\b(\g_2)) = \emptyset$ by the induction
hypothesis on $h_\b$.
If $\g_2 = \b$, then $(x_{\g_1} - h_\a(\g_1)) \cap
(x_{\g_2} - h_\a(\g_2)) = (x_{\g_1} - h_\a(\g_1)) \cap
x_{\g_2} = \emptyset$ by the definition of $h_\a(\g_1)$.
The sequence $\la x_\g - h_\a(\g) : \g \in S' \cap \a \ra$
is thus as desired.
 
If $\a$ is a limit ordinal, then as $S'$ is non-stationary
at its supremum nor has any initial segment which is
stationary at its supremum, we can let
$\la \b_\g : \g < \cof(\a) \ra$ be a strictly increasing,
continuous sequence having sup $\a$ so that for all
$\g < \cof(\a)$, $\b_\g \not\in S'$. Thus, if $\rho \in
S' \cap \a$, then $\{\b_\g : \b_\g < \rho \}$ is
bounded in $\rho$, meaning we can find some largest
$\g$ so that $\b_\g < \rho$. It is also the case that
$\rho < \b_{\g + 1}$. This allows us to define
$h_\a(\rho) = \max(\{h_{\b_{\g + 1}}(\rho),
\b_\g \})$ for the $\g$ just described. It is still the
case that $h_\a(\rho) < \rho$. And, if $\rho_1,
\rho_2 \in (\b_\g, \b_{\g + 1})$, then
$(x_{\rho_1} - h_\a(\rho_1)) \cap
 (x_{\rho_2} - h_\a(\rho_2)) \subseteq
(x_{\rho_1} - h_{\b_{\g + 1}}(\rho_1)) \cap
(x_{\rho_2} - h_{\b_{\g + 1}}(\rho_2))
= \emptyset$ by the definition of $h_{\b_{\g + 1}}$.
If $\rho_1 \in (\b_\g, \b_{\g + 1})$,
$\rho_2 \in (\b_\sigma, \b_{\sigma + 1})$ with
$\g < \sigma$, then
$(x_{\rho_1} - h_\a(\rho_1)) \cap
 (x_{\rho_2} - h_\a(\rho_2)) \subseteq
x_{\rho_1} \cap (x_{\rho_2} - \b_\sigma) \subseteq
\rho_1 - \b_\sigma \subseteq
\rho_1 - \b_{\g + 1} = \emptyset$.
Thus, the sequence $\la x_\rho - h_\a(\rho) :
\rho \in S' \cap \a \ra$ is again as desired.
This proves Lemma 2.
 
\no \hfill \finpf Lemma 2
 
At this point,
we are in a position to define in $V_1$ the partial ordering
$P^1_{\delta, \l}
[S] $ which will be used to destroy strong compactness.
$P^1_{\delta, \l} [S]$ is now the set of all $4$-tuples
$\langle w, \a, \bar r, Z \rangle $ satisfying the following properties.
 
\item{1.} $ w \in [\l^+]^{< \l}$.
\item{2.} $\a < \l$.
\item{3.} $\bar r = \langle r_i : i \in w \rangle $ is a sequence of
functions from $\a $ to $\{ 0, 1 \}$, i.e., a sequence of subsets of $\a$.
\item{4.} $Z \subseteq                 \{x_\beta : \beta \in
S \}$ is a set so that if $z \in Z$, then for some
$y \in {[w]}^\delta$, $y \subseteq z$ and $z - y$ is
bounded in the $\b$ so that $z = x_\b$.
 
           \noindent
Note that the definition of $Z$ implies $|Z| < \l$.
 
The ordering on $P^1_{\delta, \l} [S]$ is given by $\langle w^1, \alpha^1,
\bar  r^1,Z^1 \rangle \le \langle w^2, \a^2, \bar r^2, Z^2 \rangle$ iff the
following hold.
\item{1.} $w^1 \subseteq w^2$.
\item{2.} $\alpha^1 \le \a^2$.
\item{3.} If $i \in w^1$, then $r^1_i \subseteq r^2_i$.
\item{4.} $Z^1 \subseteq Z^2.$
\item{5.} If $z \in  Z^1 \cap {[w^1]}^\delta$ and $\a_1 \le \a < \a_2$, then
$| \{ i \in z: r^2_i (\a) = 0 \} | = | \{ i \in z: r^2_i (\a) = 1 \} | =
 \delta$.
 
If $W = \la \la w^\b, \a^\b, \bar r^\b, Z^\b
\ra_{\b < \g < \d} \ra$ is a directed set of elements of
$P^1_{\d, \l}[S]$, then since by the regularity of
$\d$ any $\d$ sequence from $\underset \b < \g \to{\cup}
w^\b$ must contain a $\d$ sequence from $w^\b$ for some
$\b < \g$, it can easily be verified that
$\la \underset \b < \g \to{\cup} w^\b,
     \underset \b < \g \to{\cup} \a^\b,
     \underset \b < \g \to{\cup} \bar r^\b,
     \underset \b < \g \to{\cup} Z^\b \ra$ is an
upper bound for each element of $W$. (Here, if
$\bar r^\b = \la r^\b_i : i \in w^\b \ra$, then
$r_i \in \underset \b < \g \to{\cup} \bar r^\b$ if
$i \in \underset \b < \g \to{\cup} w^\b$ and
$r_i = \underset \b < \g \to{\cup} r^\b_i$, taking
$r^\b_i = \emptyset$ if $i \not\in w^\b$.)
This means $P^1_{\d, \l}[S]$ is $\d$-directed closed.
 
At this point, a few intuitive remarks are in order.
If $\k$ is $\l$ strongly compact for $\l \ge \k$ regular,
 then it must be the case (see [SRK]) that $\l$ carries a $\k$-additive
uniform ultrafilter.
If $\delta < \k < \l$, the forcing $P^1_{\del, \l}  [S]$ has specifically been
 designed to
destroy this fact.
It has been designed, however, to destroy the $\l$ strong compactness
of $\k$ ``as lightly as possible'', making little damage. In the
case of the argument of [KiM], the non-reflecting stationary
set $S$ is added directly to $\l$ in order to kill the $\l$
strong compactness of $\k$.  In our situation, the
non-reflecting stationary set $S$, having been added to
$\l^+$ and not to $\l$, does not kill the $\l$ strong compactness
of $\k$ by itself. The additional forcing $P^1_{\delta,
\l}[S]$ is necessary to do the job. The forcing $P^1_{\delta,
\l}[S]$, however, has been designed so that if necessary,
we can resurrect the $\l$ supercompactness of $\k$ by forcing further with
$P^2_{\delta, \l} [S]$.
 
\proclaim{Lemma 3} $V^{P^1_{\delta, \l}[S]}_1 \models ``\k$ is not $\l$
strongly  compact" if \ $\delta < \k < \l$.
\endproclaim
 
\underbar{Remark}: Since we will only be concerned in general when $\k$
is strongly inaccessible and $\delta < \k < \l$, we assume
without loss of generality that this is the case
throughout the rest of the paper.
 
\demo{Proof of Lemma 3}  Assume to the contrary that
$V^{P^1_{\delta, \l}[S]}_1 \models ``\k$ is $\l$ strongly compact", and by our
earlier remarks, let $p \force `` \dot{\cal D}$ is a $\k$-additive uniform
ultrafilter over $\l$".
We show that $p$ can be extended to a condition $q$ so that for some ordinal
 $\a^q  < \l$ and some $\delta$ sequence $\langle s_i: i < \delta \rangle $
of $\cD$ measure 1 sets, $ q \force `` \underset i < \delta \to{\cap}
\dot s_i \subseteq \a^q"$, an immediate contradiction.

We use a $\Delta$-system argument to establish this.
First, for $G_1$ $V_1$-generic over $P^1_{\delta, \l} [S]$ and $i <
\lambda^+$, let $r^*_i = \cup \{ r^p_i: \exists p = \langle w^p, \a^p, \bar
r^p, Z^p \rangle \in$ $G_1 [r^p_i \in \bar r^p ] \}$.  It is the case that $
\force_{P^1_{\delta, \l}[S]} ``\dot r^*_i: \l \to \{ 0, 1 \}$ is a function
whose domain is all of $\l$".  To see this, for $p = \la w^p, \a^p, \bar r^p,
Z^p \ra$, since $|Z^p| < \l$, $w^p \in {[\l^+]}^{< \l}$, and $z \in Z^p$
implies $z \in {[\l^+]}^\d$, the condition $q = \la w^q, \a^q, \bar r^q, Z^q
\ra$ given by $\a^q = \a^p$, $Z^q = Z^p$, $w^q = w^p \cup \bigcup\{z : z \in
Z^p \}$, and $\bar r^q = \la r'_i : i \in w^q \ra$ defined by $r'_i = r_i$ if
$i \in w^p$ and $r'_i$ is the empty function if $i \in w^q - w^p$ is a
well-defined condition. (This just means we may as well assume that for $p =
\la w^p, \a^p, \bar r^p, Z^p \ra$, $z \in Z^p$ implies $z \subseteq w^p$.)
Further, since $|Z^q| < \l$, $\cup \{\b : \exists z \in Z^q[z = x_\b]\} = \g <
\l^+$. Therefore, if $\g' \in (\g, \l^+)$ and $S' \subseteq \g'$ is so that
$\sup S' = \g'$ and $S'$ is an initial segment of $S$ so that $S'$ is not
stationary at its supremum nor has any initial segment which is stationary at
its supremum, then by Lemma 2, there is a sequence $\la y_\b : \b \in S' \ra$
so that for every $\b \in S'$, $y_\b
\subseteq x_\b$, $x_\b - y_\b$ is bounded in $\b$, and if
$\b_1 \neq \b_2 \in S'$, then $y_{\b_1} \cap y_{\b_2} =
\emptyset$. This means that if $z \in Z^q$ and $z = x_\b$
for some $\b$, then $y_\b \subseteq w$.
 
Choose now for $\b \in S'$ sets $y^1_\b$ and $y^2_\b$ so that $y_\b = y^1_\b
\cup y^2_\b$, $y^1_\b \cap y^2_\b = \emptyset$, and $|y^1_\b| = |y^2_\b| =
\d$. If  $\rho \in (\a^q, \l)$, then for each $\b$ so that $x_\b \in
Z^q$ and for each $r'_i \in \bar r^q$ such that
$i \in y_\b$, we can extend $r'_i$ to $r''_i : \rho \to
\{0,1\}$ by letting $r''_i \vert \a^q = r'_i \vert \a^q$,
and for $\a \in [\a^q, \rho)$, $r''_i(\a) = 0$ if
$i \in y^1_\b$ and $r''_i(\a) = 1$ if $i \in y^2_\b$.
For $i \in w^q$ so that there is no $\b$ with $x_\b \in
Z^q$ and $i \in y_\b$, we extend $r'_i$ to $r''_i : \rho
\to \{0,1\}$ by letting $r''_i \vert \a^q = r'_i \vert
\a^q$, and for $\a \in [\a^q, \rho)$, $r''_i(\a) = 0$.
If we let $\bar s = \la r''_i : i \in w^q \ra$, then
$t = \la w^q, \rho, \bar s, Z^q \ra$ can be verified to be
such that $t$ is well-defined and $t \ge q \ge p$. We have
therefore shown by density that
$\force_{P^1_{\d, \l}[S]} ``\dot r^*_i \to \{0,1\}$ is
a function whose domain is all of $\l$''.
Thus, we can let $ r^\ell_i = \{ \a < \l: r^*_i (\a) = \ell \} $
for $\ell \in \{0,1\}$.
 
For each $i < \l^+$, pick $p_i = \langle w^{p_i}, \a^{p_i}, \bar r^{p_i},
 Z^{p_i} \rangle
 \ge p$ so that $p_i \force ``\dot r^{\ell(i)}_i \in \dot {\cal D} "$
for some $\ell(i) \in \{0,1\}$. This is possible since
$\force_{P^1_{\delta, \lambda}[S]} ``$For each $i < \lambda^+$,
$\dot r^0_i \cup \dot r^1_i = \lambda  $''.
Without loss of generality, by extending $p_i$ if necessary, we can assume that
$i \in w^{p_i}$.
Thus, since each $w^{p_i} \in [\l^+]^{< \l}$, we can find some
stationary $A \subseteq \{i < \lambda^+$ : cof$(i) = \lambda
\}$
so that $\{ w^{p_i} : i \in A \} $ forms a
$\Delta$-system, i.e., so that for $i \neq j \in A$, $w^{p_i} \cap
w\pj$ is some constant value $w$ which is an initial segment
of both. (Note we can assume that for $i \in A$, $w_i \cap i = w$,
and for some fixed $\ell(*) \in \{0,1\}$, for every $i \in A$,
$p_i \force ``\dot r^{\ell(*)}_i \in \dot{\cal D}$''.)
Also, by clause 4) of the definition of the forcing,
$|Z^{p_i}| < \l$ for each $i < \l^+$.
Therefore, $Z^{p_{i}} \in [[\l^+]^\delta]^{< \l}$, so as
$|[\l^+]^\delta| = \l^+ $ by GCH, the same sort of $\Delta$-system argument
 allows
us to assume in addition that for all $i \in A$, $Z^{p_i}
\cap {\cal P}(w)$ is  some constant value $Z$.  Further, since each
$\a{^{p_i}} < \l$, we can assume that $\a{^{p_i}}$ is some constant $\a^0$ for
$i \in A$.  Then, since any $\bar r{^{p_i}} = \langle r_j: j \in w{^{p_i}}
\rangle$ for $i \in A$ is composed of a sequence of functions from $\a_0$ to
$2$, $\a_0 < \l,$ and $|w|< \l$, GCH allows us to conclude that for $i \neq j
\in A$, $\bar r{^{p_i}} \vert w = \bar r\pj \vert w$.
And, since $i \in w{^{p_i}}$, we know that we can also assume (by thinning $A$
if necessary) that $B = \{ \sup (w{^{p_i}}): i \in A \}$ is so that $i < j \in
A$ implies $i \le \sup (w{^{p_i}}) <$ min$(w^{p_j} - w) \le \sup (w\pj)$.  We
know in addition by the choice of $X = \langle x_\beta : \beta
\in S \rangle$ that for some $\gamma \in S$, $x_\gamma
\subseteq A$. Let $x_\gamma = \{i_\beta : \beta < \delta \}$.
 
We are now in a position to define the condition $q$ referred to earlier.
We proceed by defining each of the four coordinates of $q$.
First, let $w^q	 = \underset \b < \delta \to{\cup} w^{p_{i_{\beta}}}$.
As $\l$ and $\l^+$ are regular, $\delta < \l$, and each
$w^{p_{i_\beta}} \in
[\l^+]^{<\l}$, $ w^q$
is well-defined and in ${[\lambda^+]}^{<\lambda}$.
Second, let $\a^q = \a^0$.
Third, let $\bar r^q
= \langle r^q_i : i \in w^q \rangle$ be defined by  $r^q_i =
 r^{p_{i_{\beta}}}_i$
 if $i \in w^{p_{i_{\beta}}}$.
The property of the $\Delta$-system that $i \neq j \in A$ implies $\bar
 r{^{p_i}} \vert w
= \bar r\pj \vert w$ tells us $\bar r^q$ is well defined.
Finally, to define $Z^q$,
let $Z^q = \underset  \b < \delta \to{\cup} Z^{i_{\beta}} \cup
\{ \{ i_\beta : \beta < \delta \} \}$.
By the last three sentences   in the preceding paragraph and our construction,
$ \{ i_\beta: \beta < \delta \}$ generates a new set
which can be included in $Z^q$, and $Z^q$ is well-defined.
 
We claim now that $q \ge p$ is so that $q \force ``\underset \beta <
\delta\to{\cap} \dot r^{\ell(*)}_{i_{\beta}}
\subseteq \a^q".$
To see this, assume the claim fails.  This means that for some $q^1 \ge q$ and
some $\alpha^q \le \eta < \l$, $q^1\force ``\eta \in \underset \beta < \delta
\to{\cap} \dot r^{\ell(*)}_{i_{\beta}}"$. 
Without loss of generality, since $q^1$ can always be extended if necessary,
we can assume that $\eta < \alpha^{q^1}$. But then, by the definition of
$\le$, for $\delta$ many $\beta < \delta$, $q^1 \force `` \eta \notin \dot
r^{\ell(*)}_{i_{\beta}}"$, an immediate contradiction.  Thus, $q \force
``\underset \beta < \delta \to{\cap} \dot r^{\ell(*)}_{i_{\beta}} \subseteq
\alpha^q " $, which, since $\delta < \k$, contradicts that
$q \force `` \underset \beta < \delta \to{\cap} \dot r^{\ell(*)}_{i_{\beta}}
\in \dot {\cal D}$ and  $\dot {\cal D}$ is a $\k$-additive uniform ultrafilter
over $\l$". This proves Lemma 3.
 
\noindent
\hfill $\square $  Lemma 3
 
Recall we mentioned prior to the proof of Lemma 3 that $P^1_{\delta, \l}[S]$
is designed so that a further forcing with $P^2_{\delta, \l} [S]$ will
resurrect the $\l$ supercompactness of $\k$, assuming the correct iteration
has been done.  That this is so will be shown in the next section.  In the
meantime, we give an idea of why this will happen by showing that the forcing
$P^0_{\delta, \lambda} * (P^1_{\delta, \l} [\dot S] \times P^2_{\delta, \l}
[\dot S])$ is rather nice.  Specifically, we have the following lemma.
\proclaim{Lemma 4} $P^0_{\delta, \l} * (P^1_{\delta, \l}[\dot S] \times
P^2_{\delta, \l}  [\dot S])$ is equivalent to $Q^0_\l * \dot Q^1_\l$.
\endproclaim
\demo{Proof of Lemma 4} Let $G$ be $V$-generic over $P^0_{\delta, \l}
* (P^1_{\dell} [\dot S] \times P^2_{\dell} [\dot S])$, with $G ^0_{\dell}$, $G
^1_{\dell}$, and $G ^2_{\dell}$ the projections onto $P^0_{\dell}$,
$P^1_{\dell} [S]$, and $P^2_{\dell} [S]$ respectively.  Each $G ^i_{\dell}$ is
appropriately generic.  So, since $P^1_{\dell} [S] \times P^2_{\dell} [S]$ is
a product in $V[G^0_{\dell}]$, we can rewrite the forcing in $V[G^0_{\dell}]$
as $P^2_{\dell} [S] \times P^1_{\dell} [S]$ and rewrite $V[G]$ as
$V[G^0_{\dell}] [G^2_{\dell}] [G^1_{\dell}]$.
 
It is well-known (see [MS]) that the forcing $P^0_{\dell} * P^2_{\dell} [\dot
S]$ is equivalent to $Q^0_\l$.
That this is so can be seen from the fact that $P^0_{\dell} * P^2_{\dell}
[\dot S]$ is non-trivial, has cardinality
$\l^+$, and is such that $D=\{ \langle p, q
\rangle \in P^0_{\dell} * P^2_{\dell} [\dot S]$ : For some
$\alpha$, $\dom (p) = \dom (q)  = \a + 1$, $p \force `` \a \notin
\dot S"$, and $q \force `` \a \in \dot C "  \}$ is dense
 in $P^0_{\dell} * P^2_{\dell} [\dot S]$ and is $\l$-closed.
This easily implies the desired equivalence.
Thus, $V$ and $V[G^0_{\dell}] [G^2_{\dell}]$ have the same cardinals and
cofinalities, and the proof of Lemma 4 will be complete once we show that in
$V[G^0_{\dell}] [G^2_{\dell}]$, $P^1_{\dell}  [S]$ is equivalent to $Q^1_\l$.
 
To this end, working in $V[G^0_{\dell}][G^2_{\dell}]$, we first note that as
$S \subseteq \l^+$ is now a non-stationary set all of whose initial segments
are non-stationary, by Lemma 2, for the sequence $\langle x_\beta: \beta \in S
\rangle$, there must be a sequence $\langle y_\beta: \beta \in S \rangle $ so
that for every $\beta \in S$, $y_\beta \subseteq x_\beta$, $x_\beta - y_\beta$
is bounded in $\b$, and if $\beta_1 \neq \beta_2 \in S$, then $y_{\beta_{1}}
\cap y_{\beta_{2}} = \emptyset$.  Given this fact, it is easy to observe that
$P^1 = \{ \langle w, \a, \bar r, Z \rangle \in P^1_{\dell} [S]:$ For every
$\beta \in S$, either $y_\beta \subseteq w $ or $ y_\beta \cap w =
\emptyset\}$ is dense in $P^1_{\dell} [S]$.  To show this, given $\langle w,
\a, \bar r, Z \rangle \in P^1_{\dell} [S]$, $ \bar r = \langle r_i : i \in w
\rangle$, let $Y_w = \{ y \in \langle y_\beta: \beta \in S \rangle: y \cap w
\neq \emptyset \}$.  As $|w| < \l$ and $y_{\beta_{1}} \cap y_{\beta_{2}} =
\emptyset $ for $\b_1 \neq \b_2 \in S$, $| Y_w | < \l$.  Hence, as $|y| =
\delta < \l$ for $y \in Y_w$, $|w'|<\l$ for $ w' = w \cup ( \cup Y_w)$.  This
means $\langle w', \a, \bar r', Z \rangle $ for $\bar r'= \langle r'_i: i \in
w' \rangle$ defined by $r'_i = r_i$ if $ i \in w $ and $r'_i$ is the empty
function if $i \in w' - w$ is a well-defined condition extending $\langle w,
\a, \bar r, Z \rangle$.  Thus, $P^1$ is dense in $P^1_{\dell} [S]$, so to
analyze the forcing properties of $P^1_{\dell} [S]$, it suffices to analyze
the forcing properties of $P^1$.
 
For $\b \in S$, let $Q_\b = \{ \langle w, \a, \bar r, Z \rangle \in P^1: w = 
 y_\b \}$, and let $ Q' = \{ \langle w, \a, \bar r, Z \rangle  \in P^1: w
\subseteq \l^+ - \underset \b \in S \to{\cup} y_\beta \}$.
Let $Q''$ be those elements of $\mathop{\dsize\Pi}\limits_{\b \in S} Q_\b
\times Q'$ of support $ < \l$ under the product ordering.  Adopting the
notation of Lemma 3, given $p = \langle \langle q_\b : \b \in A \rangle, q
\rangle \in Q''$ where $ A \subseteq S$ and $ |A| < \l$, as $| A | < \l$ and
$\l $ is regular, $\a = \sup \{ \a^{q_\beta} : \beta \in A \} \cup \a^q < \l$,
so without loss of generality, each $q_\b$ and $q$ can be extended to
conditions $q'_\beta$ and $q'$ so that $\a $ occurs in $q'_\beta$ and $q'$.
This means $Q = \{ p = \langle q_\beta: \beta < \gamma < \l \rangle \in Q'' :
\a^{q_\beta} = \a^{q_{\beta'}}$ for $\b $ and $\b'$ different coordinates of
$p \}$ is dense in $Q''$, so $Q$ and $Q''$ are forcing equivalent.  Then, for
$p = \langle \langle q_\beta : \beta \in A \rangle, q \rangle \in Q$ where $A
\subseteq S$ and $| A | < \l$, as $w^{q_{\b_{1}}} \cap w^{q_{\b_{2}}} =
\emptyset $ for $ \b_1 \neq \b_2 \in A$ $ (y_{\b_1} \cap y_{\b_2} =
\emptyset)$, $w^{q_{\b_{1}}} \cap w^q = \emptyset$, $\a^{q_{\b_{1}}} =
\a^{q_{\b_{2}}} = \a^q$ for $\b_1 \neq \b_2 \in A$, the domains of any two
$\bar r^{q_{\b_{1}}}$, $ \bar r^{q_{\b_{2}}}$ are disjoint for $\b_1 \neq
\b_2 \in A$, $Z^{q_{\b_{1}}} \cap Z^{q_{\b_{2}}} = \emptyset $ for
$\b_1 \neq \b_2 \in A$, the domains of $\bar r^{q_\b} $ and $\bar r^q$ are
disjoint for $\b \in A$, and $Z^{q_\b} \cap Z^q = \emptyset $ for $\b \in A$,
the function $ F(p) = \langle \underset \b \in A \to{\bigcup} w^{q_\beta} \cup
w^q, \a, \underset \b \in A \to{\bigcup} \bar r^{q_\beta} \cup \bar r^q,
\underset \b \in A \to{\bigcup} Z^{q_\b} \cup Z^q \rangle$
can easily be seen to yield an isomorphism between $Q$ and $P^1$.
Thus, over $V[G^0_{\dell}][G^2_{\dell}]$, forcing with $P^1$,
 $P^1_{\dell} [S]$, $Q$, and $Q''$ are all equivalent.
 
We examine now in more detail the exact nature of $Q''$.  For $\b \in S$, GCH
shows $|Q_\b| = \l$.  It quickly follows from the definition of $Q_\b$ that
$Q_\b$ is $< \l$-closed, so $Q_\b$ is forcing equivalent to adding a Cohen
subset to $\l$.  Since the definitions of $P^1_{\delta, \l}[S]$ and $P^1$
ensure that for $ \langle w, \a, \bar r, Z \rangle \in Q'$, $Z = \emptyset$
(for every $\b \in S$, $w \cap y_\b = \emptyset$, $y_\b
\subseteq x_\b$, and $x_\b - y_\b$ is bounded in $\delta$),
$Q'$ can easily be seen to be a re-representation of the
Cohen forcing where instead of working with functions whose domains have
cardinality $< \l$ and are subsets of $\l \times \l^+$, we work with functions
 whose domains have cardinality $< \l$ and are subsets of $\l \times (\l^+
- \underset \b\in S \to{\cup} y_\beta).$
Thus, $Q'' $ is isomorphic to a Cohen forcing using functions having
domains of cardinality $< \l $ which adds $\l^+$ many Cohen subsets to $\l$.
By the last sentence of the last paragraph, this means that over $V
 [G^0_{\dell}] [G^2_{\dell}]$, the forcings $P^1_{\dell} [S] $ and $Q^1_\l $
are equivalent. This proves Lemma 4.
 
\noindent
\hfill $\square $  Lemma 4
 
As we noted in the proof of Lemma 4, without the last coordinate $Z^p$
of a condition $p \in P^1_{\dell}[S]$ and the associated condition on the
ordering, $P^1_{\dell} [S]$ is just a re-representation of $Q^1_\l$.
This last coordinate and change in the ordering are necessary to destroy the
 $\l$ strong compactness of $\k$ when forcing with $P^1_{\dell} [S]$.
Once the fact $S$ is stationary has been destroyed by
forcing with $P^2_{\dell} [S]$, Lemma 4 shows that this last coordinate
$Z^p$ of a condition $p \in P^1_{\dell} [S]$ and change in the ordering
in a sense become irrelevant.
 
It is clear from Lemma 4 that $P^0_{\dell} * (P^1_{\dell} [\dot S]
 \times P^2_{\dell} [\dot S])$, being equivalent to $Q^0_\l  *  \dot Q^1_\l$,
preserves GCH, cardinals, and cofinalities, and has a dense
subset which is $< \l$-closed and satisfies $\l^{++}$-c.c.
Our next lemma shows that the forcing $P^0_{\dell} * P^1_{\dell} [\dot S]$ is
also rather nice.
\proclaim{Lemma 5} $P^0_{\dell} * P^1_{\dell} [\dot S]$ preserves GCH,
cardinals, and cofinalities, is $<\l$-strategically closed, and is
$\l^{++}$-c.c. 
\endproclaim
\demo{Proof of Lemma 5} Let $G' = G^0_{\dell} * G^1_{\dell}$ be $V$-generic
over $P^0_{\dell} * P^1_{\dell}[\dot S]$, and let $G ^2_{\dell}$ be
$V[G']$-generic over $P^2_{\dell} [S]$.  Thus, $G ' * G ^2_{\dell} = G$ is
$V$-generic over $P^0_{\dell} * (P^1_{\dell} [\dot S] * P^2_{\dell} [\dot S])
= P^0_{\dell} * (P^1_{\dell} [\dot S] \times P^2_{\dell} [\dot S])$.  By Lemma
4, $V[G] \models $ GCH and has the same cardinals and cofinalities as $V$, so
since $V[G'] \subseteq V[G]$, forcing with $P^0_{\dell} * P^1_{\dell} [\dot
S]$ over $V$ preserves GCH, cardinals, and cofinalities.
 
We next show the $<\l$-strategic closure of
$P^0_{\delta, \l} * P^1_{\delta, \l}[\dot S]$. We first note that
as $(P^0_{\delta, \l} * P^1_{\delta, \l}[\dot S]) *
P^2_{\delta, \l}[\dot S] = P^0_{\delta, \l} *
(P^1_{\delta, \l}[\dot S] * P^2_{\delta, \l}[\dot S])$ has by
Lemma 4 a dense subset which is $<\l$-closed, the desired fact
follows from the more general fact that if $P * \dot Q$ is a
partial ordering with a dense subset $R$ so that $R$ is
$<\l$-closed, then $P$ is $<\l$-strategically closed. To show
this more general fact, let $\gamma < \l$ be a cardinal.
Suppose I and II play to build an increasing chain of elements of
$P$, with $\langle p_\beta : \beta \le \alpha + 1 \rangle$
enumerating all plays by I and II through an odd stage
$\alpha + 1$ and $\langle \dot{q_\beta} : \beta < \alpha +
1$ and $\beta$ is even or a limit ordinal$\rangle$ enumerating a
set of auxiliary plays by II which have been chosen so that
$\langle \langle p_\beta, \dot{q_\beta} \rangle : \beta <
\alpha + 1$ and $\beta$ is even or a limit ordinal$\rangle$
enumerates an increasing chain of elements of the dense
subset $R \subseteq  P * \dot Q$. At stage $\alpha + 2$, II
chooses $\langle p_{\alpha +2}, \dot q_{\alpha + 2}
\rangle$ so that $\langle p_{\alpha + 2}, \dot q_{\alpha+2}
\rangle \in R$ and so that $\langle p_{\alpha + 2},
\dot q_{\alpha+2} \rangle \ge \langle p_{\alpha + 1},
\dot{q_\alpha} \rangle$; this makes sense, since
inductively, $\langle p_\alpha, \dot q_\alpha \rangle
\in R \subseteq P * \dot Q$, so as I has chosen
$p_{\alpha + 1} \ge p_\alpha$, $\langle p_{\alpha + 1},
\dot q_\alpha \rangle \in P * \dot Q$. By the $<\l$-closure of
$R$, at any limit stage $\eta \le \gamma$, II can choose
$\langle p_\eta, \dot{q_\eta} \rangle$ so that
$\langle p_\eta, \dot{q_\eta} \rangle$ is an upper bound to
$\langle \langle p_\beta, \dot{q_\beta} \rangle : \beta <
\eta$ and $\beta$ is even or a limit ordinal$\rangle$. The
preceding yields a winning strategy for II, so $P$ is
$<\l$-strategically closed.
 
Finally, to show $P^0_{\dell} * P^1_{\dell} [\dot S]$ is $\l^{++}$-c.c.,
 we simply note that
this follows from the general fact about iterated forcing (see
[Ba]) that if $P * \dot Q$ satisfies $\l^{++}$-c.c., then $P$
satisfies $\l^{++}$-c.c. (Here, $P = P^0_{\delta, \l} *
P^1_{\delta, \l}[\dot S]$ and $Q = P^2_{\delta, \l}
[\dot S]$.)
This proves Lemma 5.
 
\noindent
\hfill $\square $  Lemma 5
 
We remark that
$\force_{P^0_{\dell}} ``P^1_{\dell} [\dot S]$ is $\l^+$-c.c.",
for if ${\cal A} = \langle p_\a : \a < \l^+ \rangle$ were a size
$\l^+$ antichain of elements of $P^1_{\delta, \l}[S]$ in
$V[G^0_{\delta, \l}]$, then as $V[G^0_{\delta, \l}]$ and
$V[G^0_{\delta, \l}][G^2_{\delta, \l}]$ have the same cardinals,
${\cal A}$ would be a size $\l^+$ antichain of elements of
$P^1_{\delta, \l}[S]$ in $V[G^0_{\delta, \l}][G^2_{\delta, \l}]$.
By Lemma 4, in this model, a dense subset of $P^1_{\delta, \l}[S]$
is isomorphic to $Q^1_\l$, which has the same definition in
either $V[G^0_{\delta, \l}]$ or $V[G^0_{\delta, \l}]
[G^2_{\delta, \l}]$ (since $P^0_{\delta, \l}$ is $\l$-strategically
closed and $P^0_{\delta, \l} \ast P^2_{\delta, \l}[\dot S]$ is
$\l$-closed) and so is $\l^+$-c.c. in either model.
 
We conclude this section with a lemma which will be used later
in showing that it is possible to extend certain elementary
embeddings witnessing the appropriate degree of
supercompactness.
 
\proclaim{Lemma 6}
For $V_1 = V^{P^0_{\d, \l}}$, the models
$V^{P^1_{\d, \l}[S] \times P^2_{\d, \l}[S]}_1$ and
$V^{P^1_{\d, \l}[S]}_1$ contain the same $\l$ sequences
of elements of $V_1$.
\endproclaim
 
\demo{Proof of Lemma 6}
By Lemma 4, since
$P^0_{\d, \l} \ast P^2_{\d, \l}[\dot S]$
is equivalent to the forcing $Q^0_\l$ and
$V \subseteq
V^{P^0_{\d, \l}}
\subseteq
V^{P^0_{\d, \l} \ast P^2_{\d, \l}[\dot S]}
$, the models $V$, $
V^{P^0_{\d, \l}}
$, and $
V^{P^0_{\d, \l} \ast P^2_{\d, \l}[\dot S]}
$ all contain the same $\l$ sequences of elements of $V$.
Thus, since a $\l$ sequence of elements of $V_1 =
V^{P^0_{\d, \l}}
$ can be represented by a $V$-term which is
actually a function $h : \l \to V$, it immediately follows that
$
V^{P^0_{\d, \l}}
$ and $
V^{P^0_{\d, \l} \ast P^2_{\d, \l}[\dot S]}$
contain the same $\l$ sequences of elements of $
V^{P^0_{\d, \l}}
$.
 
Let now $f : \l \to V_1$ be so that $f \in (
V^{P^0_{\d, \l} \ast P^2_{\d, \l}[\dot S]}
){}^{P^1_{\d, \l}[S]} =
 V^{P^1_{\d, \l}[S] \times P^2_{\d, \l}[S]}_1$, and let
$g : \l \to V_1$, $g \in
V^{P^0_{\d, \l} \ast P^2_{\d, \l}[\dot S]}
$ be a term for $f$. By the previous paragraph, $g \in
V^{P^0_{\d, \l}}
$. Since Lemma 4 shows that $P^1_{\d, \l}[S]$ is
$\l^+$-c.c$.$ in $
V^{P^0_{\d, \l} \ast P^2_{\d, \l}[\dot S]}
$, for each $\a < \l$, the antichain $\A_\a$ defined in $
V^{P^0_{\d, \l} \ast P^2_{\d, \l}[\dot S]}
$ by $\{p \in P^1_{\d, \l}[S] : p$ decides a value for
$g(\a) \}$ is so that $
V^{P^0_{\d, \l} \ast P^2_{\d, \l}[\dot S]}
\models ``|\A_\a| \le \l$''. Hence, by the preceding
paragraph, since $\A_\a$ is a set of elements of $
V^{P^0_{\d, \l}}
$, $\A_\a \in
V^{P^0_{\d, \l}}
$ for each $\a < \l$. Therefore, again by the preceding
paragraph, the sequence $\la \A_\a : \a < \l \ra \in
V^{P^0_{\d, \l}}
$. This just means that the term $g \in
V^{P^0_{\d, \l}}
$ can be evaluated in $
V^{P^1_{\d, \l}[S]}_1
$, i.e., $f \in
V^{P^1_{\d, \l}[S]}_1
$. This proves Lemma 6.
 
\no \hfill $\square$ Lemma 6
 
\S 2 The Case of One Supercompact Cardinal with no Larger Inaccessibles
 
In this section, we give a proof of our Theorem, starting from a model
$V$ for ``ZFC + GCH +
There is one supercompact cardinal $\k$ and no $\l > \k$ is inaccessible".
Before defining the forcing conditions used in the proof of
this version of our Theorem, we first give a proof of
 the theorem of Magidor
mentioned in Section 0 which shows that if there is a supercompact cardinal,
 then there
 always must be cardinals $\delta < \lambda$ so that $\delta$ is
$\l$ strongly compact yet $\delta $ isn't $\l$ supercompact.
\proclaim{Lemma 7} (Magidor [Ma4]): Suppose $\k$ is a supercompact cardinal.
Then $B = \{  \delta < \k: \delta $ is $\l_\delta$ strongly compact for
$\l_\delta$ the least singular strong limit cardinal $> \delta $ of
cofinality $\delta$, $\delta $ is not $\l_\delta$ supercompact, yet $\delta $
is $\a$ supercompact for all $\a < \l_\delta \}$ is unbounded in $\k$.
\endproclaim
 
\demo{Proof of Lemma 7}  Let $\l_\k > \k$ be the least singular strong limit
 cardinal of cofinality $\k$, and let $j: V \to M$ be an elementary embedding
witnessing the $\l_\k$ supercompactness of $\k$ with $j (\k)$ minimal.  As
$j(\k)$ is least, $M \models ``\k$ is not $\l_\k$ supercompact".  As
$M^{\l_{\k}} \subseteq M$ and $\l_\k$ is a strong limit cardinal, $M \models
``\k $ is $\a$ supercompact for all $\a < \k"$.
 
Let $\mu \in V $ be a $\k$-additive measure over $\k$, and let $\langle \l_\a: 
\a < \l_\k \rangle$ be a sequence of cardinals cofinal in $\l_\k$ in both $V$
and $M$. As $M^{\l_{\k}} \subseteq M$ and $\l_\k$ is a strong limit cardinal,
$\mu \in M$.  Also, as $M \models ``\k$ is $\a $ supercompact for all $\a <
\l_\k$", the closure properties of $M$ allow us to find a sequence $\langle
\mu_\a : \a < \k \rangle  \in M$ so that $M \models `` \mu_\a$ is a fine,
normal, $\kappa$-additive ultrafilter over  $P_\k (\l_\a)"$.
Thus, we can define in $M$ the collection $\mu^* $ of subsets of $P_\k
(\l_\k)$ by $A \in \mu^*$ iff $\{ \a < \k$: $ A \vert \l_\a \in \mu_\a \} \in
 \mu$,  where for $A \subseteq P_\k (\l_\k)$, $ A \vert \l_\a = \{ p \cap
P_\k(\l_\a) : p \in A  \}$.
It is easily checked that $\mu^* $ defines in $M$ a $\k$-additive
 fine ultrafilter over $P_\k (\l_\k)$.
Thus, $M \models ``\k$ is $\a$ supercompact for all $\a < \l_\k$, $ \k$ is not
$\l_\k$ supercompact, yet $\k$ is $\l_\k$ strongly compact", so by reflection,
 the set $B$ of the hypothesis is unbounded in $\k$.
This proves Lemma 7.
\pbf
\hfill $\square $ Lemma 7
 
We note that the proof of Lemma 7 goes through if $\l_\delta$ becomes the
least singular strong limit cardinal $> \delta $ of cofinality $\delta^+$, of
cofinality $\delta^{++}$, etc.  To see this, observe that the closure
properties of $M$ and the strong compactness of $\k$ ensure that $\k^+$,
$\k^{++}$, etc. each carry $\k$-additive measures $\mu_{\k^+}$,
$\mu_{\k^{++}}$, etc.  which are elements of $M$.  These measures may then be
used in place of the $\mu$ of Lemma 7 to define the strongly compact measure
$\mu^*$ over $P_\k (\l_\k)$.
 
We return now to the proof of our Theorem. Let ${\bar \delta}
= \langle \delta_\a : \a \le \k \rangle$ enumerate the
inaccessibles $\le \k$, with $\delta_\k = \k$. Note that since
we are in the simple case in which $\k$ is the only supercompact
cardinal in the universe and has no inaccessibles above it,
we can assume each $\delta_\a$ isn't $\delta_{\a + 1}$
supercompact and for the least regular cardinal $\l_\a \ge
\delta_\a$ so that $V \models ``\delta_\a$ isn't $\l_\a$
supercompact'', $\l_\a < \delta_{\a + 1}$. (If $\delta$ were
the least cardinal so that $\delta$ is $<\b$ supercompact for
$\b$ the least inaccessible $> \delta$ yet $\delta$ isn't
$\b$ supercompact, then $V_\b$ would provide the desired model.)
 
We are now in a position to define the partial ordering $P$
used in the proof of the Theorem.
We define a $\k$ stage Easton support iteration $P_\k = \langle \langle P_\a,
\dot Q_\a \rangle : \a < \k \rangle$, and then define $P = P_{\k + 1} = P_\k *
\dot Q_\k
$ for a certain class partial ordering $Q_\k$ definable in $V^{P_\k}$.
The definition is as follows:
\item{1.} $P_0$ is trivial.
\item{2.} Assuming $P_\a$ has been defined for $\a < \k$,
$P_{\a + 1} = P_\a \ast {\dot Q}_\a$, with
$\dot Q_\a$ a term for the full support iteration
$\langle P^0_{\omega, \l} * (P^1_{\omega, \l} [\dot S_\l] \times P^2_{\omega,
\l} [ \dot S_\l]):
\delta^+_\a \le \l < \l_\a$ and $\l$ is regular$\rangle * \langle P^0_{\omega,
 \l_{\a}} * P^1_{\omega, \l_\a} [\dot S_{\l_{\a}}] \rangle$, where $\dot S_\l$ 
 is a term for the non-reflecting stationary subset of $\l^+$ introduced by
$P^0_{\omega, \l}$ for $\l < \l_{\a}$ and $\dot S_{\l_{\a}}$ 
is a term for the non-reflecting stationary subset of $\l^+_\a$ introduced by
 $P^0_{\omega, \l_{\a}}$.
\item{3.} $\dot Q_\k$ is a term for the Easton support iteration of
$\langle P^0_{\omega, \l} \ast (P^1_{\omega, \l} [\dot S_\l] \times
P^2_{\omega, \l} [\dot S_\l]) : \l > \k$ is a regular cardinal$\rangle$, where
as before, $\dot S_\l$ is a term for the non-reflecting stationary subset of
$\l^+$ introduced by $P^0_{\omega, \l}$.
 
 The intuitive motivation behind the above definition is that below $\k$ at
any inaccessible, we must first destroy and then resurrect all ``good"
instances of strong compactness, i.e., those which also witness
supercompactness, but then destroy the least regular ``bad" instance of strong
compactness, thus destroying all ``bad" instances of strong compactness beyond
the least ``bad" instance.  Since $\k$ is supercompact, it has no ``bad"
instances of strong compactness, so all instances of $\k$'s supercompactness
are destroyed and then resurrected.
\proclaim{Lemma 8} For $G$ a $V$-generic class over $P$, $V$ and $V[G]$ have
the same cardinals and cofinalities, and $V[G] \models $ ZFC + GCH.
\endproclaim
\demo{Proof of Lemma 8} Write $G = G_\k\ast H$, where $G_\k$ is $V$-generic
over $P_\k$, and $H$ is a $V[G_\k]$-generic class over $Q_\k$.
We show $V[G_\k][H] \models $ ZFC, and by assuming for the time being that
$V[G_\k] \models$  GCH  and has the same cardinals and cofinalities as $V$,
we show $V[G_\k] [H] \models$  GCH  and has the same cardinals and
cofinalities as $V[G_\k]$ (and hence as $V$). 
 
To do this, note that $Q_\k$ is equivalent in $V[G_\k] = V_1$ to the Easton
support iteration of $\langle Q^0_{\l} \ast \dot Q^1_\l: \l > \k$ is a regular
cardinal$\rangle$, so we assume without loss of generality that $Q_\k$ is in
fact this ordering.  Note also that as we are assuming $\k$ has no
inaccessibles above it, $Q_\k$ is in fact equivalent to the Easton support
iteration of $\langle Q^0_{\l} \ast \dot Q^1_{\l} : \l > \k$ is a successor
cardinal$\rangle$.  We first show inductively that for any successor cardinal
$\delta^+ > \k$, forcing over $V_1$ with the iteration of $\langle Q^0_\l \ast
\dot Q^1_\l: \k < \l < \delta^+$ and $\l$ is a successor cardinal$\rangle $
preserves cardinals, cofinalities, and GCH.  If $\delta $ is regular (meaning
$\delta$ is a successor cardinal since $\k$ has no inaccessibles above it),
then this iteration can be written as $Q_{<_{\delta}} \ast (\dot Q^0_{\delta}
\ast \dot Q^1_\delta)$, where $Q_{<_{\delta}}$ is the iteration of $\langle
Q^0_\l \ast \dot Q^1_\l : \k < \l < \delta$ and $\l$ is a successor cardinal$
\rangle $.  By induction, forcing over $V_1$ with $Q_{<_{\delta}} $ preserves
cardinals, cofinalities, and GCH, so since forcing over $V^{Q_{<_{\delta}}}_1$
with $\dot Q^0_\delta \ast \dot Q^1_{\delta}$ will preserve GCH and the
cardinals and cofinalities of $V^{Q_{<_{\delta}}}_1$, forcing over $V_1$ with
$Q_{<_{\delta}}\ast (\dot Q^0_\delta \ast \dot Q^1_\delta)$ preserves
cardinals, cofinalities, and GCH.
If $\delta$ is singular, let $\gamma < \delta$ be a cardinal in $V_1$, and
write the iteration of $\langle Q^0_{\l} \ast \dot Q^1_\l: \k < \l < \delta^+$
and  $\l$ is a successor cardinal$\rangle $ as $Q_{< \gamma^+} \ast \dot
Q^{\ge \gamma^+}$, where $Q_{< \gamma^+}$  is as above and $\dot Q^{\ge
\gamma^+} $ is a term in $V_1 $ for the rest of the iteration; if $\gamma <
\k$, then $Q_{< \gamma^+}$ is trivial and $\dot Q^{\ge \gamma^+}$ is a term
for the whole iteration. 
By induction, $V^{Q_{<\gamma^+}}_1 \models ``\gamma$ is a cardinal, $2^\gamma
= \gamma^+$, and $\cof (\gamma) = \cof^{V_1} (\gamma)"$, so as
$V^{Q_{< \gamma^+}}_1 \models  ``Q^{\ge \gamma^+}$ is $\gamma$-closed",
$V_1^{ Q_{< \gamma^+} \ast \dot Q^{\ge \gamma^+} }
 \models ``\gamma$ is a cardinal, $2^\gamma = \gamma^+$, and $\cof (\gamma)
 = \cof^{V_1} (\gamma)"$ , i.e., GCH, cardinals, and cofinalities
below $\delta$ are preserved when forcing over $V_1$ with
$Q_{< \gamma^+ } \ast \dot Q^{\ge \gamma^+}$.
In addition, since the last sentence shows any $f: \gamma \to \delta$ or $f:
\gamma \to \delta^+$, $f \in V^{Q_{< \gamma^+} \ast \dot Q^{\ge \gamma^+}}$
is so that $f \in V^{Q_{< \gamma^+}}_1$ for arbitrary $\gamma < \delta$,
the fact $V^{Q_{< \gamma^+}}_1$ and $V_1$ have the same cardinals and
cofinalities, together with the fact $V_1^{Q_{< \gamma^+} \ast \dot Q^{\ge
\gamma^+}} \models `` \delta$ is a singular limit of cardinals satisfying GCH"
yield that forcing over $V_1$  with $Q_{< \gamma^+} \ast \dot Q^{\ge \gamma^+}
$ preserves $\delta$ is a  singular cardinal of the same cofinality as in
$V_1$, $2^\delta = \delta^+$, and $\delta^+$ is a regular cardinal.
Finally, as GCH in $V_1$ tells  us $|Q_{< \gamma^+} \ast \dot Q^{\ge \gamma^+} 
| = \delta^+$, forcing with $Q_{< \gamma^+} \ast \dot Q^{\ge \gamma^+}$ over
$V_1$ preserves cardinals and cofinalities $\ge \delta^{++}$ and GCH $ \ge
\delta^+$. 
 
It is now easy to show $V_2 = V[G_\k] [H] \models$ ZFC + GCH and has the same
cardinals and cofinalities as $V[G_\k] = V_1$.  To show $V_2 \models $ GCH and
has the same cardinals and cofinalities as $V_1$, let again $\gamma$ be a
cardinal in $V_1$, and write $Q_\k = Q_{< \gamma^+} \ast \dot Q$, where $\dot
Q$ is a term in $V_1$ for the rest of $Q_\k$.  As before, $V^{Q_{<
\gamma^+}}_1 \models `` 2^\gamma = \gamma^+$ and $\cof (\gamma) = \cof^{V_1}
(\gamma)"$, so since $V_1^{Q_{< \gamma^+}} \models ``Q$ is $\gamma$-closed",
$V_2 \models ``2^\gamma = \gamma^+ $ and $\cof (\gamma) = \cof^{V_1}
(\gamma)$", i.e., by the arbitrariness of $\gamma$, $ V_2 \models $ GCH, and
all cardinals of $V_1$ are cardinals of the same cofinality in $V_2$.
Finally, as all functions $f: \gamma \to \delta$, $\delta \in V_1$ some
ordinal, $f \in V_2 $ are so that $f \in V^{Q_{< \gamma^+}}_1$ by the last
sentence, it is the case $V_2 \models $ Power Set, and since $V_2 \models AC $
and $Q_\k$ is an Easton support iteration, by the usual arguments, the
aforementioned fact implies $V_2 \models$ Replacement.  Thus, $V_2 \models $
ZFC.
 
It remains to show that $V[G_\k] \models $ GCH and has the same cardinals and
cofinalities as $V$.  To do this, we first note that Easton support iterations
of $\delta$-strategically closed partial orderings are $\delta$-strategically
closed for $\delta$ any regular cardinal.  The proof is via induction.  If
$R_1$ is $\delta$-strategically closed and $\force_{R_1} ``\dot R_2$ is
$\delta$-strategically closed", then let $ p \in R_1$ be so that $p \force
``\dot g$ is a strategy for player II ensuring that the game which produces an
increasing chain of elements of $\dot R_2$ of length $\delta $ can always be
continued for $\a \le \delta"$.  If II begins by picking $r_0 = \langle p_0,
\dot q_0 \rangle \in R_1 \ast
\dot R_2$ so that $p_0 \ge p$ has been chosen
according to the strategy $f$ for $R_1$ and $p_0 \force
``\dot q_0$ has been chosen according to $\dot g"$, and at even stages $\a + 2$
picks $r_{\a + 2} = \langle p_{\a + 2}, \dot q_{\a+ 2} \rangle$ so that
$p_{\a + 2}$ has been chosen according to $f$ and is so that $p_{\a + 2}
\force ``\dot q_{\a + 2}$ has been chosen according to $\dot g$",
 then at limit stages $\l \le \delta$, the chain
$r_0 = \langle p_0, \dot q_0 \rangle \le r_1 = \langle p_1, \dot q_1 \rangle
\le \cdots \le r_\a = \langle p_\a, \dot q_\a \rangle \le \cdots (\a < \l)$ is
so that II can find an upper bound $p_\l$ for
$\langle p_\a: \a < \l \rangle $ using $f$.
By construction, $p_\l \force `` \langle \dot q_\a : \a < \l \rangle  $ is so
 that at
 limit and even stages, II has played according to $\dot g$",  so for some
$\dot q_\l$, $p_\l \force `` \dot q_\l$ is an upper bound to $\langle \dot
q_\a:  \a < \l \rangle "$,  meaning the condition $\langle p_\l, \dot q_\l
\rangle$ is as  desired.
These methods, together with the usual proof at limit stages (see [Ba],
Theorem 2.5) that the Easton support iteration of $\delta$-closed partial
orderings is $\delta$-closed, yield that $\delta$-strategic closure is
preserved at limit stages of all of our Easton support iterations of
$\delta$-strategically closed partial orderings.  In addition, the ideas of
this paragraph will also show that Easton support iterations of $\prec
\d^+$-strategically closed partial orderings are $\prec \d^+$-strategically
closed for $\d$ any regular cardinal.
 
For $\a < \k$ and $P_{\a + 1} = P_\a \ast \dot Q_\a$, since
$\l_\a < \delta_{\a + 1}$,
the definition of $Q_\a$ in $V^{P_\a}$ implies $V^{P_\a} \models ``|Q_\a| <
 \delta_{
\a + 1}"$.
This fact, together with Lemma 5 and the definition of $Q_\a$ in $V^{P_\a}$,
now yield the proof that $V^{P_{\a + 1}} \models $ GCH and has the same
cardinals and cofinalities as $V$ is virtually identical to the proof given in
the first part of this lemma that $V_2 \models $ GCH and has the same
cardinals and cofinalities as $V_1$, replacing $\gamma$-closure with
$\gamma$-strategic closure, which also implies that the forcing adds no new
functions from $\gamma$ to the ground model.
 
If $\l$ is a limit ordinal so that $\bar \l = \sup (\{ \delta_\a : \a < \l \}
)$ is singular, then again, the proof that $V^{P_{\l}} \models $ GCH and has
the same cardinals and cofinalities as $V$ is virtually the same as the just
referred to proof of the first part of this lemma for virtually identical
reasons as in the previous sentence, keeping in mind that since $|P_\a| <
\delta_\a$ inductively for $\a < \l$, $| P_\l| = \bar \l^+$.  If $\l \le \k$
is a limit ordinal so that $\bar \l = \l$, then for cardinals $\gamma \le \l$,
the proof that $V^{P_\l} \models ``\gamma $ is a cardinal and $\cof (\gamma) =
\cof^V (\gamma)$"  is once more as before, as is the proof that $V^{P_\l}
 \models ``2^\gamma = \gamma^+"$ for $\gamma < \l$.
As again $|P_a|< \delta_\a < \l$ for $\a < \l$,
$|P_\l| = \l$, so $V^{P_\l} \models ``\gamma$ is a cardinal, $\cof (\gamma) =
\cof^V (\gamma)$, and $2^\gamma = \gamma^+"$ for $\gamma \ge \l$ a cardinal.
Thus, $V[G_\k] \models $  GCH  and has the same cardinals and cofinalities
 as $V.$ This proves Lemma 8.
\pbf \hfill $\square $ Lemma 8
 
We now show that the intuitive motivation for the definition of $P$ as set
forth in the paragraph immediately preceding the statement of Lemma 8 actually
works.
\proclaim{Lemma 9} If $\delta < \gamma $ and $V \models ``\delta $ is $\gamma$
 supercompact
and $\gamma $ is regular'', then for
$G$ $V$-generic
 over $P$, $V[G] \models ``\delta$ is $\gamma $ supercompact".
\endproclaim
\demo{Proof of Lemma 9} Let
$j: V \to M$ be an elementary embedding witnessing the $\gamma $
 supercompactness
 of $\delta $ so that $M \models ``\delta$ is not $\gamma$ supercompact".
For the $\alpha_0$ so that $\delta = \delta_{\a_0}$, let $P =
P_{\a_0} \ast \dot Q'_{\a_0} \ast \dot T_{\a_0} \ast \dot R$, where $\dot
 Q'_{\a_0}$
 is a term for the full support iteration of
$\langle P^0_{\omega, \l} \ast (P^1_{\omega, \l} [ \dot S_\l] \times
P^2_{\omega, \l} [\dot S_\l]): \delta^+ \le       \l \le \gamma$ and $\l$ is
regular$\rangle$, 
$\dot T_{\a_0}$ is a term for the rest of $Q_{\a_0}$, and $\dot R$ is a term
for the rest of $P$.
We show that $V^{P_{\a_{0}}  \ast \dot Q'_{\a_0}}\models ``\delta $ is $\gamma
 $ supercompact".
This will suffice, since $\force_{P_{\a_0} \ast \dot Q'_\a} ``\dot T_{\a_0}
\ast \dot R$ is $\gamma$-strategically closed", so as the regularity of
$\gamma$ and GCH in $V^{ P_{\a_{0}} \ast \dot Q'_{\a_0}}$ imply $V^{ P_{\a_0}
\ast \dot Q'_{\a_{0}}} \models ``| [ \gamma]^{< \delta}| = \gamma"$, if
$V^{ P_{\a_{0}} \ast \dot Q'_{\a_{0}}} \models `` \delta $ is $\gamma $
supercompact'', then $V^{ P_{\a_0} \ast \dot Q'_{\a_0} \ast \dot T_{\a_0} \ast
\dot R} = V^P \models `` \delta $ is $\gamma $ supercompact via any
ultrafilter ${\cal U} \in V^{P_{\a_0} \ast \dot Q'_{\a_0}}$''.
 
To this end, we first note we will actually show that for $G_{\a_0} \ast
G'_{\a_0}$ the portion of $G$ $V$-generic over $P_{\a_0} \ast \dot Q'_{\a_0}$,
the embedding $j$ extends to $ k : V [G_{\a_0} \ast G '_{\a_0}] \to M [H]$ for
some $ H \subseteq j(P)$.  As $\langle j (\a): \a < \gamma \rangle \in M$,
this will be enough to allow the definition of the ultrafilter $x \in {\cal
U}$ iff $\langle j (\a ) : \a < \gamma \rangle \in  k(x) $ to be given in
$V[G_{\a_0} \ast G'_{\a_0}]$.
 
We construct $H$ in stages.  In $M$, as $\d = \d_{\a_0}$ is the critical point
of $j$, $j(P_{\a_0} \ast \dot Q'_{\a_0}) = P_{\a_0} \ast \dot R'_{\a_0} \ast
\dot R''_{\a_0} \ast \dot R'''_{\a_0}$, where $\dot R'_{\a_0} $ will be a term
for the full support iteration of $\langle P^0_{\omega, \lambda} \ast
(P^1_{\omega, \l} [\dot S_\l] \times P^2_{\omega, \l} [\dot S_\l]) :
\delta^+ \le \l < \gamma $ and $\l$ is regular$\rangle \ast \langle
P^0_{\omega, \gamma} \ast P^1_{\omega, \gamma} [\dot S_\gamma] \rangle $
(note that as $M^\gamma \subseteq M$, GCH implies that $ M \models ``\delta $
is $\l$ supercompact" if $\l < \gamma $ is regular, so since $M \models ``
\delta$ is not $\gamma$ supercompact", $\dot R'_{\a_0}$ is indeed as just
stated), $ \dot R''_{\a_0}$ will be a term for the rest of the portion of
$j(P_{\a_0})$ defined below $j(\delta)$, and $\dot R'''\ao$ will be a term for
$j (\dot Q'_{\a_0})$.  This will allow us to define $ H$ as $H _{\a_0} \ast H
'_{\a_0} \ast H ''_{\a_0}\ast  H'''_{\a_0}$.
Factoring $G'_{\a_0} $ as $\langle G^0_{\omega, \l} \ast (G^1_{\omega, \l}
\times G ^2_{\omega, \l}): \delta^+ \le \l \le \gamma$ and $\l$ is
regular$\rangle$, we let $  H  _{\a_0} =    G _{\a_0}$ and
$  H  '_{\a_0} = \langle G^0_{ \omega, \l} \ast (    G ^1_{\omega, \l} \times
G ^2_{\omega, \l} ) :  \delta^+
 \le \l < \gamma $ and $\l $ is regular$ \rangle *
\langle G^0_{\omega, \gamma} *    G ^1_{\omega, \gamma} \rangle$.
Thus, $  H  '_{\a_0}$ is the same as $G    '_{\a_0}$, except, since
$M \models ``\delta$ is not $\gamma$ supercompact", we omit the generic object
$G    ^2_{\omega, \gamma}$.
 
To construct $H_{\alpha_0}''$, we first note that the
definition of $P$ ensures $|P_{\alpha_0}| = \delta$ and, since
$\delta$ is necessarily Mahlo, $P_{\alpha_0}$ is
$\delta$-c.c. As $V[G_{\alpha_0}]$ and $M[G_{\alpha_0}]$ are
both models of GCH, the definition of $R_{\alpha_0}'$ in
$M[H_{\alpha_0}]$, Lemmas 4, 5, and 8, and the remark
immediately following Lemma 5 then ensure that
$M[H_{\alpha_0}] \models ``$The portion of
$R_{\a_0}'$ below $\g$ is $\g^+$-c.c$.$ and the portion of
$R_{\a_0}'$ at $\g$ is a $\g$-strategically closed partial
ordering followed by a $\g^+$-c.c$.$ partial ordering''.
Since $M^\gamma \subseteq M$
implies ${(\gamma^+)}^V =
{(\gamma^+)}^M$ and
$P_{\alpha_0}$ is $\delta$-c.c., Lemma 6.4 of [Ba] shows
$V[G_{\alpha_0}]$ satisfies these facts as well.
This means applying the argument of Lemma 6.4 of
[Ba] twice, in concert with an application of the fact
a portion of $R_{\a_0}'$ at $\g$ is $\g$-strategically
closed, shows
$M[H_{\alpha_0} *
H_{\alpha_0}'] = M[G_{\alpha_0} * H_{\alpha_0}']$ is
closed under $\gamma$ sequences with respect to
$V[G_{\alpha_0} * H_{\alpha_0}']$, i.e., if $f:\gamma \to
M[H_{\alpha_0} * H_{\alpha_0}']$, $f \in V[G_{\alpha_0} *
H_{\alpha_0}']$, then $f \in M[H_{\alpha_0} *
H_{\alpha_0}']$. Therefore, as $M[H_{\alpha_0} *
H_{\alpha_0}'] \models ``R_{\alpha_0}''$ is both
$\gamma$-strategically closed and $\prec
\gamma^+$-strategically closed'', these facts are true in
$V[G_{\alpha_0} * H_{\alpha_0}']$ as well.
 
Observe now that GCH allows us to assume $\gamma^+ < j (\delta) < j(\delta^+)
< \gamma^{++}$.
Since $M[H_{\a_0} \ast    H '_{\a_0} ] \models `` |R''_{a_0} | = j (\delta)
$ and $|{\cal P} (R''_{\a_0}) | = j(\delta^+)$" (this last fact follows from
GCH in $M[H_{\a_0 } \ast    H '_{\a_{0}}]),$ in $V[G_{\a_0} \ast    H
'_{\a_0}]$, we can let $\langle D_\a: \a < \gamma^+ \rangle $ be an
enumeration of the dense open subsets of $R''_{\a_0} $ present in $M[H_{\a_0}
\ast    H '_{\a_0}]$. The $\prec\gamma^+$-strategic closure of $R''_{\a_0}$
in both $M[H_{\a_0} \ast    H '_{\a_0}]$ and $V[G_{\a_0} \ast    H '_{\a_0}]$
 now allows us to meet all of these dense subsets as follows.
Work in $V[G_{\a_0} \ast    H '_{\a_0}]$.
Player I picks $p_\a \in D_\a $ extending $\sup ( \langle q_\b: \b < \a
\rangle )$ (initially, $q_{-1}$ is the trivial condition), and player II
responds by picking
$q_\a \ge p_\a$ (so $q_\a \in D_\a$).
By the $\prec\gamma^+$-strategic closure of
$R''_{\a_0}$ in $V[G{}_{\a_0} \ast    H
{}'_{\a_0}]$,
player II has a winning strategy for this game,
so $\langle q_\a: \a < \gamma^+ \rangle$ can be taken as
an increasing sequence of conditions with $q_\a \in D_\a$ for
$\a <
 \gamma^+$.
Clearly, $  H  ''_{\a_0} = \{ p \in R''_{\a_0} : \exists \a < \gamma^+ [q_\a
\ge p] \}$ is our $M [H_{\a_0} \ast    H '_{\a_0}]$-generic object over
$R''_{\a_0}$ which has been constructed in $V[G_{\a_0} \ast H '_{\a_0}]
\subseteq V[G_{\a_0} \ast G '_{\a_0}]$, so $ H ''_{\a_0} \in V[G_{\a_0} \ast G
'_{\a_0}]$.
 
To construct $H    '''_{\a_0}$, we note first that as in our remarks in Lemma
8, since $\gamma $ must be below the least inaccessible $> \delta$ and
$\gamma$ is regular, $\gamma = \sigma^+$ for some $\sigma$.
This allows us to write in $V[G_{\a_0}]$ $Q'_{\a_0} = Q''_{\a_0} \ast \dot
 Q'''_{\a_0},$
 where $Q''_{a_0}$ is the full support iteration of $\langle P^0_{\omega, \l}
\ast (P^1_{\omega, \l} [ \dot S_\l] \times P^2_{\omega, \l} [\dot S_\l ] :
\delta^+ \le \l \le \sigma $ and $\l$ is regular$\rangle$ 
 and $\dot Q'''_{\a_0}$ is a term for $P^0_{\omega, \gamma} \ast (P^1_{\omega,
\gamma} [\dot S_\gamma]
\times P^2_{\omega, \gamma} [\dot S_\gamma])$. This factorization of
$Q'_{\a_0}$ induces through $j$ in $M[H_{\a_0} \ast H '_{\a_0} \ast    H
''_{\a_0}]$ a factorization of $R'''_{\a_0}$ into $R^4_{a_0} \ast \dot
R^5_{\a_0} = \langle$ the full support iteration of $\langle P^0_{\omega, \l}
\ast (P^1_{\omega, \l} [\dot S_\l] \times P^2_{\omega, \l} [\dot S_\l]) :
j(\delta^+) \le \l \le j (\sigma)$ and $\l$ is regular $\rangle \ast \langle
\dot P^0_{\omega, j(\gamma)} \ast (P^1_{\omega, j(\gamma)}
[\dot S_{j(\gamma)} ] \times P^2_{\omega, j(\gamma)} [\dot S_{j(\gamma)
}])\rangle$.

Work now in $V[G{}_{\a_0} \ast H'_{\a_0}].$
In $M[H_{\a_0} \ast H'_{\a_0}]$, as previously noted,
$R''_{\a_0}$ is $\gamma$-strategically closed.
Since $M[H_{\a_0} \ast H'_{\a_0}]$ has already been observed to be closed
under $\gamma$ sequences with respect to $V[G{}_{\a_0} \ast H'_{a_0}]$,
and since any $\gamma$ sequence of elements of $M[H_{\a_0} \ast
H'_{\a_0} \ast H''_{\a_0}]$ can be represented, in
$M[H_{\a_0} \ast H'_{\a_0}]$, by a term which is actually a
function $f: \gamma \to M[H_{\a_0} \ast H'_{\a_0}]$,
$M[H_{\a_0} \ast H'_{\a_0} \ast H''_{\a_0}]$
is closed under $\gamma$ sequences with respect to
$V[G{}_{\a_0} \ast H'_{\a_0}]$, i.e., if $f: \gamma \to M[H_{\a_0} \ast
H'_{\a_0} \ast H''_{\a_0}]$, $f \in V[G{}_{\a_0} \ast H'_{\a_0}]$, then $f \in
M [H_{\a_0} \ast H'_{\a_0} \ast H''_{\a_0}]$.
 
Factor (in $V[G{}\ao \ast G '\ao])$ $ G '\ao$ as $G''\ao \ast G'''\ao$, with 
$G ''\ao =\langle G{}^0_{\omega,\l} \ast ( G ^1_{\omega, \l} \times G
^2_{\omega, \l}): \delta^+ \le \l \le \sigma$ and $\l$ is regular$\rangle$ and
$ G '''\ao = G ^0_{\omega, \gamma} \ast ( G ^1_{\omega, \gamma}
\times G{}^2_{\omega, \gamma})$, where $G ''\ao$ is the projection of $G
'\ao$ onto $Q''\ao$ and $G '''\ao$ is the projection of $G '\ao$ onto
$Q'''\ao$. By our definitions, $Q''\ao \in V[G{}\ao]$ and $G ''\ao \in
V[G{}\ao \ast H'\ao]$.  Also, our construction to this point guarantees that
in $V[G{}\ao \ast H'\ao]$, the embedding $j$ extends to $j^*$: $V[G{}\ao] \to
M[H\ao \ast H'\ao \ast H''\ao ].$ Thus, as GCH in $V[G{}\ao \ast H'\ao]$
implies $V[G{}\ao \ast H'\ao] \models ``|Q''\ao | = |G{}''\ao | = \gamma$",
the last paragraph implies $\{ j^*(p) : p \in G ''\ao \} \in M[H\ao \ast H'\ao
\ast H''\ao]$.  Since $\{ j^\ast (p) : p \in G ''\ao \} \subseteq R^4\ao$, $
M[H\ao \ast H'\ao \ast H''\ao ] \models ``R^4\ao$ is equivalent to a $j^\ast
(\del) = j(\del)$-directed closed partial ordering", and $j (\del) > \gamma$,
$ q = \hbox{\rm sup} \{ j^\ast (p): p \in G ''_{\a_0} \}$ can be taken as a
condition in $R^4_{\a_0}$.
 
Note that GCH in $M[H_{\a_0} \ast H'_{\a_0} \ast H''\ao ]$ implies $M[H_{\a_0}
\ast H'_{\a_0} \ast H''_{\a_0}] \models ``|R^4_{\a_0} | = j(\gamma)"$, and by
choice of $j : V \to M$, $V[G{}_{\a_0} \ast H'_{\a_0}] \models ``|j(\gamma) |
= \gamma^+$ and $|j(\gamma^+)| = \gamma^+$''.  Hence, as the number of dense
open subsets of $R^4_{\a_0}$ in $M[H_{\a_0} \ast H'_{\a_0} \ast H''_{\a_0}]$
is $(2^{ j(\gamma)})^{M[H_{a_0} \ast H'_{a_0} \ast H''_{\a_0}]} = (j
(\gamma)^+)^{M[H_{\a_0} \ast H'_{\a_0} \ast H''_{\a_0}]}$ which has
cardinality $(\gamma^+)^V = (\gamma^+)^{V [G{}\ao \ast H'\ao]}$, we can let
$\langle D_\a : \a < \gamma^+ \rangle \in V[G_{\a_0 } \ast H'_{\a_0}]$
enumerate all dense open subsets of $R^4_{\a_0} $ in $M[H_{\a_0} \ast
H'_{\a_0} \ast H''_{\a_0}]$.  The $\gamma$-closure of $R^4_{\a_0}$ in
$M[H_{\a_0} \ast H'_{\a_0} \ast H''_{\a_0}]$ and hence in $V[G{}_{\a_0} \ast
H'_{\a_0}]$ now allows an $M[H_{a_0} \ast H'_{\a_0} \ast H''_{\a_0}]$-generic
object $H^4_{\a_0}$ over $R^4_{\a_0} $ containing $q$ to be constructed in the
standard way in $V[G{}_{\a_0} \ast H'_{\a_0}]$, namely let $q_0 \in D_0$ be so
that $q_0 \ge q$, and at stage $\a < \gamma^+$, by the $\gamma$-closure of
$R^4_{\a_0}$ in $V[G{}_{\a_0} \ast H'_{\a_0}]$, let $q_\a \in D_\a$ be so that
$q_\a \ge \sup ( \langle q_\b: \b < \a \rangle )$.  As before, $H^4\ao = \{ p
\in R^4\ao:\exists \a < \gamma^+ [q_\a \ge p ] \} \in V[G{}\ao \ast H'\ao]
\subseteq V[G{}\ao \ast G{}'\ao] $ is clearly our desired generic object.
 
By the above construction, in $V[G{}\ao \ast G '\ao]$, the embedding $j^* :
V[G{}\ao] \to M[H\ao \ast H'\ao \ast H''\ao]$ extends to an embedding $j^{**}:
V[G{}\ao \ast G ''\ao] \to M[H\ao \ast H'\ao \ast H''\ao \ast H^4\ao]$.  We
will be done once we have constructed in $V[G{}\ao \ast G '\ao]$ the
appropriate generic object for $R^5\ao = P^0_{\omega, j(\gamma)} *
(P^1_{\omega, j (\gamma)} [ \dot S_{j(\gamma)}] \times P^2_{\omega, j(\gamma)}
[\dot S_{j(\gamma)} ]) = (P^0_{\omega, j(\gamma)} \ast P^2_{\omega, j(\gamma)}
[\dot S_{j(\gamma)}]) \ast P^1_{\omega, j (\gamma)} [\dot S_{j (\gamma)}]$.
To do this, first rewrite $G '''\ao$ as $(G{}^0_{\omega, \gamma} \ast G
^2_{\omega, \gamma}) \ast G ^1_{\omega , \gamma}$.  By the nature of the
forcings, $G ^0_{\omega, \gamma} \ast G ^2_{\omega, \gamma}$ is $V[G_{\a_0}
\ast G''_{\a_0}]$-generic over a partial ordering which is $(\g,
\infty)$-distributive. Thus, by a general fact about transference of generics
via elementary embeddings (see [C], Section 1.2, Fact 2, pp$.$ 5-6), since
$j^{**} : V[G_{\a_0} \ast G''_{\a_0}] \to M[H_{\a_0} \ast H'_{\a_0} \ast
H''_{\a_0} \ast H^4_{\a_0}] $ is so that every element of $ M[H_{\a_0} \ast
H'_{\a_0} \ast H''_{\a_0} \ast H^4_{\a_0}]$ can be written $j^{**}(F)(a)$ with
$\dom(F)$ having cardinality $\g$, ${j^{**}} '' G^0_{\omega, \g} \ast
G^2_{\omega, \g}$ generates an $ M[H_{\a_0} \ast H'_{\a_0} \ast H''_{\a_0}
\ast H^4_{\a_0}]$-generic set $H^5_{\a_0}$.
 
It remains to construct $H^6\ao$, our $ M [H\ao \ast H'\ao \ast H''\ao
\ast H^4\ao \ast H^5 \ao]$-generic object over $P^1_{\omega, j(\gamma)}
 [S_{j(\gamma)}]$.
To do this, first note that $H^4\ao$ (which was constructed in $V[G\ao \ast
 H'\ao] )$ is $M[H\ao \ast H'\ao \ast H''\ao]$-generic over
$R^4\ao$, a partial ordering which in $M[H\ao \ast H'\ao \ast H''\ao]$ is
$j(\delta)$-closed. Since $j(\delta) > \gamma$ and $
M[H\ao \ast H'\ao \ast H''\ao]$ is closed under $\gamma$ sequences with respect
 to $V[G\ao \ast H'\ao]$, we can apply earlier reasoning to infer $M[H\ao \ast
H'\ao \ast H''\ao \ast H^4\ao]$ is closed under $\gamma$ sequences with
respect to $V [G\ao \ast H'\ao ]$, i.e., if
$f: \gamma \to M[H\ao \ast H'\ao \ast H''\ao \ast H^4\ao]$,
$ f \in V[G\ao \ast H'\ao ]$ then $f \in M[H\ao \ast H'\ao \ast H''\ao \ast
 H^4 \ao]$.
 
Choose in $V[G\ao \ast G'\ao]$ an enumeration $\langle p_\a: \a < \gamma^+
\rangle$ of $G^1_{\omega, \gamma}$.  Working now in $V[G\ao \ast G'\ao]$, let
$f$ be an isomorphism between (a dense subset of) $P^1_{\omega, \gamma}
[S_\gamma]$ and $Q^1_\gamma$.  This gives us a sequence $\langle f(p_\a): \a <
\gamma^+ \rangle$ of $\gamma^+$ many compatible elements of $Q^1_\gamma$.
Letting $p'_\a = f(p_\a)$, we may hence assume that $I = \la p'_\a : \a < \g^+
\ra$ is an appropriately generic object for $Q^1_\g$.  By Lemma 6, $V[G_{\a_0}
\ast G''_{\a_0} \ast G^0_{\omega, \g} \ast G^1_{\omega, \g} \ast G^2_{\omega,
\g}] = V[G_{\a_0} \ast G'_{\a_0}]$ and $V[G_{\a_0} \ast G''_{\a_0} \ast
G^0_{\omega, \g} \ast G^1_{\omega, \g}] = V[G_{\a_0} \ast H'_{\a_0}]$ have the
same $\g$ sequences of elements of $V[G_{\a_0} \ast G''_{\a_0}]$ and hence of
$V[G_{\a_0} \ast H'_{\a_0}]$. Thus, any $\gamma$ sequence of elements of
$M[H\ao \ast H'\ao \ast H''\ao \ast H^4\ao ]$ present in $V[G\ao \ast G'\ao]$
is actually an element of $V[G\ao \ast H'\ao]$ (so $M[H\ao \ast H'\ao \ast H''
\ao \ast H^4\ao]$ is really closed under $\gamma $ sequences with respect to
$V[G\ao \ast G'\ao])$.
 
For $\a \in (\g, \g^+)$ and $p \in Q^1_\g$, let
$p \vert \a = \{\la \la \rho, \sigma \ra, \eta \ra \in
Q^1_\g : \sigma < \a \}$ and $I \vert \a = \{ p \vert \a :
p \in I \}$. It is clear $V[G\ao \ast G'\ao] \models
``|I \vert \a| = \g$ for all $\a \in (\g, \g^+)$''.
Thus, since $Q^1_{j(\g)} \in
M[H\ao \ast H'\ao \ast H''\ao \ast H^4\ao]$ and $
M[H\ao \ast H'\ao \ast H''\ao \ast H^4\ao] \models
``Q^1_{j(\g)}$ is $j(\g)$-directed closed'', the facts $
M[H\ao \ast H'\ao \ast H''\ao \ast H^4\ao]
$ is closed under $\g$ sequences with respect to
$V[G\ao \ast G'\ao]$ and $I$ is compatible imply that
$q_\a = \cup \{j^{**}(p) : p \in I \vert \a \}$ for
$\a \in (\g, \g^+)$ is well-defined and is an element of
$Q^1_{j(\g)}$. Further, if
$\la \rho, \sigma \ra \in \dom(q_\a) -
\dom(\underset \b < \a \to{\cup} q_\b)$
($\underset \b < \a \to{\cup} q_\b \in
Q^1_{j(\g)}$ as $
M[H\ao \ast H'\ao \ast H''\ao \ast H^4\ao]
$ is closed under $\g$ sequences with respect to
$V[G\ao \ast G'\ao]$), then $\sigma \in
[\underset \b < \a \to{\cup} j(\b), j(\a))$.
(If $\sigma < \underset \b < \a \to{\cup} j(\b)$, then
let $\b$ be minimal so that $\sigma < j(\b)$, and let
$\rho$ and $\sigma$ be so that $\la \rho, \sigma \ra
\in \dom(q_\a)$. It must thus be the case that for some
$p \in I \vert \a$, $\la \rho, \sigma \ra \in
\dom(j^{**}(p))$. Since by elementarity and the definitions of
$I \vert \b$ and $I \vert \a$, for $p \vert \b = q
\in I \vert \b$, $j^{**}(q) = j^{**}(p) \vert j(\b) =
j^{**}(p \vert \b)$, it must be the case that
$\la \rho, \sigma \ra \in \dom(j^{**}(q))$. This means
$\la \rho, \sigma \ra \in \dom(q_\b)$, a contradiction.)
 
We define now an $ M [H\ao \ast H'\ao \ast H''\ao \ast H^4\ao
 \ast
 H^5\ao]$-generic object $H^{6,0}\ao $ over $Q^1_{j(\gamma)}$ so that $p
 \in f''
G^1_{\omega, \gamma}$ implies $j^{**} (p) \in H^{6,0}\ao $.
First, for $\b \in (j(\gamma), j(\gamma^+))$, let
$Q^{1,\b}_{j(\gamma)} \in
M[H_{\a_0} \ast H'_{\a_0} \ast H''_{\a_0} \ast H^4_{\a_0}]$ be
the forcing for adding $\b$ many Cohen subsets to $j(\gamma)$,
i.e., $Q^{1,\b}_{j(\gamma)} = \{g: j(\gamma) \times \b \to
\{0,1\} : g$ is a function so that $|$dom$(g)|<
j(\gamma)\}$, ordered by inclusion.
Next, note that since $
M[H\ao \ast H'\ao \ast H''\ao \ast H^4\ao \ast H^5\ao]
\models {\hbox{\rm GCH}}$, $
M[H\ao \ast H'\ao \ast H''\ao \ast H^4\ao \ast H^5\ao]
\models ``Q^1_{j(\g)}$ is $j(\g^+)$-c.c. and
$Q^1_{j(\g)}$ has $j(\g^+)$ many maximal antichains''.
This means that if $\A \in
M[H\ao \ast H'\ao \ast H''\ao \ast H^4\ao \ast H^5\ao]
$ is a maximal antichain of $Q^1_{j(\g)}$, then
$\A \subseteq Q^{1, \b}_{j(\g)}$ for some
$\b \in (j(\g), j(\g^+))$. Also, since
$V \subseteq V[G\ao \ast G''\ao] \subseteq
V[G\ao \ast H'\ao] \subseteq V[G\ao \ast G'\ao]$ are
all models of GCH containing the same cardinals and
cofinalities, $V[G\ao \ast G'\ao] \models
``|j(\g^+)| = \g^+$''. The preceding thus means we can let
$\la \A_\a : \a < \g^+ \ra \in
V[G\ao \ast G'\ao]$ be an enumeration of the maximal
antichains of $Q^1_{j(\g)}$ present in $
M[H\ao \ast H'\ao \ast H''\ao \ast H^4\ao \ast H^5\ao]
$.
 
Working in $V[G\ao \ast G'\ao]$, we define now an
increasing sequence $\la r_\a : \a \in
(\g, \g^+) \ra$ of elements of $Q^1_{j(\g)}$ so that
$\forall \a < \g^+[r_\a \ge q_\a$ and
$r_\a \in Q^{1, j(\a)}_{j(\g)}]$ and so that
$\forall \A \in \la \A_\a : \a \in (\g, \g^+) \ra
\exists \b \in (\g, \g^+)
\exists r \in \A[r_\b \ge r]$.
Assuming we have such a sequence,
$H^{6,0}_{\a_0} = \{p \in Q^1_{j(\g)} :
\exists r \in \la r_\a : \a \in (\g, \g^+) \ra
[r \ge p]\}$ is our desired generic object. To define
$\la r_\a : \a \in (\g, \g^+) \ra$, if $\a$ is a limit,
we let $r_\a = \underset \b < \a \to{\cup} r_\b$.
By the facts $\la q_\b : \b \in (\g, \g^+) \ra$
is (strictly) increasing and $
M[H\ao \ast H'\ao \ast H''\ao \ast H^4\ao]
$ is closed under $\g$ sequences with respect to
$V[G\ao \ast G'\ao]$, this definition is valid.
Assuming now $r_\a$ has been defined and we wish to define
$r_{\a + 1}$, let $\la \B_\b : \b < \eta \le \g \ra$ be
the subsequence of $\la \A_\b : \b \le \a + 1 \ra$
containing each antichain $\A$ so
that $\A \subseteq Q^{1, j(\a + 1)}_{j(\g)}$.
Since $q_\a, r_\a \in Q^{1, j(\a)}_{j(\g)}$,
$q_{\a + 1} \in Q^{1, j(\a + 1)}_{j(\g)}$, and
$j(\a) < j(\a + 1)$, the condition
$r'_{\a + 1} = r_\a \cup q_{\a + 1}$ is well-defined,
as by our earlier observations, any new elements of
$\dom(q_{\a + 1})$ won't be present in either
$\dom(q_\a)$ or $\dom(r_\a)$. We can thus using the fact $
M[H\ao \ast H'\ao \ast H''\ao \ast H^4\ao]
$ is closed under $\g$ sequences with respect to
$V[G\ao \ast G'\ao]$ define by induction an
increasing sequence $\la s_\b : \b < \eta \ra$ so that
$s_0 \ge r'_{\a + 1}$, $s_\rho =
\underset \b < \rho \to{\cup} s_\b$ if $\rho$ is a limit,
and $s_{\b + 1} \ge s_\b$ is so that
$s_{\b + 1}$ extends some element of $\B_\b$. The
just mentioned closure fact implies
$r_{\a + 1} = \underset \b < \eta \to{\cup} s_\b$ is
a well-defined condition.
 
In order to show $H^{6,0}_{\a_0}$ is $
M[H\ao \ast H'\ao \ast H''\ao \ast H^4\ao \ast H^5\ao]
$-generic over $Q^1_{j(\g)}$, we must show that
$\forall \A \in \la \A_\a : \a \in (\g, \g^+) \ra
\exists \b \in (\g, \g^+)
\exists r \in \A[r_\b \ge r]$.
To do this, we first note that
$\la j(\a) : \a < \g^+ \ra$ is unbounded in
$j(\g^+)$. To see this, if $\b < j(\g^+)$ is an ordinal,
then for some $g : \g \to M$ representing $\b$, we can
assume that for $\l < \g$, $g(\l) < \g^+$. Thus, by the
regularity of $\g^+$ in $V$, $\b_0 =
\underset \l < \g \to{\cup} g(\l) < \g^+$, and
$j(\b_0) > \b$.
This means by our earlier remarks that if
$\A \in \la \A_\a : \a < \g^+ \ra$,
$\A = \A_\rho$, then we can let $\b \in (\g, \g^+)$
be so that $\A \subseteq Q^{1, j(\b)}_{j(\g)}$. By
construction, for $\eta > \max(\b, \rho)$, there is
some $r \in \A$ so that $r_\eta \ge r$. Finally, since
any $p \in Q^1_\g$ is so that for some $\a \in
(\g, \g^+)$, $p = p \vert \a$,
$H^{6,0}_{\a_0}$ is so that if $p \in
f''G^1_{\omega, \g}$, then $j^{**}(p) \in
H^{6,0}_{\a_0}$.
 
Note now that our earlier work ensures $j^{**}$
extends to $j^{***} : V [G\ao \ast G''\ao \ast G^0_{\omega, \gamma}
\ast G^2_{\omega, \gamma}] \to M [H\ao \ast H'\ao \ast H''\ao \ast
H^4\ao \ast H^5\ao ]$.  By Lemma 4, the isomorphism $f$ is definable over $ V
[G\ao \ast G''\ao \ast G^0_{\omega, \gamma} \ast G^2_{\omega, \gamma}]$.  This
means the notions $j^{***} (f) $ and $j^{***} (f^{-1})$ make sense, so
$j^{***} (f)$ is a definable isomorphism over $ M [H\ao \ast H'\ao \ast H''\ao
\ast H^4\ao \ast H^5\ao]$ between (a dense subset of) $P^1_{\omega, j(\gamma)}
[S_{j(\gamma)}]$ and $Q^1_{j(\gamma)}$, and $j^{***} (f^{-1})$ is its inverse.
If $H^6\ao = \{ j^{***} (f^{-1}) (p) : p \in H^{6,0}\ao \}$, then it is now
easy to verify that $H^6\ao$ is an $ M [H\ao \ast H'\ao \ast H''\ao \ast
H^4\ao \ast H^5\ao]$-generic object over (a dense subset of) $P^1_{\omega,
j(\gamma)} [S_{j(\gamma)}]$ so that $p \in $ (a dense subset of) $P^1_{\omega,
\gamma } [S_{ \gamma }]$ implies $j^{***} (p) \in H^6\ao$.  Therefore, for
$H''' = H^4\ao \ast H^5\ao \ast H^6\ao$ and $H = H\ao \ast H'\ao \ast H''\ao
\ast H'''\ao$, $j: V \to M$ extends to $k : V [G\ao \ast G'\ao] \to M[H]$, so
$V[G] \models ``\d$ is $\gamma$ supercompact'' if $\gamma$ is regular.  This
proves Lemma 9.
\pbf
\hfill $\square $ Lemma 9
\proclaim{Lemma 10} For $\gamma$ regular, $V[G] \models ``\delta$ is $\gamma$
strongly compact iff $\delta $ is $\gamma$ supercompact''. \endproclaim
 
\demo{Proof of Lemma 10} Assume towards a contradiction the lemma is false,
 and let $\delta < \gamma$ be so that $ V[G] \models ``\delta $ is $\gamma$
strongly compact, $\delta$ isn't $\gamma$ supercompact, $\gamma$ is regular,
and $\gamma$ is the least such cardinal".  As before, let $\delta =
\delta_\a$, i.e., $\delta$ is the $
\a$th inaccessible cardinal.
If $V \models ``\delta_\a$ is $\gamma$ supercompact", then Lemma 9 implies
$V[G] \models ``\delta_\a$ is $\gamma$ supercompact", so it must be the case
 that $V\models ``\delta_\a$ isn't $\gamma$ supercompact".
We therefore have $\lambda_\a \le \gamma$ for $\l_\a$ the least regular
cardinal so that $V \models  ``\delta_\a $ isn't $ \l_\a$ supercompact".
 
In the manner of Lemma 9, write $P = P_\a \ast \dot Q_\a \ast \dot Q'_\a$,
where $P_\a$ is the iteration through stage $\a $, $\dot Q_\a$ is a term for
the full support iteration of $\langle P^0_{\omega, \l} \ast (P^1_{\omega, \l}
[\dot S_\l] \times P^2_{\omega, \l} [\dot S_\l ] ): \delta^+ \le \l < \l_{\a}$
and $\l$ is regular$\rangle
\ast \langle \dot P^0_{\omega, \l_\a} \ast P^1_{\omega, \l_{\a}} [ \dot
 S_{\l_{\a}}]
 \rangle$, and $\dot Q'_\a$ is a term for the rest of $P$.
By our previous results, $V^{P_\a \ast \dot Q_\a} \models ``\delta_\a$ isn't
 $\l_\a$
strongly compact", and $\force_{P_\a \ast \dot Q_\a} ``\dot Q'_\a$ is
$\delta_{\a + 1}$-strategically closed'' (where $\delta_{\a + 1}$ is the least
 inaccessible
 $> \delta_\a)$.
It must thus be the case that $V^{P_\a \ast \dot Q_\a \ast \dot Q'_\a}
= V^P \models ``\delta_\a$ isn't $\l_\a $ strongly compact", so of course,
as $\l_\a \le \gamma$, $V [G] \models `` \delta_a$ isn't $\gamma$ strongly
 compact".
This proves Lemma 10.
\pbf
\hfill $\square $ Lemma 10
 
\proclaim{Lemma 11} For $\gamma$ regular, $V[G] \models ``\delta$ is
$\gamma$ supercompact'' iff $V \models ``\delta$ is $\gamma$
supercompact''. \endproclaim
 
\demo{Proof of Lemma 11} By Lemma 9, if $V \models ``\delta$ is
$\gamma$ supercompact and $\gamma$ is regular'', then $V[G]
\models ``\delta$ is $\gamma$ supercompact''. If $V[G] \models
``\delta$ is $\gamma$ supercompact and $\gamma$ is regular'' but
$V \models ``\delta$ is not $\gamma$ supercompact'', then as in
Lemma 10, for the $\a$ so that $\delta = \delta_\a$, $\lambda_\a
\le \gamma$ for $\lambda_\a$ the least regular cardinal so that
$V \models ``\delta_\a$ isn't $\lambda_\a$ supercompact''. The
proof of Lemma 10 then immediately yields that $V[G] \models
``\delta_\a$ isn't $\lambda_\a \le \gamma$ strongly compact''.
This proves Lemma 11.
\pbf\hfill $\square$ Lemma 11
 
The proof of Lemma 11 completes the proof of our Theorem in the case $\kappa$
is the unique supercompact cardinal in the universe and has no inaccessibles
 above
it. This guarantees the Theorem to hold non-trivially.
\pbf
\hfill $\square $ Theorem
 
\S 3 The General Case
 
We will now prove our Theorem under the assumption that there may be more than
one supercompact cardinal in the universe (including a proper class of
supercompact cardinals) and that the large cardinal structure above any given
supercompact can be rather complicated, including possibly many inaccessibles,
measurables, etc.  Before defining the forcing conditions, a few intuitive
remarks are in order.  We will proceed using the same general paradigm as in
the last section, namely iterating the forcings of Section 1 using Easton
supports so as to destroy those ``bad" instances of strong compactness which
can be destroyed and so as to resurrect and preserve all instances of
supercompactness.  For each inaccessible $\delta_i$, a certain coding ordinal
$\theta_i < \delta_i$ will be chosen when possible which we will use to define
$P^0_{\theta_i, \lambda}$, $P^1_{\theta_i, \lambda} [S_{\theta_i, \lambda}]$,
and $P^2_{\theta_i, \lambda} [S_{\theta_i, \lambda}]$, where $S_{\theta_i,
\lambda}$ is the non-reflecting stationary set of ordinals of cofinality
$\theta_i$ added to $\lambda^+$ by $P^0_{\theta_i, \lambda}$.  We will need to
have different values of $\theta_i$, instead of having $\theta_i = \omega $ as
in the last section, so as to destroy the $\lambda$ strong compactness of some
$\delta$ and yet preserve the $\lambda$ supercompactness of a $\delta' \neq
\delta $ when necessary.  When $\theta_i$ can't be defined, we won't
necessarily be able to destroy the $\lambda$ strong compactness of $\delta_i$,
although we will be able to preserve the $\lambda $ supercompactness of
$\delta_i$ if appropriate.  This will happen when instances of the results of
[Me] and [A] occur, i.e., when there are certain limits of supercompactness.
 
Getting specific, let $\langle \delta_i: i \in {\hbox{\rm Ord}} \rangle $
enumerate the inaccessibles of $V \models$ GCH, and let $\lambda_i > \delta_i$
be the least regular cardinal so that $V \models ``\delta_i$ isn't $\lambda_i$
supercompact" if such a $\lambda_i$ exists.  If no such $\lambda_i$ exists,
i.e., if $\delta_i$ is supercompact, then let $\lambda_i = \Omega$, where we
think of $\Omega$ as some giant ``ordinal" larger than any $\alpha \in
{\hbox{\rm Ord}}$.  If possible, choose $\theta_i < \delta_i$ as the least
regular cardinal so that $\theta_i < \delta_j < \delta_i$ implies $\lambda_j <
\delta_i$ (whenever $j < i$).  Note that $\theta_i$ is undefined for
$\delta_i$ iff $\delta_i$ is a limit of cardinals which are $< \delta_i$
supercompact because for $j < i$, if $\delta_j$ is $< \delta_i$ supercompact,
then $\l_j \ge \delta_i$.
 
We define now a class Easton support iteration $P = \langle \langle P_\a, \dot
Q_\a \rangle :  \a \in {\hbox{\rm Ord}} \rangle $ as follows: \hfil\break
\noindent 1. $P_0$ is trivial. \hfil\break
\noindent 2. Assuming $P_\a$ has been defined, $P_{\a + 1} = P_\a \ast
\dot Q_\a$, where
$ \dot Q_\a$ is a term for the trivial partial ordering unless $\a$  is regular
and for some inaccessible $\delta = \delta_i < \alpha$ with $\theta_i$ defined,
either $\delta_i$ is $\a$ supercompact or $\a = \lambda_i$.
Under these circumstances $\dot Q_\a$ is a term for $(\underset
\{ i < \a: \delta_i {\hbox{\rm is }}\;  \a \; {\hbox{\rm supercompact}} \}
 \to{\Pi}
 (P^0_{\theta_i, \a} \ast P^2_{\theta_i, \a} [\dot S_{\theta_i, \a} ]) \ast
\underset \{ i < \a:  \delta_i \; \; \hbox{\rm is} \; \;  \a \; \;
\hbox{\rm supercompact}  \} \to{\Pi}
P^1_{\theta_i, \a} [\dot S_{\theta_i, \alpha}])
\times (\underset  \{ i < \a: \a = \lambda_i \} \to{\Pi} P^0_{\theta_i, \alpha}
\ast \underset \{ i < \a: \a = \lambda_i \} \to{ \Pi}
 P^1_{\theta_i, \alpha} [\dot S_{\theta_i, \alpha}]) = (\dot P^0_\a
 \ast \dot P^1_\a) \times (\dot P^2_\a \ast \dot P^3_\a)$,
 with the proviso that elements of $\dot P^0_\a$ and $\dot P^2_\a$ will have
 full
 support, and elements of $\dot P^1_\a$ and $\dot P^3_\a$ will have support
$< \alpha$. \hfil\break\noindent
Note that unless $| \{ i < \a : \delta_i$ is $< \a$ supercompact$\} | = \a$,
 the elements of $\dot P^i_\a$ will have full support for $i = 0, 1, 2, 3.$
 
The following lemma is the natural analogue to Lemma 8.
\proclaim{Lemma 12} For $G$ a $V$-generic class over $P$, $V$ and $V[G]$ have
the same cardinals and cofinalities, and $V[G] \models $ ZFC + GCH.
\endproclaim
 
\demo{Proof of Lemma 12} We show inductively that for any $\a$, $V$ and
 $V^{P_\a}$
 have the same cardinals and cofinalities, and $V^{P_\a} \models$  GCH.
This will suffice to show $V[G] \models $ GCH  and has the same cardinals
 and cofinalities as $V$, since if $\dot R$ is a term so that $P_\a
\ast \dot R = P$, then $\force_{P_\a}$ ``The iteration $\dot R$ is
 $< \a$-strategically
 closed", meaning $V^{P_\a \ast \dot R}$ and $V^{P_\a}$ have the same cardinals
and cofinalities $\le \a$ and GCH holds in both of these models for
cardinals $<   \a$.
 
Assume now $V$ and $V^{P_\a}$ have the same cardinals and cofinalities, and
$V^{P_\a} \models  $ GCH.  We show $V$ and $V^{P_{\a + 1}} =
V^{P_\a \ast \dot Q_\a}$ have the same cardinals and cofinalities,
 and $V^{P_{\a + 1}} \models $ GCH.
If $\dot Q_\a$ is a term for the trivial  partial ordering,
this is clearly the case, so we assume $\dot Q_\a$ is not a term for the
trivial partial ordering.
Let then $\dot Q'_\a$ be a term for
$(\dot P^0_\a \ast \dot P^1_\a) \times (\underset
 \{ i < \a: \a = \l_i \} \to{\Pi} (\dot P^0_{\theta_i, \a} \ast
P^2_{\theta_i, \a} [\dot S_{\theta_i, \a}]) \ast \dot P^3_\a) =
(\dot P^0_\a \ast \dot P^1_\a) \times (\dot P^4_\a \ast \dot P^3_\a)$,
 where as earlier, the elements of $\dot P^0_\a$ and $\dot P^4_\a$ will
  have full support, and the elements of $\dot P^1_\a$ and $\dot P^3_\a$ will
 have
support $< \a$.
We are now able to rewrite $\dot Q'_\a$ as $(\underset
\{ i < \a: \delta_i \; \; \hbox{\rm is} \; \; \a \; \;
\hbox{\rm supercompact or} \; \;  \a = \l_i \} \to{\Pi}
(P^0_{\theta_i, \a} \ast P^2_{\theta_i, \a}
 [\dot S_{\theta_i, \a}])) \ast (\underset \{ i < \a: \delta_i \; \;
 \hbox{\rm is } \; \; \a \; \; \hbox{\rm supercompact or} \; \; \a = \l_i \}
\to{\Pi} P^1_{\theta_i, \a} [\dot S_{\theta_i, \a}]) =
\dot P^5_\a \ast \dot P^6_\a$, where the elements of $\dot P^5_\a
$ will have full support, and the elements of $\dot P^6_\a$ will have support
$< \a$.
By Lemma 4, in $V^{P_\a}$, each $P^0_{\theta_i, \a} \ast (P^1_{\theta_i, \a}
[\dot S_{\theta_i, \a}] \times P^2_{\theta_i, \a}
[\dot S_{\theta_i, \a}])$ is equivalent to $Q^0_\a \ast \dot Q^1_\a$.
We therefore have that in $V^{P_\a}$, $Q'_\a$ is equivalent to
$( \underset \beta < \gamma \to{\Pi} Q^0_\a) \ast (\underset \beta < \gamma
\to{\Pi} \dot Q^1_\a)$, where $\gamma = | \{ i < \a : \delta_i$ is $\a$
 supercompact or $\a = \l_i \} | $ $(\gamma$ is a cardinal in both $V$ and
$V^{P_\a}$ by induction), i.e., the full support product of $\gamma$ copies of
$Q^0_\a$ followed by the $< \a$ support product of $\gamma$ copies of
$Q^1_\a$.  Since $\gamma \le \a$, $ \underset \beta < \gamma \to{\Pi} Q^0_\a$
is isomorphic to the usual ordering for adding $\gamma $ many Cohen subsets to
$\a^+$ using conditions of support $< \a^+$, and since $\underset \beta <
\gamma \to{\Pi} Q^1_\a$ is composed of elements having support $< \a$,
$\underset \beta < \gamma \to{\Pi} Q^1_\a$ is isomorphic to a single partial
ordering for adding $\a^+$ many Cohen subsets to $\a$ using conditions of
support $<\a$.  Hence, $V^{P_\a \ast \dot Q'_\a}$ and $V^{P_\a}$ have the same
cardinals and cofinalities, and $V^{P_\a \ast \dot Q'_\a}
\models $ GCH, so $V^{P_\a \ast \dot Q'_\a}$ and $V$ have the same cardinals
and cofinalities.
And, for $G_\a$ the projection of $G$ onto $P_\a$, if $H$ is $V[G_\a]$-generic
 over $Q'_\a$, for any $i < \a$ so that $\a = \l_i$, we can omit the
portion of $H$ generic over $P^2_{\theta_i, \alpha}[S_{\theta_i, \alpha}]$ and
 thus obtain
a $V[G_\a]$-generic object $H'$ for $Q_\a$.
Since $V \subseteq V[G_\a][H'] \subseteq  V[G_\a][H]$, as in Lemma 5, it must
therefore be the case that $V,$ $V^{P_\a \ast \dot Q_\a} = V^{P_{\a + 1}}$, and
$V^{P_\a \ast \dot Q'_\a} $ all have the same
cardinals and cofinalities and satisfy GCH.
 
To complete the proof of Lemma 12, if now $\a$ is a limit ordinal, the proof
that $V$ and $V^{P_\a}$ have the same cardinals and cofinalities and $V^{P_\a}
\models $ GCH is the same as the proof given in the last paragraph of Lemma 8,
since the iteration still has enough strategic closure and can easily be seen
by GCH to be so that for any $\beta < \a$, $|P_\beta| < \a$.  And, since for
any $\a$, $\force_{P_\a} ``\dot Q_\a$ is $< \a$-strategically closed", all
functions $f : \gamma \to \beta$ for $\gamma < \a$ and $\beta$ any ordinal in
$V[G]$ are so that $f \in V^{P_\a}$.  Thus, since $P$ is an Easton support
iteration, as in Lemma 8, $V[G]$ satisfies Power Set and Replacement.  This
proves Lemma 12.
\pbf
\hfill $\square $ Lemma 12
 
We remark that if we rewrite $\dot Q_\a$ as
$(\dot P^0_\a \times \dot P^2_\a) \ast
(\dot P^1_\a \times \dot P^3_\a)$, then the ideas used in the
proof of Lemma 12 combined with an argument analogous
to the one in the remark following the proof of Lemma 5 show
$\force_{P_\a \ast (\dot P^0_\a \times \dot P^2_\a)}
``\dot P^1_\a \times \dot P^3_\a$ is $\a^+$-c.c.''
Also, by their definitions, $\force_{P_\a}
``\dot P^0_\a \times \dot P^2_\a$ is $\a$-strategically
closed''. These observations will be used in the proof
of the following lemma, which is the natural analogue to
Lemma 9.
 
\proclaim{Lemma 13} If $\delta < \gamma$ and $V \models ``\delta $ is $\gamma$
 supercompact and $\gamma$ is regular'',
 then for $G$ $V$-generic over $P$, $V[G] \models ``\delta$ is $\gamma$
 supercompact".
\endproclaim
\demo{Proof of Lemma 13} We mimic the proof of Lemma 9.
Let $j: V \to M$ be an elementary embedding witnessing
the $\gamma $ supercompactness of $\delta$ so that $M \models ``\delta$ is not
$\gamma$  supercompact", and let $\a_0$ be so that $\delta = \delta_{\a_0}.$
 
Let $P = P_\delta \ast \dot Q'_\delta \ast \dot R$, where $P_\delta$ is the
iteration through stage $\delta$, $\dot Q'_\delta$ is a term for the iteration
$\langle \langle P_\a / P_\delta, \dot Q_\a \rangle: \delta \le \a \le \gamma
\rangle$, and $\dot R $ is a term for the rest of $P$.  As before, since
$\force_{P_\delta \ast \dot Q'_\delta} ``\dot R$ is $\gamma$-strategically
closed", the regularity of $\gamma$ and GCH in $V^{P_\delta
\ast \dot Q'_\delta}$ mean it suffices to show
$V^{P_\delta \ast \dot Q'_\delta} \models ``\delta$ is $\gamma$ supercompact".
 
We will again show that $j : V \to M$ extends to
$k : V[G_\d \ast G'_\d] \to M[H]$ for some
$H \subseteq j(P)$.
In $M$, $j(P_\delta \ast \dot Q'_\delta) = P_\delta \ast \dot
 R'_\delta
 \ast \dot R''_\delta \ast \dot R'''_\delta$, where
 $\dot R'_\delta$ will be a term for the iteration (as defined in
$M^{P_\d}$)
$ \langle \langle P_\a / P_\delta    , \dot Q_\a \rangle : \delta \le \a \le
\gamma \rangle$, $ \dot R'' _\delta $ will be a term for the iteration (as
defined in $M^{P_\delta \ast \dot R'_\delta})$  $ \langle \langle P_\a /
P_{\gamma+1}, \dot Q_\a \rangle: \gamma + 1 \le \a < j (\delta)  \rangle$,
and $\dot R'''_\delta$ will be a term for the iteration (as defined in
$M^{P_\delta \ast \dot R'_\delta \ast \dot R''_{\delta}})$ $\langle \langle
 P_\a / P_{j(\delta)}, \dot Q_\a \rangle : j (\delta) \le \a \le j(\gamma)
\rangle $. 
By the facts that GCH holds in both $V$ and $M$, $M^\gamma \subseteq M$, and
M $\models ``\delta $ is $< \gamma $ supercompact but $\delta $ is not $\gamma$
supercompact", $\dot R'_\delta$ will actually be a term for the iteration
$\langle \langle P_\a / P_\delta, \dot Q_\a \rangle$: $\delta \le \a < \gamma
\rangle \ast \langle (\dot P^0_{\gamma} \ast \dot P^1_\gamma) \times (\dot
P^2_\gamma \ast \dot P^3_\gamma) \rangle$, where the term for the iteration
$\langle \langle P_\a / P_\delta, \dot Q_\a \rangle : \delta \le \a < \gamma
\rangle $ is the same as in $V$, any term of the form $( \dot P^0_{\theta_i,
\gamma  } \ast P^2_{\theta_i, \gamma}
[\dot S_{\theta_i, \gamma}]) \ast P^1_{\theta_i, \gamma} [\dot S_{\theta_i,
\gamma}]$ appearing in $\dot R'_\delta$ (more specifically, in $\dot
P^0_{\gamma} \ast \dot P^1_\gamma)$ is identical to one appearing in $\dot
Q'_\delta$, and if $\dot P^0_{\theta_i, \gamma} \ast P^1_{\theta_i, \gamma}
[\dot S_{\theta_i, \gamma}] $ appears in $\dot R'_\delta$ (more specifically,
in $\dot P^2_\gamma \ast \dot P^3_\gamma)$, then either it appears as an
identical term in $\dot Q'_\delta$, or (as is the case, e.g., when $i = \a_0$
and $\theta_i$ is defined) it appears as the term $(\dot P^0_{\theta_i,
\gamma} \ast P^2_{\theta_i, \gamma} [\dot S_{\theta_i, \gamma}]) \ast
P^1_{\theta_i, \gamma} [\dot S_{\theta_i, \gamma}]$ in $\dot Q'_\delta$.  This
allows us to define $H_\delta = G_\delta$, where $G_\delta $ is the portion of
$G$ $V$-generic over $P_\delta$, and $H'_\delta = K \ast K'$, where $K$ is the
projection of $G$ onto $\langle \langle P_\a / P_\delta, \dot Q_\a \rangle$:
$\delta \le \a < \gamma
\rangle $ and $K'$ is the projection of $G$
onto
 $(P^0_\gamma \ast \dot P^1_\gamma) \times (P^2_\gamma \ast \dot P^3_\gamma)$
 as defined in $M$.
 
To construct the next portion of the generic object $H''_\delta$, note that as
in Lemma 9, the definition of $P_\delta $ ensures $|P_\delta| = \delta $ and
$P_\delta $ is $\delta$-c.c.  Thus, as before, GCH in $V[G_\delta] $ and
$M[G_\delta], $ the definition of $\dot R'_\delta$, the fact $M^\gamma
\subseteq M$, and some applications of Lemma 6.4 of [Ba] allow us to conclude
that $M[H_\delta \ast H'_\delta] = M[G_\delta \ast H'_\delta]$ is closed under
$\gamma$ sequences with respect to $V[G_\delta \ast H'_\delta]$.  Thus, any
partial ordering which is $\prec \gamma^+$-strategically closed in $M[H_\delta
\ast H'_\delta ]$ is actually $\prec \gamma^+$-strategically closed in
$V[G_\delta \ast H'_\delta ]$.
 
Observe now that if $\langle T_\a: \a < \eta \rangle $ is so that each $T_\a$
is $\prec \rho^+$-strategically closed for some cardinal $\rho$, then
$\underset\a < \eta \to{\Pi} T_\a$ is also $\prec \rho^+$-strategically
closed, for if $\langle f_\a: \a < \eta \rangle $ is so that each $f_\a$ is a
winning strategy for player II for $T_\a$, then $\underset \a < \eta \to{\Pi}
f_\a, $ i.e., pick the $\a$th coordinate according to $f_\a$, is a winning
strategy for player II for $\underset \a < \eta \to{\Pi} T_\a$.  This
observation easily implies $\force_{P_\delta \ast \dot R'_\delta} ``\dot
R^{''}_\delta$ is $\prec \gamma^+$-strategically closed'' in either
$V[G_\delta \ast H'_\delta]$ or $M[H_\delta \ast H'_\delta ]$.  The definition
of the iteration $R''_\delta$ then allows us, as in Lemma 9, to construct in
$V[G_\delta \ast H'_\delta] \subseteq V[G_\delta \ast G'_\delta] $ an
$M[H_\delta \ast H'_\delta]$-generic object $H''_\delta $ over $R''_\delta$.
As in Lemma 9, $M[H_\delta \ast H'_\delta \ast H''_\delta]$ is closed under
$\gamma$ sequences with respect to $V[G_\delta \ast H'_\delta]$.
 
Write $\dot R'''_\delta$ as $\dot R^4_\delta \ast \dot R^5_\delta$, where
$\dot R^4_\delta$ is a term for the iteration $\langle \langle P_\a /
P_{j(\delta)},\dot Q_\a \rangle: j(\delta) \le \a < j(\gamma) \rangle $ and
$\dot R^5_\delta$ is a term for $\dot Q_{j(\gamma)}$.  Also, write in $V$
$\dot Q'_\delta = \dot Q''_\delta \ast \dot Q'''_\delta$, where $\dot
Q''_\delta$ is a term for the iteration $\langle \langle P_\a / P_\delta,
\dot Q_\a \rangle : \delta \le \a < \gamma \rangle$ and $\dot Q'''_\delta $
 is a term for $\dot Q_\gamma$, and let $G'_\delta = G''_\delta \ast
 G'''_\delta$ be the
corresponding factorization of $G'_\delta$.
For any non-trivial term $\dot Q_\a = (\dot P^0_\alpha \ast \dot P^1_\alpha)
 \times (\dot P^2_\a \ast \dot P^3_\alpha)$ appearing in $\dot R^4_\delta$,
Lemma 4 and the fact elements of $\dot P^0_\alpha$ will have full support
and elements of $\dot P^1_\alpha$ will have support $< \alpha$ imply that in
$M$, for $T = P_\delta \ast \dot R'_\delta \ast \dot R''_\delta \ast \langle
\langle P_\beta / P_{j(\delta)},  \dot Q_\beta \rangle: j (\delta) \le \beta <
\alpha \rangle$, 
$\force_T $ ``(a dense subset of) $\dot P^0_\a \ast \dot P^1_\a$ is
 $\gamma^+$-directed
 closed".
Further, if $\a \in [j(\delta), j(\gamma)]$ is so that for some $i$, $\alpha =
 \lambda_i$,
 then it must be the case that $j(\delta) < \delta_i$, for if $\delta_i \le
 j(\delta)$,
 then by a theorem of Magidor [Ma2], since $M \models ``\delta_i$ is $<
 j(\delta)$
supercompact and $j (\delta)$ is $j(\gamma) $ supercompact", $M \models
``\delta_i$ is $j(\gamma)$ supercompact", a contradiction to the fact $M
\models ``\alpha = \lambda_i < j(\gamma)"$.
Hence, by the definition of $\theta_i$, it must be the case that $j(\delta)
\le \theta_i$, i.e., since $j (\delta) > \gamma$, $ \theta_i > \gamma$.
This means $\force_T ``\dot P^0_{\theta_i, \alpha}$ and $P^1_{\theta_i, \alpha}
 [\dot S_{\theta_i, \alpha}]$ are $\gamma^+$-directed closed",
 so as elements of $\dot P^2_\alpha$ will have full support and elements of
$\dot P^3_\alpha$ will have support $< \alpha$, $\force_T ``\dot P^2_\alpha
\ast \dot P^3_\alpha$ is $\gamma^+$-directed
 closed", i.e., $\force_T $ ``(A dense subset of) $(\dot P^0_\alpha \ast \dot
P^1_\alpha) \times (\dot P^2_\alpha \ast \dot P^3_\alpha)$ is
$\gamma^+$-directed closed".  Thus, in $M$, $\force_{P_\delta \ast \dot
R'_\delta \ast \dot R''_\delta}$ ``(A dense subset of) $\dot R^4_\delta$ is
$\gamma^+$-directed closed".  Therefore, using the extension of $j$, $j^\ast:
V[G_\delta] \to M[H_\delta \ast H'_\delta \ast H''_\delta]$ which we have
produced in $V[G_\delta \ast H'_\delta]$, the fact that GCH in $M[H_\delta
\ast H'_\delta \ast H''_\delta ]$ implies $M[H_\delta \ast H'_\delta \ast
H''_\delta] \models `` |R^4_\delta | = j(\gamma)$ and $2^{j(\gamma)} =
j(\gamma^+)"$, $V[G_\delta \ast H'_\delta] \models ``|j(\gamma^+) | =
(\gamma^+)^V = \gamma^+"$, and the closure properties of $M[H_\delta \ast
H'_\delta \ast H''_\delta]$, we can produce in $V[G_\delta \ast H'_\delta]$ as
in Lemma 9 an upper bound $q$ for $\{ j^\ast (p) : p \in G''_\delta \}$ and an
$M[H_\delta \ast H'_\delta \ast H''_\delta]$-generic object $H^4_\delta$ for
$R^4_\delta$ so that $q \in H^4_\delta$.  Again, as in Lemma 9, $ M[H_\d \ast
H'_\d \ast H''_\d \ast H^4_\d] $ is closed under $\g$-sequences with respect
to $V[G_\d \ast H'_\d]$. Therefore, by the remarks after the proof of Lemma 12
and the proof of Lemma 6, $ M[H_\d \ast H'_\d \ast H''_\d \ast H^4_\d] $ is
closed under $\g$-sequences with respect to $V[G_\d \ast G'_\d]$.
 
Rewrite $\dot R^5_\delta$ as $( \underset \{ i < j(\gamma): \delta_i \; \;
\hbox{\rm is} \; \; j (\gamma) \; \; \hbox{\rm supercompact} \} \to{\Pi}$
\  $(\dot P^0_{\theta_i, j(\gamma)} \ast P^2_{\theta_i, j(\gamma)}
[ \dot S_{\theta_i, j(\gamma)}]) \break\times
\underset \{ i < j (\gamma): j(\gamma) = \lambda_i \} \to{\Pi}$ \
 $\dot P^0_{\theta_i, j(\gamma)}) \ast 
      ( \underset \{ i < j(\gamma): \delta_i \; \; \hbox{\rm is} \; \;
j(\gamma) \; \; \hbox{\rm supercompact or}\; \; j(\gamma) = \lambda_i \}
 \to{\Pi} $
\ $ \dot P^1_{\theta_i, j(\gamma)} [\dot S_{\theta_i, j (\gamma)}])
\break
 = \dot R^6_\delta \ast \dot R^7_\delta$, where all elements of $\dot
 R^6_\delta$ will have full support, and all elements of $\dot R^7_\delta$
will have support $< j(\gamma)$.
By our earlier observation that products of
(appropriately) strategically closed partial orderings
retain the same amount of strategic closure, it is
clearly the case that $Q^*_\g$, the portion of
$Q_\g$ corresponding to $R^6_\d$, i.e., $Q^*_\g =
\underset \{ i < \gamma :\delta_i \; \; \hbox{\rm is} \; \; \gamma
 \; \; \hbox{\rm supercompact} \} \to{\Pi}
$ \ $ (P^0_{\theta_i, \gamma}  \ast
P^2_{\theta_i, \gamma} [\dot S_{\theta_i, \gamma}]) \times \underset
 \{ i < \gamma: \gamma = \lambda_i \} \to{\Pi} P^0_{\theta_i, \g}$,
is $\g$-strategically closed and therefore is
$(\g, \infty)$-distributive.
Hence, as we again
have that in $V[G_\delta \ast H'_\delta ]$, $j^\ast $ extends to
$j^{**}: V[G_\delta \ast G''_\delta] \to M[H_\delta \ast H'_\delta
 \ast H''_\delta \ast H^4_\delta ]$, we can use $j^{**}
$ as in the proof of Lemma 9 to transfer
$G^4_\d$, the projection of $G'''_\d$ onto $Q^*_\g$,
via the general transference principle of [C],
Section 1.2, Fact 2, pp$.$ 5-6 to an
$M[H_\d \ast H'_\d \ast H''_\d \ast H^4_\d]$-generic
object $H^5_\d$ over $R^6_\d$.
 
By its construction, since $p \in G^4_\delta$ implies $j^{**} (p) \in H^{5
}_\delta$, $ j^{**}$ extends in $V[G_\delta \ast G'_\delta]$ to $j^{***} : V
[G_\delta \ast G''_\delta \ast G^4_\delta] \to M [H_\delta \ast H'_\delta \ast
H''_\delta \ast H^4_\delta \ast H^5_\delta]$.  And, since $R^6_\d$ is
$\gamma$-strategically closed, $M [H_\delta \ast H'_\delta \ast H''_\delta
\ast H^4_\delta \ast H^5_\delta]$ and $M [H_\delta \ast H'_\delta \ast
H''_\delta \ast H^4_\delta]$ contain the same $\gamma$ sequences of elements
of $M[H_\d
\ast H'_\delta \ast H''_\delta \ast H^4_\delta]$ with respect to $V[G_\delta
 \ast
 G'_\delta ]$.
As any $\gamma$ sequence of elements of $M[H_\d \ast
 H'_\delta
 \ast H''_\delta \ast H^4_\delta \ast H^5_\delta]$
can be represented, in $M[H_\d \ast H'_\d \ast
 H''_\delta \ast
 H^4_\delta]$, by a term  which is actually a function
$f: \gamma \to M[H_\d \ast H'_\d \ast H''_\d
\ast H^4_\delta]$, and as $M[H_\d \ast H'_\d \ast
 H''_\delta \ast H^4_\delta ]$ is closed under $\gamma $ sequences with respect
 to
$V[G_\delta \ast G'_\delta]$, $M[H_\d \ast H'_\d \ast
 H''_\delta
 \ast H^4_\delta \ast H^5_\delta ]$ is closed under $\gamma$ sequences
with respect to $V[G_\delta \ast G'_\delta]$.
 
It remains to construct the $M[H_\d \ast H'_\d \ast
 H''_\delta
\ast H^4_\delta \ast H^5_\delta]$-generic object $H^6_\delta$ over
$R^7_\delta$. To do this, take
$Q^{**}_\gamma$ to be the portion of $Q_\gamma$ corresponding to
 $R^7_\delta$,
 i.e., $Q^{**}_\gamma$ is the $< \gamma$ support product
$\underset \{ i < \gamma:
\delta_i \; \; \hbox{\rm  is } \; \; \gamma \; \hbox{\rm supercompact or } \;
  \gamma
 = \l_i \} \to{\Pi}   P^1_{\theta_i, \gamma} [S_{\theta_i, \gamma}]$,
 with $G^5_\delta$ the projection of $G'''_\delta$ onto $Q^{**}_\gamma$.
Next, for the purpose of the remainder of the proof of this lemma,
if $p \in R^6_\d$ and $i < j(\g)$ is an ordinal, say that
$i \in \supp(p)$ iff for some non-trivial component $\bar p$
of $p$, $\bar p \in
P^0_{\theta_i, j(\g)}
$.
Analogously, it is clear what $i \in \supp(p)$
for $p \in R^7_\d$ means.
Now, let $A = \{i < j(\g)$ : For some
$p \in {j^{**}}''G^4_\d$, $i \in \supp(p) \}$, and let
$B = \{i < j(\g)$ : For some $q \in R^7_\d$,
$i \in \supp(q)$ but $i \not\in \supp(p)$ for any
$p \in {j^{**}}''G^4_\d \}$. Write
$A = A_0 \cup A_1$, where $A_0 = \{i \in
A : j(\g) = \l_i \}$ and
$A_1 = \{i \in
A : j(\g) \neq \l_i \}$.
Note that since $H^5_\d = \{ q \in R^6_\d :
\exists p \in {j^{**}}''G^4_\d[q \le p] \}$,
$A, A_0, A_1, B \in
M[H_\d \ast H'_\d \ast H''_\d \ast H^4_\d \ast H^5_\d]
$.
 
If $i \in A_1$, then by the genericity of $H^5_\d$,
$
P^1_{\theta_i, j(\g)}[S_{\theta_i, j(\g)}]
$ contains a dense subordering $P^*_i$ given by Lemma 4
which is isomorphic to $Q^1_{j(\g)}$. Hence, we can infer
that the ($< j(\g)$ support)
product $\underset i \in A_1 \to{\Pi} P^*_i$ is dense in the
($< j(\g)$ support)
product $\underset i \in A_1 \to{\Pi}
P^1_{\theta_i, j(\g)}[S_{\theta_i, j(\g)}]
$. We thus without loss of generality consider
$\underset i \in A_1 \to{\Pi} P^*_i$ instead of
$\underset i \in A_1 \to{\Pi}
P^1_{\theta_i, j(\g)}[S_{\theta_i, j(\g)}]
$. Further, if $i \in A_0$, then since
$j(\g) = \l_i$, by our earlier remarks,
$\theta_i > \g$. This means $
P^1_{\theta_i, j(\g)}[S_{\theta_i, j(\g)}]
$ is $\g^+$-directed closed.
 
As we observed in the proof of Lemma 4, for any
$i \in A$ and any $\la w^i, \a^i, \bar r^i, Z^i \ra \in
P^1_{\theta_i, j(\g)}[S_{\theta_i, j(\g)}]
$, the first three coordinates $\la w^i, \a^i,
\bar r^i \ra$ are a re-representation of an element of
$Q^1_{j(\g)}$. Since the $< j(\g)$ support product of
$j(\g)$ many copies of $Q^1_{j(\g)}$ is isomorphic to
$Q^1_{j(\g)}$, for any condition $p =
\la {\la w^i, \a^i, \bar r^i, Z^i \ra}_{i <
\ell_0 < j(\g)}, {\la w^i, \a^i, \bar r^i, Z^i \ra}_{i
< \ell_1 < j(\g)} \ra
\in
\underset i \in A_0 \to{\Pi}
P^1_{\theta_i, j(\g)}[S_{\theta_i, j(\g)}]
\times
\underset i \in A_1 \to{\Pi} P^*_i$,
we can in a unique and canonical way write $p$ as
$\la \bar p, \bar Z \ra$, where $\bar p \in Q^1_{j(\g)}$
and $\bar Z =
\la \la Z^i : i < \ell_0 < j(\g) \ra,
\la Z^i : i < \ell_1 < j(\g) \ra \ra$.
Further, this rearrangement can be taken so as to
preserve the order relation on
$\underset i \in A_0 \to{\Pi}
P^1_{\theta_i, j(\g)}[S_{\theta_i, j(\g)}]
\times
\underset i \in A_1 \to{\Pi} P^*_i$.
Therefore, since our remarks in the last paragraph imply
$\underset i \in A_0 \to{\Pi}
P^1_{\theta_i, j(\g)}[S_{\theta_i, j(\g)}]
\times
\underset i \in A_1 \to{\Pi} P^*_i$ is
$\g^+$-directed closed, the fact $
M[H_\d \ast H'_\d \ast H''_\d \ast H^4_\d \ast H^5_\d]
$ is closed under $\g$ sequences with respect to
$V[G_\d \ast G'_\d]$ means that we can in essence
ignore each sequence $\bar Z$ as above and apply
the arguments used in Lemma 9 to construct the generic
object for $Q^1_{j(\g)}$ to construct an $
M[H_\d \ast H'_\d \ast H''_\d \ast H^4_\d \ast H^5_\d]
$-generic object $H^{6,0}_\d$ for
$\underset i \in A_0 \to{\Pi}
P^1_{\theta_i, j(\g)}[S_{\theta_i, j(\g)}]
\times
\underset i \in A_1 \to{\Pi} P^*_i$. As before, since
$\underset i \in A_0 \to{\Pi}
P^1_{\theta_i, j(\g)}[S_{\theta_i, j(\g)}]
\times
\underset i \in A_1 \to{\Pi} P^*_i$ is
$\g^+$-directed closed,
$
M[H_\d \ast H'_\d \ast H''_\d \ast H^4_\d \ast H^5_\d \ast
H^{6,0}_\d]$ is closed under $\g$ sequences with respect to
$V[G_\d \ast G'_\d]$.
 
By our remarks following the proof of Lemma 12
and the ideas used in the remark following the proof
of Lemma 5,
$\underset i \in B \to{\Pi}
P^1_{\theta_i, j(\g)}[S_{\theta_i, j(\g)}]
$ is $j(\g^+)$-c.c. in
$
M[H_\d \ast H'_\d \ast H''_\d \ast H^4_\d \ast H^5_\d]
$ and $
M[H_\d \ast H'_\d \ast H''_\d \ast H^4_\d \ast H^5_\d \ast
H^{6,0}_\d]$. Since
$\underset i \in B \to{\Pi}
P^1_{\theta_i, j(\g)}[S_{\theta_i, j(\g)}]
$ is a $< j(\g)$ support product and
$
P^1_{\theta_i, j(\g)}[S_{\theta_i, j(\g)}]
$ has cardinality $j(\g^+)$ in $
M[H_\d \ast H'_\d \ast H''_\d \ast H^4_\d \ast H^5_\d \ast
H^{6,0}_\d]$ for any $i < j(\g)$,
$\underset i \in B \to{\Pi}
P^1_{\theta_i, j(\g)}[S_{\theta_i, j(\g)}]
$ has cardinality $j(\g^+)$ in $
M[H_\d \ast H'_\d \ast H''_\d \ast H^4_\d \ast H^5_\d \ast
H^{6,0}_\d]$. We can thus as in Lemma 9
let $\langle \A_\a: \a < \gamma^+ \rangle $
enumerate in $V[G_\del \ast G'_\del]$ the maximal antichains
of $\underset i \in B \to{\Pi}
P^1_{\theta_i, j(\g)}[S_{\theta_i, j(\g)}]
$ with respect to $
M[H_\d \ast H'_\d \ast H''_\d \ast H^4_\d \ast H^5_\d \ast
H^{6,0}_\d]$, and we can once more mimic the
construction in Lemma 9 of $H''_{\a_0}$ to produce in
$V[G_\d \ast G'_\d]$ an
$
M[H_\d \ast H'_\d \ast H''_\d \ast H^4_\d \ast H^5_\d \ast
H^{6,0}_\d]$-generic object $H^{6,1}_\d$ over
$\underset i \in B \to{\Pi}
P^1_{\theta_i, j(\g)}[S_{\theta_i, j(\g)}]
$.
If we now let $H^6_\del = H^{6,0}_\delta \ast
H^{6,1}_\delta$
and $H = H_\del \ast H'_\del \ast H''_\d \ast H^4_\delta
\ast H^5_\del \ast H^6_\del$,
 then our construction guarantees
$j : V \to M$
extends to
$k : V[G_\d \ast G'_\d] \to M[H]$, so $V[G]
\models ``\del $ is $\gamma$ supercompact".
This proves Lemma 13.
\pbf
\hfill $\square $ Lemma 13
 
We remark that the proof of Lemma 13 will work, regardless if $\theta\ao$
 is defined.
 
We prove now the natural analogue of Lemma 10.
\proclaim{Lemma 14} For $\gamma$ regular, $V[G] \models ``\delta$ is $\gamma$
strongly compact iff 
$\delta$ is $\gamma$ supercompact, except possibly if
for the $i$ so that $\delta = \delta_i$, $\theta_i$
 is undefined".
\endproclaim
 
\demo{Proof of Lemma 14}  As in Lemma 10, we assume towards a contradiction
the lemma is false, and let $\delta = \delta_{i_0} < \gamma$ be so that $V[G]
\models``\delta$ is $\gamma$ strongly compact, $\del$ isn't $\gam $
supercompact, $\theta_{i_0}$ is defined, $\gamma$ is regular, and $\gam$ is
the least such cardinal".  Since Lemma 13 implies that if $V \models ``\delta$
is $\gamma$ supercompact", then $V[G] \models `` \del$ is $\gamma$
supercompact", as in Lemma 10, it must be the case that $\lambda_{i_0} \le
\gamma$.
 
Write $P = P_{\lambda_{i_0}} \ast \dot Q_{\lambda_{i_0}} \ast \dot R$,
where $P_{\lambda_{i_0}} $ is the forcing through stage $\lambda_{i_0}$,
$\dot Q_{\lambda_{i_0}}$ is a term for the forcing at stage $\lambda_{i_{0}}$,
 and $\dot R$ is a term for the rest of the forcing.
In $V^{P_{\lambda_{i_{0}}}}$, since $V \models ``\del = \del_{i_{0}}$ isn't
$\lambda_{i_{0}}$ supercompact", we can write
$Q_{\lambda_{i_{0}}}$ as $T_0 \times T_1$, where $T_1$ is
$P^0_{\theta_{i_{0}}, \lambda_{i_{0}}} \ast
P^1_{\theta_{i_{0}}, \lambda_{i_{0}}}
[\dot S_{\theta_{i_{0}}, \lambda_{i_{0}}}]$, and $T_0$ is the rest of
 $Q_{\lambda_{i_{0}}}$.
Since $V^{P_{\lambda_{i_{0}}}} \models ``T_0 \times P^0_{\theta_{i_{0}},
 \lambda_{i_{0}}}$
 is $< \lambda_{i_{0}}$-strategically closed" (and hence adds no new bounded
 subsets
 of $\lambda_{i_{0}}$ when forcing over $V^{P_{\lambda_{i_{0}}}})$, the
 arguments of
 Lemma 3 apply in
$V^{P_{\lambda_{i_{0}}} \ast (\dot T_0 \times \dot P^0_{\theta_{i_{0}},
\lambda_{i_{0}}})}$ to show
  $V^{(P_{\lambda_{i_{0}}} \ast (\dot T_0 \times \dot P^0_{\theta_{i_{0}},
\lambda_{i_{0}}}))
\ast P^1_{\theta_{i_{0}}, \lambda_{i_{0}}} [\dot S_{\theta_{i_{0}},
 \lambda_{i_{0}}}]}
= V^{P_{\lambda_{i_{0}}} \ast \dot Q_{\lambda_{i_{0}}}} \models
``\delta_{i_{0}}$ isn't $\lambda_{i_{0}}$ strongly compact since
 $\lambda_{i_{0}}$
 doesn't carry a $\delta_{i_{0}}$-additive uniform ultrafilter".
 
It remains to show that $V^{P_{\lambda_{i_{0}}} \ast \dot Q_{\lambda_{i_{0}}}
\ast \dot R} = V^P \models  ``\delta_{i_{0}}$ isn't $\lambda_{i_{0}}$ strongly
compact".
If this weren't the case, then let $\dot {\cal U}$ be a term in
$V^{P_{\lambda_{i_{0}}}\ast \dot Q_{\lambda_{i_{0}}}}$ so that $\force_R
``\dot {\cal U}$ is a $\delta_{i_{0}}$-additive
 uniform ultrafilter over $\lambda_{i_{0}}"$.
Since $\force_{P_{\lambda_{i_{0}}} \ast \dot Q_{\lambda_{i_{"0}}}}
``\dot R$ is $\prec \lambda_{i_{0}}^+$-strategically closed" and
$V^{P_{\lambda_{i_{0}}}\ast \dot Q_{\lambda_{i_{0}}}}
 \models $ GCH,  if we let $\langle x_\a : \a < \l^+_{i_0} \rangle$ be in
$ V^{P_{\l_{i_{0}}} \ast \dot Q_{\l_{i_{0}}}}$
 a listing of all of the subsets
 of $\lambda_{i_{0}}$, as in the construction of
$H''\ao$ in Lemma 9, we can let $\langle r_\a : \a < \lambda^+_{i_{0}} \rangle
$ be an increasing sequence of elements of $R$ so that
$r_\a \| ``x_\a \in \dot {\cal U}"$.
If we now in
$V^{P_{\lambda_{i_{0}}} \ast \dot Q_{\lambda_{i_{0}}}}$ define ${\cal U}'$ by
$x_\a \in {\cal U}'$
 iff $  r_\a \force ``x_\a \in \dot {\cal U}"$,
 then it is routine to check ${\cal U}'$ is a $\del_{i_{0}}$-additive uniform
 ultrafilter
 over $\lambda_{i_{0}}$ in $V^{P_{\lambda_{i_{0}}} \ast \dot
 Q_{\lambda_{i_{0}}}}$,
 which contradicts that there is no such ultrafilter in $V^{P_{\lambda_{i_{0}}}
\ast \dot Q_{\lambda_{i_{0}}}}$.
Thus, $V^P \models ``\delta_{i_{0}}$ isn't $\lambda_{i_{0}}$
strongly compact", a contradiction to $V[G] \models ``\delta  $ is $\gamma$
 strongly
 compact''.
This proves Lemma 14.
\pbf
\hfill $\square $ Lemma 14
 
Note that the analogue to Lemma 11 holds if $\delta =
\delta_i$ and $\theta_i$ is defined, i.e., for $\gamma$
regular, $V[G] \models ``\delta$ is $\gamma$ supercompact''
iff $V \models ``\delta$ is $\gamma$ supercompact'' if $\delta
= \delta_i$ and $\theta_i$ is defined. The proof uses Lemmas 13
and 14 and is exactly the same as the proof of Lemma 11.
 
Lemmas 12--14 complete the proof of our Theorem in the general case.
\pbf
\hfill $\square$ Theorem
 
\S 4 Concluding Remarks
 
In conclusion, we would like to mention that it is possible to use
generalizations of the methods of this paper to answer some further questions
concerning the possible relationships amongst strongly compact, supercompact,
and measurable cardinals.  In particular, it is possible to show, using
generalizations of the methods of this paper, that the result of [Me] which
states that the least measurable cardinal $\kappa$ which is the limit of
strongly compact or supercompact cardinals is not $2^\kappa$ supercompact is
best possible.  Specifically, if $V \models ``$ZFC + GCH + $\k$ is the least
supercompact limit of supercompact cardinals + $\l > \k^+$ is a regular
cardinal which is either inaccessible or is the successor of a cardinal of
cofinality $> \k$ + $h : \k \to \k$ is a function so that for some elementary
embedding $j : V \to M$ witnessing the $< \lambda$ supercompactness of $\k$,
$j(h)(\k) = \lambda$'', then there is some generic extension $V[G] \models
``$ZFC + For every cardinal $\d < \k$ which is an inaccessible limit of
supercompact cardinals and every cardinal $\g \in [\d, h(\d))$, $2^\g = h(\d)$
+ For every cardinal $\g \in [\k, \l)$, $2^\g = \l$ + $\k$ is $< \lambda$
supercompact + $\k$ is the least measurable limit of supercompact cardinals''.
 
It is also possible to show using generalizations of the
methods of this paper that if $V \models ``$ZFC +
GCH + $\k < \l$ are such that $\k$ is $< \l$
supercompact, $\l > \k^+$ is a regular cardinal which is
either inaccessible or is the successor of a cardinal of
cofinality $> \k$ + $h : \k \to \k$ is a function so that
for some elementary embedding $j : V \to M$ witnessing the
$< \l$ supercompactness of $\k$, $j(h)(\k) = \l$'',
then there is some cardinal and cofinality preserving
generic extension $V[G] \models ``$ZFC + For every
inaccessible $\d < \k$ and every cardinal $\g \in
[\d, h(\d))$, $2^\g = h(\d)$ + For every cardinal
$\g \in [\k, \l)$, $2^\g = \l$ + $\k$ is $< \l$
supercompact + $\k$ is the least measurable cardinal''.
This generalizes a result of Woodin (see [CW]), who
showed, in response to a question posed to him by the
first author, that it was possible to start from a model for
``ZFC + GCH + $\k < \l$ are such that $\k$ is
$\l^+$ supercompact and $\l$ is regular'' and use Radin
forcing to produce a model for ``ZFC +
$2^\k = \l$ + $\k$ is $\d$ supercompact for all regular
$\d < \l$ + $\k$ is the least measurable cardinal''.
In addition, it is possible to iterate the forcing used in the
construction of the above model to show, for instance, that if
$V \models ``$ZFC + GCH + There is a proper class of cardinals
$\k$ so that $\k$ is $\k^+$ supercompact'', then there is some
cardinal and cofinality preserving generic extension
$V[G] \models ``$ZFC + $2^\k = \k^{++}$ iff $\k$ is
inaccessible + There is a proper class of measurable cardinals +
$\forall \k[\k$ is measurable iff $\k$ is $\k^+$ strongly
compact iff $\k$ is $\k^+$ supercompact$]$ + No cardinal $\k$
is $\k^{++}$ strongly compact''.
In this result, there is nothing special about $\kappa^+  $, and each $\kappa$
 can be $\lambda  $ supercompact for $\lambda   = \kappa^{++}$, $\lambda =
 \kappa^{+++}  $,
 or $\lambda  $ essentially any ``reasonable" value below $2^\k$.
The proof of these results will appear in [AS].
\hfill\vskip .25in
\noindent
Acknowledgement: The authors wish to thank Menachem Magidor for several
helpful conversations on the subject matter of this paper.
In addition, the authors would like to express their
gratitude to the referee for his thorough and
careful reading of the manuscript for this paper.
The referee's many corrections and helpful suggestions
considerably improved the presentation of the material
contained herein and have been incorporated into this
version of the paper.
\hfill\break\vfill\eject\frenchspacing
\centerline{References}\vskip .5in
 
\item{[A]} A. Apter, {\sl ``On the Least Strongly Compact Cardinal"},
\underbar{Israel J. Math.} {\bf 35}, 1980, 225--233.
\item{[AS]} A. Apter, S. Shelah, {\sl ``Menas' Result is
Best Possible''}, in preparation.
\item{[Ba]} J. Baumgartner, {\sl ``Iterated Forcing"}, in: A. Mathias, ed.,
\underbar{Surveys in Set Theory},
 Cambridge University Press, Cambridge, England, 1--59.
\item{[Bu]} J. Burgess, {\sl ``Forcing"}, in: J. Barwise, ed.,
\underbar{Handbook of Mathematical Logic}, North-Holland, Amsterdam, 1977,
403--452.
\item{[C]} J. Cummings, {\sl ``A Model in which GCH Holds
at Successors but Fails at Limits''},
\underbar{Transactions AMS} {\bf 329}, 1992, 1--39.
\item{[CW]} J. Cummings, H. Woodin, {\sl Generalised Prikry Forcings},
circulated manuscript of a forthcoming book.
\item{[J]} T. Jech, {\sl Set Theory}, Academic Press, New York, 1978.
\item{[KaM]} A. Kanamori, M. Magidor, {\sl ``The Evolution of Large
Cardinal Axioms in Set Theory"}, in: \underbar{Lecture Notes in Mathematics}
 {\bf 669}, Springer-Verlag, Berlin, 1978, 99--275.
\item{[KiM]} Y. Kimchi, M. Magidor, {\sl ``The Independence between
the Concepts of Compactness and Supercompactness''}, circulated
manuscript.
\item{[Ma1]}  M. Magidor, {\sl ``How Large is the First Strongly Compact
Cardinal?"},
\item{} \underbar{Annals Math. Logic} {\bf 10}, 1976, 33--57.
\item{[Ma2]} M. Magidor, {\sl ``On the Role of Supercompact and Extendible
 Cardinals
 in Logic"},
\underbar{Israel J. Math.} {\bf 10}, 1971, 147--157.
\item{[Ma3]} M. Magidor, {\sl ``There are Many Normal Ultrafilters
Corresponding to a Supercompact Cardinal"}, \underbar{Israel J. Math.} {\bf
9}, 1971, 186--192.
\item{[Ma4]} M. Magidor, unpublished; personal communication.
\item{[Me]} T. Menas, {\sl ``On Strong Compactness and Supercompactness"},
 \underbar{Annals Math. Logic} {\bf 7}, 1975, 327--359.
\item{[MS]} A. Mekler, S. Shelah, {\sl ``Does $\k$-Free Imply Strongly
$\k$-Free?"}, in:
\item{} \underbar{Proceedings of the Third Conference on Abelian Group Theory},
 Gordon and
Breach, Salzburg, 1987, 137--148.
\item{[SRK]} R. Solovay, W. Reinhardt, A. Kanamori, {\sl ``Strong Axioms of
 Infinity
 and Elementary Embeddings"},  \underbar{Annals Math. Logic} {\bf 13}, 1978,
 73--116.
 
\bye